\documentclass[12pt]{amsart}
\usepackage{amsmath}
 
\usepackage{amsxtra}
\usepackage{amscd}
\usepackage{amsthm}
\usepackage[T1]{fontenc}
\usepackage{amsfonts}
\usepackage{amssymb}
\usepackage{eucal}
\usepackage[all]{xypic}
\usepackage{color}
\newtheorem{thm}[subsection]{Theorem}
\newtheorem{cor}[subsection]{Corollary}
\newtheorem{lem}[subsection]{Lemma}
\newtheorem{prop}[subsection]{Proposition}

\theoremstyle{definition}

\newtheorem{defn}[subsection]{Definition}
\newtheorem{notation}[subsection]{Notation}
\newtheorem{rem}[subsection]{Remark}
\newtheorem{rems}[subsection]{Remarks}
\newtheorem{remark}[subsection]{Remark}
\newtheorem{question}[subsection]{Question}
\newtheorem{example}[subsection]{Example}

\usepackage{yhmath}
\DeclareSymbolFont{largesymbols}{OMX}{yhex}{m}{n}
\DeclareMathAccent{\widetilde}{\mathord}{largesymbols}{"65}

\newcommand{\thmref}[1]{Theorem~\ref{#1}}
\newcommand{\secref}[1]{\S~\ref{#1}}
\newcommand{\lemref}[1]{Lemma~\ref{#1}}
\newcommand{\defref}[1]{Definition~\ref{#1}}
\newcommand{\propref}[1]{Proposition~\ref{#1}}
\newcommand{\corref}[1]{Corollary~\ref{#1}}

\newcommand{\remref}[1]{Remark~\ref{#1}}
\newcommand{\exref}[1]{Example~\ref{#1}}

\newcommand{\nc}{\newcommand}
\nc{\renc}{\renewcommand}
\nc{\ssec}{\subsection}
\nc{\sssec}{\subsubsection}
\nc{\on}{\operatorname}

\nc\ol{\overline}
\nc\wt{\widetilde}
\nc\wh{\widehat}
\nc\tboxtimes{\wt{\boxtimes}}

\nc{\nint}[1]{\langle {#1} \rangle}
\renc{\d}{{\delta}}
\nc{\Aa}{{\mathbb{A}}}
 \nc{\Gg}{{\mathbb{G}}}  
\nc{\Hh}{{\mathbb{H}}}
 \nc{\Nn}{{\mathbb{N}}}
\nc{\Pp}{{\mathbb{P}}}
\nc{\vdiv}{ V_0}
\nc{\Rr}{{\mathbb{R}}}
\nc{\BV}{{\mathbb{V}}}
\nc{\BW}{{\mathbb{W}}}
\nc{\Zz}{{\mathbb{Z}}}

\nc{\Qq}{{\mathbb{Q}}}
\nc{\Ss}{{\mathbb{S}}}
\nc{\Cc}{{\mathbb{C}}}
\nc{\Ff}{{\mathbb{F}}}
\nc{\Ii}{{\mathcal{I}}}
\nc{\Mm}{{\mathcal{M}}}
\nc{\tT}{{\mathcal{T}}}
\nc{\CA}{{\mathcal{A}}}
\nc{\CB}{{\mathcal{B}}}
\def\te{\tilde{e}}
\nc{\CE}{{\mathcal{E}}}
\nc{\CF}{{\mathcal{F}}}
\nc{\CG}{{\mathcal{G}}}
\nc{\lL}{{\mathcal{L}}}
\nc{\CC}{{\mathcal{C}}}
\nc{\CM}{{\mathcal{M}}}
\nc{\mM}{{\mathcal{M}}}
\nc{\CN}{{\mathcal{N}}}
\nc{\Oog}{{\mathbb{O}}}
\nc{\Oo}{{\mathcal{O}}}
\nc{\pP}{{\mathcal{P}}}
  
\def\pos{{\mathcal P}}
\nc{\CQ}{{\mathcal{Q}}}
\nc{\CR}{{\mathcal{R}}}
\nc{\CS}{{\mathcal{S}}}
  
\nc{\CU}{{\mathcal{U}}}
\nc{\CV}{{\mathcal{V}}}
\nc{\CK}{{\mathcal{K}}}
\nc{\CW}{{\mathcal{W}}}
\nc{\CZ}{{\mathcal{Z}}}

\nc{\csM}{{\check{\mathcal A}}{}}
\nc{\oM}{{\overset{\circ}{\mathcal M}}{}}
\nc{\obM}{{\overset{\circ}{\mathbf M}}{}}
\nc{\oCA}{{\overset{\circ}{\mathcal A}}{}}
\nc{\obA}{{\overset{\circ}{\mathbf A}}{}}
\nc{\ooM}{{\overset{\circ}{M}}{}}
 \nc{\vM}{{\overset{\bullet}{\mathcal M}}{}}
\nc{\nM}{{\underset{\bullet}{\mathcal M}}{}}
\nc{\oD}{{\overset{\circ}{\mathcal D}}{}}
\nc{\obD}{{\overset{\circ}{\mathbf D}}{}}
\nc{\oA}{{\overset{\circ}{\mathbb A}}{}}
\nc{\op}{{\overset{\bullet}{\mathbf p}}{}}
 
\nc{\oU}{{\overset{\bullet}{\mathcal U}}{}}
\nc{\oZ}{{\overset{\circ}{\mathcal Z}}{}}
\nc{\ofZ}{{\overset{\circ}{\mathfrak Z}}{}}
\nc{\oF}{{\overset{\circ}{\fF}}}

\nc{\fa}{{\mathfrak{a}}}
\nc{\fb}{{\mathfrak{b}}}
\nc{\fg}{{\mathfrak{g}}}
\nc{\fgl}{{\mathfrak{gl}}}
\nc{\fh}{{\mathfrak{h}}}
\nc{\fj}{{\mathfrak{j}}}
\nc{\fm}{{\mathfrak{m}}}
\nc{\fn}{{\mathfrak{n}}}
\nc{\fu}{{\mathfrak{u}}}
\nc{\fp}{{\mathfrak{p}}}
\nc{\fr}{{\mathfrak{r}}}
\nc{\fs}{{\mathfrak{s}}}
\nc{\fsl}{{\mathfrak{sl}}}
\nc{\hsl}{{\widehat{\mathfrak{sl}}}}
\nc{\hgl}{{\widehat{\mathfrak{gl}}}}
\nc{\hg}{{\widehat{\mathfrak{g}}}}
\nc{\chg}{{\widehat{\mathfrak{g}}}{}^\vee}
\nc{\hn}{{\widehat{\mathfrak{n}}}}
\nc{\chn}{{\widehat{\mathfrak{n}}}{}^\vee}

\nc{\Xx}{{\mathbb X}}

\nc{\bb}{{\mathbf{b}}}
\nc{\bc}{{\mathbf{c}}}
\nc{\bd}{\partial}
\nc{\be}{{\mathbf{e}}}
\nc{\bj}{{\mathbf{j}}}
\nc{\bn}{{\mathbf{n}}}
\nc{\bp}{{\mathbf{p}}}
\nc{\bq}{{\mathbf{q}}}
\nc{\bF}{{\mathbf{F}}}
\nc{\bu}{{\mathbf{u}}}
\nc{\ou}{{\overline{u}}}
\nc{\bv}{{\mathbf{v}}}
\nc{\bx}{{\mathbf{x}}}
\nc{\bs}{{\mathbf{s}}}
\nc{\by}{{\mathbf{y}}}
\nc{\bw}{{\mathbf{w}}}
\nc{\bA}{{\mathbf{A}}}
\nc{\bK}{{\mathbf{K}}}
\nc{\bI}{{\mathbf{I}}}
\nc{\bB}{{\mathbf{B}}}
\nc{\bG}{{\mathbf{G}}}
\nc{\bC}{{\mathbf{C}}}
\nc{\bD}{{\mathbf{D}}}
\nc{\bP}{{\mathbf{P}}}
\nc{\bH}{{\mathbf{H}}}
\nc{\bM}{{\mathbf{M}}}
\nc{\bN}{{\mathbf{N}}}
\nc{\bV}{{\mathbf{V}}}
\nc{\uV}{{\overline{V}}}
\nc{\bU}{{\mathbf{U}}}
\nc{\bL}{{\mathbf{L}}}
\nc{\bT}{{\mathbf{T}}}
\nc{\bW}{{\mathbf{W}}}
\nc{\bX}{{\mathbf{X}}}
\nc{\bY}{{\mathbf{Y}}}
\nc{\bZ}{{\mathbf{Z}}}
\nc{\bS}{{\mathbf{S}}}

\nc{\sA}{{\mathsf{A}}}
\nc{\sB}{{\mathsf{B}}}
\nc{\sC}{{\mathsf{C}}}
\nc{\sD}{{\mathsf{D}}}
\nc{\sF}{{\mathsf{F}}}
\nc{\sG}{{\mathsf{G}}}
\nc{\sK}{{\mathsf{K}}}
\nc{\sM}{{\mathsf{M}}}
\nc{\sO}{{\mathsf{O}}}
\nc{\sQ}{{\mathsf{Q}}}
\nc{\sP}{{\mathsf{P}}}
\nc{\sZ}{{\mathsf{Z}}}
\nc{\sfp}{{\mathsf{p}}}
\nc{\sr}{{\mathsf{r}}}
\nc{\sg}{{\mathsf{g}}}
\nc{\sff}{{\mathsf{f}}}
\nc{\sfb}{{\mathsf{b}}}
\nc{\sfc}{{\mathsf{c}}}
\nc{\sd}{{\ltimes}}

\nc{\tH}{{\widetilde{H}}}
\nc{\tA}{{\widetilde{\mathbf{A}}}}
\nc{\tB}{{\widetilde{\mathcal{B}}}}
\nc{\tg}{{\widetilde{\mathfrak{g}}}}
\nc{\tG}{{\widetilde{G}}}
\nc{\TM}{{\widetilde{\mathbb{M}}}{}}
\nc{\tO}{{\widetilde{\mathsf{O}}}{}}
\nc{\tU}{\widetilde{U}}
\nc{\TZ}{{\tilde{Z}}}
\nc{\tx}{{\tilde{x}}}
\nc{\tX}{{\widetilde{X}}}
\nc{\tY}{{\widetilde{Y}}}
\nc{\tq}{{\tilde{q}}}

\nc{\tfP}{{\widetilde{\mathfrak{P}}}{}}
\nc{\tz}{{\tilde{\zeta}}}
\nc{\tmu}{{\tilde{\mu}}}

 \def\e{\epsilon}

  \nc{\vol}{{\mathop{\operatorname{\rm vol\,}}}}
   \nc{\cp}{{{\mathop{\operatorname{\rm cap}}}}}
  \nc{\disc}{{\mathop{\operatorname{\rm disc}}}}
  \nc{\Sym}{{\mathop{\operatorname{\rm Sym}}}}
   \nc{\Aut}{{\mathop{\operatorname{\rm Aut}}}}
 \nc{\Spec}{{\mathop{\operatorname{\rm Spec}}}}
  \nc{\spec}{{\mathop{\operatorname{\rm Spec}}}}
\nc{\Ker}{{\mathop{\operatorname{\rm Ker}}}}
 \nc{\dom}{{\mathop{\operatorname{\rm dom}}}}
\nc{\End}{{\mathop{\operatorname{\rm End}}}}
 \nc{\Hom}{\on{\Hom}}
 \nc{\GL}{{\mathop{\operatorname{\rm GL}}}}
 \nc{\Id}{{\mathop{\operatorname{\rm Id}}}}
 \nc{\rk}{{\mathop{\operatorname{\rm rk}}}}
 \nc{\length}{{\mathop{\operatorname{\rm length}}}}
\nc{\supp}{{\mathop{\operatorname{\rm supp}}}}
\nc{\val}{{\rm val}}
\nc{\res}{{\mathop{\operatorname{\rm res}}}}

\def\Ind#1#2#3{{#1} {\downarrow}_{#3} {#2} }

\def\tensor{{\otimes}}
\def\meet{\cap}
\def\union{\cup}

\def\si{\sigma}
\def\g{\gamma}
\def\G{\Gamma}

\def\<{\begin}
 \def\>{\end}

\def\m{\smallsetminus}

\nc{\seq}[1]{\stackrel{#1}{\sim}}

\def\inv{^{-1}}
\def\claim#1{{\noindent \bf Claim #1.\ }}

\def\beq#1{\begin{equation} \label{#1}}
\def\eeq{\end{equation}}

\def\prf{\begin{proof}}
\def\pv{\end{proof} }
 \def\eprf{\end{proof} }

\def\lbl#1{     \label{#1}  }
\def\a{\alpha}

 \renc{\b}{{\beta}}

  \def\hX{{\widehat{X}}}
    \def\hK{{\widehat{K}}}

\def\Ind#1#2{#1\setbox0=\hbox{$#1x$}\kern\wd0\hbox to 0pt{\hss$#1\mid$\hss}
\lower.9\ht0\hbox to 0pt{\hss$#1\smile$\hss}\kern\wd0}

 \def\Lam{\Lambda}          \def\lam{\lambda}

\def\Lam{\Lambda}

\def\tS{\widetilde{S}}

 \def\Om{\Omega}

 \usepackage{comment}

\nc{\FM}{}

 \def\pic{\widehat{Pic}}
 \def\psf{N^1_+}
 \def\wW{\mathcal{W}}

 \def\lam{\lambda}
 
 \def\cha{Ch}
\newcommand{\cupdot}{\mathbin{\mathaccent\cdot\cup}}
 \def\du{\cupdot}
 \title{Globally valued function fields:  existential closure}
 \author{Ita\"i ben Yaacov} 
 \author{Ehud Hrushovski}
 
\thanks{E.H.:  Research   funded in part by the European Research Council under the European Unions Seventh Framework Programme (FP7/2007- 2013)/ERC Grant Agreement No. 291111.} 
\begin{document}
\maketitle
 
The material presented in these notes forms part of a joint research project on the logic of fields with many valuations, connected by a product formula; we define such structures and name them {\em globally valued fields} (GVFs).  
  This text aims primarily at a proof that {\em the canonical GVF structure on $k(t)^{alg}$ is existentially closed}.  This can be read as saying that a variety {\em with a distinguished curve class} is a good approximation for a formula in the language of GVFs, in the same way
  that a variety is close to a formula for the theory ACF of algebraically closed fields.   
  
  The text is 
 based on classes given by the second author in Spring 2015 Jerusalem and  in Spring 2016 in Paris (FSMP).  It has a
 large expository element, especially concerning the   algebraic geometry needed for this.     This material is based on  (short initial sections of) \cite{debarre}, \cite{lazarsfeld}, 
\cite{gromov}, \cite{okounkov},    \cite{bfj} and \cite{yuan-zhang}.    However in  a number of cases, modifications of results there were needed.  
   One critical statement  is
a strong form of BDPP duality (\cite{bdpp}), where elements of  the interior of the dual cone to the effective
cone are shown to have   $n-1$'st roots in some limit of blowups; the usual version in the literature presents
them only as convex combinations of powers of ample divisors.   We work in arbitrary characteristic, whereas some of the literature restricts to characteristic zero.       Since these notes were first written (and perhaps before), many of these improvements have 
 undoubtedly appeared in the literature; we would be grateful for further references.    
 
   In these notes,  we restrict attention to the function field case.  Thus our globally valued fields should properly be called
   `purely non-archimedean globally valued fields'.  
There are nevertheless some applications to number fields, see \corref{agw}.  

Globally valued fields are defined in \S 1 via globalizing measures (up to renormalization),  and also via an explicit axiomatization in an appropriate language for real-valued logic.  
The equivalence of the two was proved in notes for the previous semester; in the function field case, it also
follows from the geometric characterization of quantifier-free types in \secref{qft-s}.

  We assume some basic algebraic geometry, including mainly blow-ups, Cartier divisors, Weil divisors;
  some of the material is reviewed in \secref{divisors}.   
 
\secref{convexity-s} collects some results related to convex subsets of $\Rr^n$.  The most important are \thmref{convexity},  
a multiplicative version of Legendre duality, and a finiteness statement for finitely generated 
groups
or graded rings, due to Khovanskii and Okounkov.     These convexity results will later be used to amplify a very modest
geometric input to yield a great deal of information.

  \secref{hodge-s} presents the Hodge index theorem.  
 
Volume is defined, via section growth, in \secref{volume-s}.

Okounkov bodies are presented in   \secref{okounkov-s}, thereby giving another definition of volume for big divisors.
They are used to prove log-concavity of volume   for big divisors.   
\secref{fujita-s} proves the Fujita approximation theorem.

\secref{bfj-s} presents the positive intersection product of \cite{bfj}.

The space of quantifier-free GVF types over a given constant field is described in \secref{qft-s}.

The existential closure of $k(x)^{alg}$ is proved in \secref{xc-s}.

In a final section, we  include a proof of the non-archimedean Yau theorem 
 of \cite{bfj2} and \cite{kt}, and more recently \cite{gubler+}.   The geometric content of   existential closure
 is quite close - but intriguingly not identical - to this statement, and we give a proof taking place in a finite dimensional setting,
 based on the ingredients of the proof of the main theorem.

 \ \,
 
 \par{\bf Notation for induced maps.}  
 Given a morphism of spaces or varieties $f: X \to U$, there are many induced   morphisms
 on associated objects.  They  will all be denoted in the same 
way, 
namely   $f_*$ in the 
 contravariant case and  $f^*$ if covariant.  In case of ambiguity, we make the domain explicit  as follows:
 $f_* | N_1(X)$ (a map from $N_1(X) $ to $N_1(U)$), or $f^* | N^1(U)$ (a map   $N^1(U) \to N^1(X)$. )
 
A similar convention applies to various groups and spaces formed out of algebraic cycles.  
For instance consider a Weil divisor represented by an irreducible subvariety $D$ of $X$.  We will be interested
in the image of $D$   in various
  groups, such as $Pic(X)$ and $N_1(X) = \Rr \tensor Pic(X) / Pic^o(X)$.   The class of $D$ in any such group
 will be denoted $[D]$, but if we write $[D]=[D'] \in N_1(X)$ we mean that the equality holds specifically for the classes in 
 $N_1(X)$.

\begin{section}{Globally valued fields}  

 A globally valued field is a field with a measure space of valuations and absolute values, satisfying the product formula.
 Here we consider only the geometric case and do not use absolute values.  Our valuations will all be valued in $\Rr$.
 
  Let $K$ be a field.   Let $\hK$ be  the space of $\Rr$-valued valuations of $K$.   
The topology on $\hK$ is induced from the Tychonoff topology on the functions from $K \m (0)$ to $\Rr$.

If $k$ is a subfield of $K$, let $\hK/k$ be the subspace of valuations trivial on $k$.  When $K=k(X)$ for some variety $X$,
this is the Berkovich space of $X$, from which the Berkovich spaces of all proper $k$-subvarieties have been removed.

$v_{triv}$ denotes the trivial valuation.

\begin{defn}
A {\em purely non-archimedean  GVF} is  a field $K$ along with a Baire (or regular Borel) measure $\mu$ on 
some Borel set of representatives for equivalence classes of nontrivial valuations;  such that
for any $f \in K \m (0)$,  
$v \mapsto v(f)$ is integrable, and
\[ \int v(f) d\mu(v) = 0 \]
\end{defn}

 In these notes, we will simply write GVF for purely non-archimedean  GVF.
 
 If $K'$ is  a subfield of $K$, we obtain a GVF structure on $K'$ by taking the pushforward measure $\mu'=r_* \mu | (\hK \m \hK/K')$,
 where $r$ is the restriction map from  measures on $\hK$, nontrivial on $K'$ to measures on $\hK'$.  
  
 If $g$ is an $L^\infty$ function on $\hK$ whose multiplicative inverse $g \inv$ is also $L^\infty$,
 we have a map $rn(g): \hK \to \hK$ multiplying each valuation $v$ by $g(v)$; 
 we let $\mu^g = g \inv rn(g)_*\mu$.   If $\mu$ concentrates on some set of representatives
 for the valuations up to equivalence, then $\mu^g$ corresponds on a different set, where each $v$ is replaced
 by the equivalent valuation $g(v) v$ and given weight $1/g(v)$ of the $\mu$-weight of $v$.  We view
 $(K,\mu)$ and $(K,\mu^g)$ as `the same' GVF structure.  (In particular they give all formulas the same value.)
 
 
 A {\em morphism} of globally valued fields is the composition of an inclusion with a renormalization.  
 

\ssec{A language for global fields}.

A real-valued language.   A {\em formula $\phi(x_1,\ldots,x_n)$, evaluated on a structure $A$}, gives a real-valued function
on $A^n$.   

 Basic formulas:  

 Let $Tr_n$ be the set of $n$-ary terms
in the language $+,\min,  \alpha \cdot x$ of divisible ordered Abelian groups (here $\alpha \cdot$ denotes scalar
multiplication by the rational number $\alpha$.)   We will refer to these as {\em $\Qq$-tropical  polynomials}.

Each $t \in Tr_n$ defines a positively 1-homogeneous continuous function $t$ on $\Rr^k$.

Let $S_n= \{x: \max_i |x_i|=1\}$.  Then $\{t|S_n: t \in Tr_n\}$   contains the constant
function $1$, is closed under addition,   scalar multiplication and $\min$, and separates points. By \lemref{stonew}
these are uniformly dense in $S_n$.   

For any 
$t \in Tr_n$, a symbol $R_{t}(x_1,\ldots,x_n)$ interpreted as:
$x \mapsto \int t(v   x_1 ,\ldots,v   x_n) $
(the domain being $\{x: \Pi_{i=1}^k x_i \neq 0\}$.)

Example:   $\int v(x)$.  The main axiom, the (logarithmic) {\em product formula}, asserts that $\int v(x)=0$
for any $x\neq 0$.

Example:  $\int  v(x)^+$ where $x^+ = \max(x,0)$.   By definition, this is the {\em height} of $x$, denoted $ht(x)$.  
 
 (In fact we will see later that the projective height $ht(x_1,\ldots,x_n) = \int \max(0,v(x_1),\ldots,v(x_n))$ suffices to generate the entire language, at least in the function field case.)
 
The analogue of quantifiers in real-valued logic are inf and sup operators.  Let $\psi$ be a continuous function on $\Rr$ with compact support.  Let $\phi(x,y)$ be a formula.  Then so is $\sup_x \psi(ht(x)) \phi(x,y) $.  Thus quantification is only
available over bounded height subsets.

\ssec{Universal axioms}  GVF.  (Function field case.)  

We have an integral domain $F$, along with a function $F_t: (F \m (0))^n\to \Rr$ written as:
\[ F_t(a_1,\ldots,a_n) = \int t(v(a_1),\ldots,v(a_n)) dv  \]
given for each $n$ and each   $t \in Tr_n$.   The $F_t$ are compatible with permutations of variables and dummy variables.

 Let $POS$ be the set of pairs $(\phi,t)$ of  formulas $\phi(x_1,\ldots,x_n)$ in the language of rings, and
 $t \in Tr_n$, such that the theory of valued fields $VF$ implies $ (\forall x)(\phi(x) \implies t(v(x_1),\ldots,v(x_n)) \geq 0 )$.
The universal axioms are:
\<{enumerate}
\item (Linearity:)  $F_{t_1+t_2} = F_{t_1} + F_{t_2}$.
\item (Local-global principle for  positivity)  
 If $(\phi,t) \in POS$ and $\phi(a_1,\ldots,a_n)$ then $  \int t(v(a_1),\ldots,v(a_n)) dv \geq 0$.
 \item (Product formula)  $F_{x}(a) =0$ where $x$ is the identity on $\Rr$.
\>{enumerate}

 \begin{rem} \rm If we wish to axiomatize GVF fields containing a given field $k$ as a constant field,
 we may use the same axioms with respect to the collection $POS_{k}$, defined in the same way
 for the theory $VF_k$ of valued fields containing $k$ as a trivially valued field.  To see the equivalence,
 a lemma of valued field theory is required; namely (writing $x=(x_1,\ldots,x_n)$, $y=(y_1,\ldots,y_m)$),
  if \[ VF_k \models   (\forall x)(\phi(x) \implies t(v(x_1),\ldots,v(x_n)) \geq 0) \] 
  then for some rational $\a>0$, some $\phi'(x,y)$, and some $a_1,\ldots,a_m \in k$  we have:  
 \[VF_k \models (\forall x)(\phi(x) \implies \phi'(x,a_1,\ldots,a_m))\] and
 \[ VF \models    (\forall x,y)(\phi'(x,y) \implies t(v(x_1),\ldots,v(x_n)) \geq - \a \min |v(y_1)|,\ldots,|v(y_m)|)\]
  \end{rem}

\end{section}

\begin{section}{Classical structures}.

(a)  $k(x)^{alg}$, $h(x)=\alpha$.   

(b)  $\Qq^{alg}[r]$ ($ h(2)=r$.)  This is not a purely non-archimedean GVF for fixed $r$; but it is so asymptotically, 
i.e.   any axiom is true for small enough $r$.  

(c)  $\Mm[l,r]$ value distribution theory.  Not a GVF at all for fixed $r$, but is a purely non
archimedean GVF for fixed $r$; but it is so asymptotically. 

We elaborate a little on (c).   
  Let $\Mm$ be a countably generated  subfield of the field of meromorphic functions.     Fix a  function
$\eta(r)$ (diverging to $\infty$; such as $\log(r)$ or $r^d$), and also an ultrafilter $u$ on $\Rr^{>0}$, avoiding finite measure sets.


Let $\mu_r$ be the measure space on $\{a: 0< |a|\leq r\}$  giving mass $\log(r/|a|) / \eta(r)$ to each point $0<|a|<r$, 
and the  uniform measure of mass $1/\eta(r)$ to the circle $|t|=r$.  Define 
\[ v_a(f)=ord_a f  \hbox{ for } |a|<r, \ \  v_t(f)=- \log |f(t)| \]
\[ ht_{\eta,u}(f) = \lim_{r \to u} \int  \max ( v_a  f  ,0) d\mu_r a\]
\[  \Mm[\eta,u] = \{f \in \Mm:   ht_{\eta,u}(f) < \infty \} \]
  \[{R_t}(f_1,\ldots,f_n) :=\lim_{r \to u}  \int t( v_a  f_1 ,\ldots,v_a f_n ) d\mu_r a\].


The product axiom is Jensen's formula:  
\[ 
  \sum_{0<|a|<r} \log \frac{r}{|a|}  ord_a(f) 
  + \frac{1}{2\pi} \int_0^{2\pi} -\log|f(re^{i\theta})| \, d\theta  =O(1)\]

The $O(1)$ term is actually $-\log |f(0)|$, which does not depend on $r$.
 
We will see that in the language described above, $\Mm[\eta]^{alg}$ has the same {\em universal} theory as the ultraproduct of the $\Qq^{alg}[r]$,
and also as $\Cc(t)^{alg}[1]$.     This
formalizes a part of  Vojta's dictionary \cite{vojta} between 
Nevanlinna theory and number theory.   As
the GVF language as we defined it does not  permit 
discussing the support of a function, the most interesting parts of the dictionary are not currently accessible.

\begin{thm} \label{ec}  $k(x)^{alg}[r] $ is existentially closed, for any $r>0$.  \end{thm}

%

\begin{rem}     It follows from \thmref{ec} that the asymptotic universal theory of $\Qq^{alg}[r]$ is precisely
GVF (+ char = 0 + h(2)=0.  ).    Also it is existentially closed over any single parameter.   \end{rem}

Let $X$ be a   projective variety over a field $k$.  
Let $Y$ be a subscheme, cut out by specified homogeneous polynomials $g_1,\ldots,g_{k_Y}$; and let $v$ be a valuation of a valued field $L \supseteq k$.  For $x \in X(L)$ we define
\[ \d_v(x,Y) = \min_{g} v(g_i(\bar{x}))  \]
where $\bar{x}$ is a representative of $x$ in homogeneous coordinates, chosen so that the minimal valuation of a coordinate is $0$. 
 
Then $\d=\d_v(x,Y)$ 
 (or rather $2^{-\d}$)
  is a notion of distance from $x$ to $Y$.   If $L$ is a number field and $x \in X(L)$, let
$\rho(x,Y)^L = \int \delta_v(x,Y) dv$ be the weighted sum of the local distances from $x$ to $Y$.    

Note that $\rho(x,Y)^L$ is the $L$-value of a certain quantifier-free formula $\phi_Y(x)$  in the language of GVF's.

Let $Y_0,\ldots,Y_m$ be subschemes of $X$.  Let $\rho_Y(a)=(\rho_{Y_0}: \cdots : \rho_{Y_m}) \in \Pp^m(\Rr)$.  

For a curve $C$ on $X$ let   $i(C,Y_m)$ be the  intersection number of $C$ with $Y_m$.
 Similarly let $\bar{i} (a) = (i(C,Y_0): \cdots : i(C,Y_m)) $.  


\begin{cor} \label{agw} Let $X \subset \Pp^n$ be a smooth projective variety over $\Qq$.  Let $Y_0,\ldots,Y_m$ be subschemes of $X$.  Let $e \in \Pp^m(\Rr)$.   
   Assume there exists a Zariski dense sequence $a_i \in V(\Qq)$ 
with $\rho_Y(a_i) \to e$ in $\Pp^m(\Rr)$.  Then there exists    a curve $C$ on $V$ with $\bar{i}(C,Y_k)  $ as close as one wishes to $e$.\end{cor}

In fact the same conclusion is true if the $a_i$ are only assumed to be of bounded degree in $V(\Qq^{alg})$; or just,  of  height approaching infinity.

One can also state, given $X,Y_0,\ldots,Y_m$, that there exists a proper  subvariety $W$ of $X$ (not necessarily irreducible), such that for $e \in \Pp^m(\Rr)$ and any sequence of elements $a_i \in V(\Qq^a)$ of unbounded
height, the same conclusion holds. 

\prf  We may assume $Y_0$ is an ample divisor, and use height with respect to $Y_0$.
Let $a_i \in V(\Qq^{alg})$ ($i=1,2,\cdots$) have height approaching $\infty$,  such that $\rho_Y(a_i) \to e$.  
Choose $r_i=ht(2)/ht(a_i)$ so that $\Qq^{alg}[r_i]$ gives $a_i$ height $1$.  Consider any   non-principal ultrafilter $u$
on the index set $\Nn$, and let 
$(L,a)$ be the GVF ultraproduct of $(\Qq^{alg}[r_i],a_i)$.  Then $(L,a)$ is a purely non-archimedean GVF,
$a \in V(L)$,
and    $\rho_Y(a) = \phi(a)^L = e$.    Let $\e>0$ and fix some metric on $\Pp^m(\Rr)$.    Then by \thmref{ec}
there exists $a' \in V(K)$, $K=k(t)^{alg}$ with $e'=\phi(a')^{K} $ satisfying $|e'-e|<\epsilon$.    In fact $a' \in V(k(C))$
for some curve $C$, so $a'$ corresponds to a morphism $g: C \to V$.   We may choose $a'$
so that $g$ avoids any given proper subvariety of $V$.  By computing the meaning
of $\phi$ in $k(t)^{alg}$ we see that $\bar{i}(C,Y_k)  =e'$.
\eprf

\begin{rem}     Conversely, if   $C$ is a curve on $X$ defined over $\Qq^{alg}$, then for any sequence of distinct $a_i \in C(\Qq^{alg})$ of bounded degree over $\Qq$,
$\rho_Y(a) \to   i_Y(C)$.    \rm This follows upon taking normalized ultraproducts as above,  from the uniqueness of the GVF structure on $k(C)$, \lemref{whaples}.
\end{rem}

By applying the above remark along with \corref{agw}, we obtain a purely arithmetic corollary:

\begin{cor}  Let $X \subset \Pp^n$ be a smooth projective variety over $\Qq$.  Let $Y_0,\ldots,Y_m$ be subschemes of $X$.  Let $e \in \Pp^m(\Rr)$.  
   Assume there exist    $a_i \in V(\Qq^{alg})$ of height approaching $\infty$,
with $\rho_Y(a_i) \to e$.  Then there exists such a sequence $a_i'$ of bounded degree over $\Qq$, 
$[\Qq[a_i']:\Qq] <d$.\end{cor}

\ssec{Fekete-Szeg\H{o}}
 
Let $C$ be a compact subset of $\Cc$.  The {\em energy} of a measure $\mu$ on $C$ is defined to be:
\[  \int  -\log |x-y| d\mu(x) d \mu(y)   \]
In particular, a measure of finite energy cannot assign positive mass to any point.  
If $C$ supports some measure of finite energy, there exists a unique probability measure 
attaining the minimum possible energy $m_C$;   it is called the {\em equilibrium measure} $\mu_C$, and one defines
the {\em capacity} of $C$ to be $\cp(C)=exp(-m_C)$.  

The classical Fekete-Szeg\H{o} theorem asserts that if the capacity is $\geq 1$, there exist infinitely many Galois orbits
of algebraic integers
contained in any given neighborhood $U$ of $C$; if it is $<1$,   some neighborhood $U$ contains  only finitely many such orbits.    Note in the first case that  the height of
the Galois orbits is bounded, so their degree must approach $\infty$.
   
This has adelic generalizations, due to David Cantor; the complex field is replaced by similar conditions at finitely many places $S$, while the algebraic integers are replaced by $S$-integers.  

Restriction to a neighborhood of a closed set is not a GVF condition; however, insisting that
for any neighborhood $U$ of $C$, 
{\em most} of the orbit lie in $U$, can be viewed as such.  Similarly we relax the notion of {\em  $S$-integrality}
to assert that a natural measure of non-integrality, even when summed over all places outside $S$, approaches
zero.    This may be considered an $L^1$-analogue of the Fekete-Szeg\H{o} question.  
 In the next section we prove such a   dichotomy  over an existentially closed GVF, 
 for varieties of any dimension.      In particular it applies to $k(t)[1]$ and its finite extensions.

\begin{lem}  Let $A_i$ be an infinite sequence of Galois orbits of algebraic integers, $|A_i| \to \infty$.  Assume 
$A_i \subset C$ or at least that $|A_i| / |A_i \meet C|   \to 1$.  Then $C$ has capacity $\cp(C) \geq 1$. 
If $\cp(C)=1$, then for any such sequence, the sets $A_i$ equi-distribute to the equilibrium measure on $C$.  
 \end{lem}
\prf  Let $a_i \in A_i$.  Take an ultrapower $K$ of $\Qq^{alg}$, viewed as a GVF.
Let $a \in K$ be the
limit of the $a_i$.   Let $K'$ be a further ultrapower, embed $K$ into $K'$ and let $a'$ be the new limit of the $a_i$.  

Let $b = a-a'$.   Note $a \neq a'$, so $b \neq 0$.  Note that for any $v$, and all algebraic integers $c$, we have
$\int   v(c)^- v(2)^+ dv  =0$, since either $v(2)<0$ or $v(c)=0$.  Thus the same holds for $a,a'$ and $b$;
so $v(b) \geq 0$ for all non-archimedean $v$, i.e. all $v$ with $v(2) > 0$.   Hence
\[ 0= \int v(b) dv \geq  \int_{v(2)>0} v(b) = \int_{\alpha}  - \log |\alpha(a-a')| \]
where $\alpha$ ranges over the complex embeddings of $\Qq(a,a')$.   Let $\mu$ be the measure corresponding
to the globalizing measure of $\Qq(a)$ (or $\Qq(a')$) over $\Qq$, restricted to the archimedean primes. 
Since $tp(a'/\Qq^a,a)$ is finitely satisfiable in $\Qq^a$, the archimedean part of the measure of $\Qq(a,a')$
is the product measure $\mu \times \mu$ (say by $L^2$ exhaustion ).  We have 
$ 0 \geq \int  - \log(x-x') d\mu_x d \mu_{x'} \geq  - \log \cp(C) $, so $\cp(C) \geq 1$.

In case $\cp(C)=1$, $\mu$ must equal the equilibrium measure $\mu_C$, by uniqueness of the latter.  
  \eprf
   
 \end{section}

\begin{section}{Fekete-Szeg\"o }
 \def\uX{\underline{X}}  \def\uD{\underline{D}}

Let $K$ be a GVF, with globalizing measure $\mu_K$   concentrating on a Borel set
$\Pi$ of representatives for the nontrivial valuations of $K$; including in the archimedean case
$-\log_2 |\si(x)|$, $\si: K \to \Cc$ a place.   For simplicity we can take $K$ to be countable, so that 
(away from $\mu_K$-measure $0$) $\si$ is a complex embedding.
Let $\widehat{K}=Val_K$.

     Let 
$X$ be a projective variety over $K$, with a very ample divisor $D$. 
  For $p \in {\Pi}$, let $X_p$ be the Berkovich space of $X$ over 
the Henselization of $K$ at $p$.  

It is really just as well to think of the affine variety $X_{aff} = X \m D$, whose polynomial ring is the graded ring
$\oplus L(X,mD)$; so  $D$ is the hyperplane at $\infty$ of the projectivization of $X_{aff}$.

We first consider a family of partial types that we call {\em $D$-adelic}.  
 Such a type is associated with
a  positively homogeneous function $s: \widehat{K} \to \Rr$, with $\int s < \infty$.  
We assume $s$ is given by $s(v) = v(t)$ for some term $t$ over $K$, or as a uniform limit of such functions.
 Classically, $K$ is a global field, $\widehat{K}$ is discrete and $s$ has finite support on $\Pi$.

Recall that $ ht_D(x) =  \int \max_{i=1}^n v^-(x_i) dv$.  
Define the partial type $P_{D,s}$ by:
 \[  ht_D(x) =  \int \min(s(v), \max_{i=1}^n v^-(x_i) )  \] 
 This partial type asserts that the $D$-height of $x$ is entirely attributable to the primes in the support of $S$,
 and moreover the local heights are bounded by $s$.  
 Any partial type extending some $P_{D,s}$ will be called $D$-adelic.  

To define $D$-adelic subsets of $X$.

We first consider compact sets $S \subset \widehat{K}$, such that   $S \subset \Pi$.  We may assume (adjusting $\Pi,\mu_K$ if necessary) that for some terms $t_S \geq 0$, we have $S=\{v: t(v)= 1 \}$.  Typically we could have $t_S = \max_{i=1}^l v(a_i)$;
for instance over $\Qq$, for $a_1=1,a_2=2,a_3=3$ we have $S= \{2,3,\infty\}$.  

   Note in general 
that $S$ has finite $\mu_K$-volume.  
 These sets $S$ generalize finite sets of primes.

A {\em $D$-adelic set}   is a compact subset of  of $\widehat{K(X)}$, lying above $S$.   Thus $A$ is the union of compact subsets $A_p \subset X_p$ for each $p \in S$.
{\em We assume that each $A_p $ is disjoint from $D$, i.e. lies in the Berkovich space of the affine variety $X_{aff}$.}

Alternative formulation:  
 we assume that $A$ is cut out by inequalities of term functions, including:  $h_D \leq m$; where
 $h_D(v) = \min_i v(x_i)$, $\{x_i\}$ being a basis for the global sections of $D$.

If we pass to an algebraic closure $K_p^{alg}$ of the valued field $K_p = (K,p)$, we can identify 
$A_p$ with a Galois-invariant compact subset of the Berkovich space  $X$ over $K_p^{alg}$.


The {\em Chebyshev constant}  \cha($A$)  is defined by:

\[ - \log \hbox{\cha} (A) = \lim_{n \to \infty}  \frac{1}{n} \sup_{f \in L(nD)} \int_{p \in \Pi} \inf_{v \in A_p} v(f)  \]
 
 Where $L(nD)=L_K(X,nD)$ is the $K$-space of rational functions in $K(X)$ whose poles are at most $nD$.
 
 Since $L(mD)L(nD) \subset L((m+n)D)$, 
  the terms $a_n= \sup_{f \in L(nD)} \int_{p \in \Pi} \inf_{v \in A_p} v(f) $ satisfy $a_{n+m} \geq a_n+a_m$.
  By Fekete's lemma, $\lim a_n/n$ exists and equals $\sup_n a_n/n$.      Thus
  \cha $(A)  \geq 1$ iff each $a_n \leq 0$.

\begin{example}  Let $K=\Qq$, $X=\Pp^1$, $D=\{\infty\}$.   Let $A_\infty = A$ be a compact subset of $\Rr$, and $A_p= \Oo_p$ otherwise.  
 Note that for $e \in K$ we have $ \int_{p \in \Pi} \inf_{v \in A_p} v(fe) 
 = \int_p \inf_{v \in A_p} v(f)  $.  This uses the fact that $v(e)=p(e)$ for $v \in A_p$, and that $\int_p p(e)=0$.
 Thus in the supremum it sufices to consider monic polynomials.  In this case $\inf_{v \in A_p} v(f) = 0$
 for finite $p$; so 
 \[ - \log \hbox{\cha} (A) = \lim_{n \to \infty}  \frac{1}{n} \sup_{f \in L(nD) ^{\hbox{\\ monic}}}   \inf_{v \in A_\infty} v(f)  \]
 which is the usual definition of the Chebyshev constant, see e.g. \cite{chinburg-capacity-varieties}.    It was shown by Fekete-Szeg\"o that this constant 
 equals the capacity of $A$.
\end{example}

 Let $X,Y$ be projective varieties over $K$, with very ample divisors $D_X,D_Y$; let $D = \pi_X ^* D_X + \pi_Y^* D_Y$.  Given adelic sets $A,B$ on $X,Y$ respectively, we have the obvious definition of $A \times B$ on $X \times Y$; 
 namely the pullback to the Berkovich space of $X \times Y$ of $A \times B$.
It is an adelic set.  We have $- \log \cha(A \times B) \geq -\log( \cha(A) ) + -\log (\cha(B))$, or $\cha(A \times B) \leq \cha(A) \cha(B)$.  In particular, $\cha(A^n) \leq \cha(A)^n$.
 
 \ssec{The partial type associated with an adelic set}
Recall that  a Zariski dense type in $X$ over $K$ can be   identified with a measure on valuations of $K(X)$, satisfying certain
properties.  The set of measures concentrating above each $p$ on $A_p$ forms a closed set, and thus describes a partial type; but one of unbounded height.    We can further specify that the Weil height of a realization $(x_1,\ldots,x_n)$ (in affine coordinates), 
$  \int \max_{i=1}^n v^-(x_i)$  
 be attained  adelically, indeed above $S$:   
  \[  ht_D(x) =  \int \max_{i=1}^n v^-(x_i) = \int_{p \in S}    \int_{v | p}  \max_{i=1}^n v^-(x_i)\]
With the above assumptions on $S$ and $A$, this is a closed condition and thus determines a partial type $P_{A,D}$. 
 
We call it the {\em adelic type} associated with $A$.   
  
Model theoretic questions arise immediately:  consistency; completeness; extensions to $K^{a}$; stable emeddedness and minimality.   We will see that these relate to classical definitions and theorems: 
consistency to an $L^1$-version of Fekete-Szeg\H{o} conditions, adelic completeness to  low height generic equidistribution.  Full qf-completeness will follow, in case $NS(X)$ is one-dimensional.

 Before giving a consistency criterion for $P_A$ (\lemref{fs6}), we consider adelic measures on $A$.  Let $\mu$ be
 a measure on the Berkovich space of $K(X)$ concentrating on valuations {\em nontrivial} on $K$, and
 let $(\mu_p)_{p \in P}$ be the decomposition into probability measures on the $X_p$.    Let $P_{\mu}$ be the set of adelic qf types 
 whose restriction above $p \in S$ is precisely $\mu_p$.  Again by adelic (for $(X,D)$) we mean
 that 
   \[  ht_D(x) =    \int \max_{i=1}^n v^-(x_i) d \mu(v)\]
   (where $x$ is again a basis for the global sections of $L$; in particular $K(X)=K(x)$.)
 this is a closed set of types and will be
 viewed as a partial type.  (We will consider $\mu$ such that $\mu_p$ is supported on $A_p$ for $p \in S$.)
 
 \begin{lem} \label{fs5}   $P_{\mu}$ is consistent iff for any $m \in \Nn$ and $f \in L(X,mD)$,
\[ \int v(f) d\mu  (v) \leq 0 \]
 \end{lem}

\prf  Assume $P_{\mu_{S}}$ is consistent, and so realized by $(a_1,\ldots,a_n) \in X(K) \leq \Aa^n(K)$.
Let $f \in L(X,mD)$.  Then $(f)+mD \geq 0$ (as Cartier divisors over $K$) so for $v$ trivial on $K$ we have $v(f(a)) + m \max_{i=1}^n v^-(a_i) \geq 0$.
Now since $tp(a/K)$ is adelic, for almost every $v$ that is not over some $p \in \Pi$, we have $ \max_{i=1}^n v^-(a_i) =0$; so $f(a) \geq 0$.   Since $\int v(f(a)) = 0$, 
it follows that $\int v(f)   d \mu_{S} (v) \leq 0$.

\begin{rem} \label{chbshf+} \rm  We can relax the notion of $D$-adelic by allowing the divisorial valuation trivial on $K$
corresponding to $D$, to carry weight $\alpha>0$.  Thus let $P_{\mu,\alpha}$ be the set of adelic qf types 
 whose restriction above $p \in S$ is precisely $\mu_p$, and with  
 $ht_D(x) =  \alpha+  \int \max_{i=1}^n v^-(x_i) d \mu(v)$.    We obtain in the same way:  $P_{\mu,\alpha}$
  is consistent iff for any $m \in \Nn$ and $f \in L(X,mD)$,
\[ \frac{1}{m} \int v(f) d\mu  (v) \leq  \alpha  \] 
\end{rem}

The other direction follows from \lemref{globalize}, with $Z$ consisting of valuations nontrivial on $K$,
and all valuations trivial on $K$ and centered on $D$.  The measure is extended so as to be zero on the latter.  Note that the ring of elements of $K(X)$
that are non-negative at every $v \notin Z$ is precisely $\union_m L(mD)$.  
\eprf

\begin{lem}  \label{fs6}  
$P_A$ is consistent iff   \cha$(A) \geq  1$.   \end{lem}

\prf  First, assume $P_A$ is consistent.  A complete type extending $P_A$ determines measures $\mu_p$
above each $p \in \Pi$.  Let $\mu=(\mu_p)_{p \in S}$.  Then $P_\mu$ is consistent with a $D$-adelic type $P$.
Thus for any $m \in \Nn$ and $f \in L(X,mD)$,
 $ \int v(f) d\mu(v) \leq 0 $; and as $P$ is $D$-adelic,   $\int_{p \in S} \inf_{v \in A_p} v(f)   \leq 0$.  It follows that
  $\sup_{f \in L(nD)} \int_{p \in S} \inf_{v \in A_p} v(f)  \leq 0$ and so also  
\[\lim_{n \to \infty}  \frac{1}{n} \sup_{f \in L(nD)} \int_{p \in S} \inf_{v \in A_p} v(f)  \leq 0\]
as required.

Conversely, assume \cha$(A) \geq  1$.
    We seek a measure $\mu$ on  the locally compact space $\du_p A_p$, 
satisfying the conditions of \lemref{fs5}.    Let $C(A)$ be the space of continuous functions with compact support
on $\du_p A_p$.  We are looking for a linear functional $\lambda$ on $C(A)$ such that:
\begin{enumerate}
\item $\lambda(1_{A_p}) = m(p) $.  In other words $\lambda$ extends $m$ on $C(S)$.  
\item  $\lambda( c) \geq 0$ if $c$ is a positive function.
\item  For $f \in L(X,mD)$,   $\lambda(\hat{f}) \leq 0 $; where $\hat{f}(v) = v(f)$.
\end{enumerate}

Note that any element of $C(A)$ is bounded by an element of $C(S)$ (with finite support).  

By \lemref{mriesz}, it suffices to show that if $\sum a_i p_i \in C(S)$ satisfies $\sum a_i p_i \leq \hat{f}$
with $f \in L(X,mD)$, then $\sum a_i m(p_i) \leq 0$.    Assume $\sum a_i p_i \leq \hat{f}$.  Then 
$a_i \leq \inf_{v \in A_p} v(f)$.  So $\sum a_i p_i \leq \int \inf_{v \in A_p} v(f) \leq 0$.  

\eprf

\begin{rem}  \rm The proof of \lemref{fs6} essentially gives a characterization of   uniqueness of the measure $\mu$,
or of the restriction of $\mu$ to $S$.  Namely,   any continuous function on the Berkovich space (above $S$),
with $\mu$-integral $= 0$, should be uniformly approximable (above $S$) by $\val(f)$ with 
$f \in L(X,mD)$ for some $m$.  This clearly implies  that   \cha$(A) =  1$.   Over a GVF $K$ with existentially
closed symmetric extension $K^{alg}$, one can hope for the converse.   \end{rem}

\begin{cor}  \label{fs8} Assume the {\em symmetric extension} $K^{alg}$ of $K$ is existentially closed as a GVF,
and $X$ remains irreducible over $K^{alg}$.   
Then $\cha(A)\geq 1$ iff $\cha(A^n) \geq 1$.  \end{cor}
\prf  We saw that $\cha(A^n) \leq \cha(A)^n$, so one direction is clear.  Conversely assume $\cha(A)\geq 1$.
Then $P_A$ is consistent; by taking symmetric extensions, noting that 
the restriction map from the absolute Galois group of $K(X)$ to that of $K$ is surjective, we see that
$P_A$  extends to a qf type $p_A$ over $K^{alg}$.  As $K^{alg}$ is existentially closed, $p_A$ is 
finitely satisfiable in $K$; we can find $c_1,\ldots,c_n \models p_A$ within an ultrapower of $K^{alg}$,
with $tp(c_i/K^{alg}(c_1,\ldots,c_{i-1}))$ finitely satisfiable in $K$.  It follows that $P_{A^n}$ is consistent,
and so by \lemref{fs6}, we have $\cha(A^n) \geq 1$.  \eprf


\<{lem} \lbl{mriesz} [M. Riesz] Let $V$ be a $\Qq$- vector space.  Let $U$ be a subspace of $V$, $P$ a cone, and let $h: U \to \Rr$
be a $\Qq$-linear map, positive with respect to $P \meet U$, i.e. $h(x) \geq 0$ for $x \in P \meet U$.     Assume:
    for any $v \in V$, for some $w \in U$,
$w - v \in P$. 
Then $h$ extends to a positive linear map $H$ on $(V,P)$.  \>{lem}
  
  Note that if $V$ is actually an $\Rr$-space, and $P$ an $\Rr$-cone, and any element of $V$ has the form $a-b$, with $a,b \in P$
  then $H$ will be $\Rr$-linear.  It suffices to show for $\a \in \Rr$, $a \in P$
  that $H(\a a) = \a H(a)$.  This is clear since for rational $\a' \in \Qq$ with $\a' \leq \alpha$ we have $(\a - \a') a \in P$ so 
  $H( (\a - \a') a ) =H(\a a) -  \a' H(a) \geq 0$, or $H(\a a) \geq \a' H(a)$; and if $\alpha < \alpha'$, then $(\a' - \a) a \in P$,
  so $\a' H(a) \geq H(\a a)$.     
  
  \bigskip
%

 Let $F$  be a field.  There is a  Galois correspondence $W\mapsto R^W$, $R \mapsto W^R$ between integrally closed subrings $R$ of $F$, and  certain closed subsets of $\Om_F$ consisting of non-archimedean   valuations and  invariant  under positive scalar multiplication.    
 Namely, let 
 \[  R^W =     \{r: (\forall v \in W) \, v(r) \geq 0 \} \union \{0\} \]
 Conversely given an integrally closed subring subring $R$ of $F$, let 
 \[   W^R= \{v: (\forall a \in R \m (0)) \, v(a) \geq 0 \} \]
 It is clear that the correspondence is order-preserving, that $R^W$ is integrally closed,   that $W^R$ is a closed set consisting of non-archimedean valuations,
 and closed under equivalence.  Also, $R^W = R^{cl(W)}$, and $W^R = W^{R^{int}}$ where $R^{int}$ is the integral closure of $R$ in $F$.
 We need not have $W^{R^W} = W$ in general.

 We also denote $R_Z = R^{Val_F \m Z}$,   the subring of $F$ consisting of elements whose poles lie on $Z$.  
 
\<{lem}  \label{globalize}  [See better version in file:  class1.tex]  Let $F$ be a countable field.
Let   $Z \subset Val_F$ be a Borel subset including all archimedean valuations, and invariant under multiplication
by positive scalars.  Let 
$R =R_Z= \{r: v \notin Z \implies v(r) \geq 0 \} \union \{0\}$ be the subring of $F$ consisting of elements whose poles lie on $Z$.  
  Assume the field of fractions of $R$ is $F$.
 Let $\nu$ be a measure on $Z$, such that for any nonzero $r \in R$,
 \[ \int v(r) d\nu(r) \leq 0 \]
Then there exists a GVF structure on $F$, such that for any $t \in Term_F$ 
supported on $Z$ we have: $\int t = \int t(v) d\nu(v)$.  

\>{lem}

\prf  By the discussion above, We may replace $Z$ by the interior $Z^o$ of $Z$: we have  $R_{Z^o} =R_Z$, 
and  if $t$ is supported on $Z$, it is supported on $Z^o$.  Thus we will assume that $Z$ is open.

 Let $V$ be the $\Qq$-space of term functions on   $Val_F$. 
Among these, we have the functions $(f)$ for $f \in F \m (0)$, defined by: $f(v) = v(f)$.
Let ${Pr}$ be the $\Qq$-span of  $ \{(f): f \in F\}$.  
 Let $U$ be the subspace of functions
supported on $Z$.     For $\alpha \in U$,  let
$\phi( \alpha) = \int \alpha(v) d\nu(v) $. 

\claim{1}  Let $\alpha_1,\ldots,\alpha_n \in \Rr$ be linearly independent over $\Qq$, and $f_1,\ldots,f_n \in F^*$.  If
$\sum \alpha_i (f_i) \in U$, then each $(f_i) \in U$.   
 
 So $f \in R$ by definition of $R$, and hence by assumption, $\int v(f) d\nu(v) \leq 0$,
and hence $\int v(\a - (f)) d \nu(v) \geq 0$.

  For $\alpha \in {Pr}$, let $\phi(\alpha)=0$.  It follows from the Claim that $\phi$ can be extended to
  a linear map on $U + {Pr}$ with $\phi(\alpha)=0$ for $\alpha \in {Pr}$.   We again denote this by $\phi$.  
 
We define a cone  $\pos$ on $V$:  ${\pos}$ is the set of non-negative functions in $V$.      
It is clear  that $\phi$ is positive on $U+{Pr}$, i.e. $\phi(g) \geq 0$ if $g \in U+{Pr} \meet P$.

Let us check the boundedness condition of M. Riesz (\ref{mriesz} (2)). 
  Any element of $V$ is bounded by a finite sum of elements $h^+$,
with $h = f/g$ for some $f, g \in R$, so it suffices to check that the image of $(f/g)^+$ in $V$ is bounded by an element of $U$.  
 We  have
\[ (f/g)^+ \leq (f)^+ + (g)^-  = (f)^- + (g)^-  +(f)\]
now $(f)^- + (g)^-$ is supported on $Z$, so lies in $U$, and hence $(f)^- + (g)^-  +(f) $lies in $U+ {Pr}$.
 
By \lemref{mriesz}, $\phi$ can be extended to a positive functional on $\bV$.  

Hence there exists a GVF structure on $F$, agreeing with $\nu$ on terms supported on $Z$.  


\eprf

For simplicity, we state an approximation criterion over global fields (though it could easily be generalized);
thus we assume $S$ is finite.  
Let $a_l \in K^{alg}$.  We say that the $a_l$ approximate $P_\mu$ if  
for any formula $\phi$ over $K$, if for some subsequence the values $\phi(a_l)$ converge to $\alpha$,
then $P_\mu$ is consistent with $\phi(x)=\alpha$.  Equivalently, for any subsequence $a_j$, if $qftp(a_j/K)$
converges to a qf type $q$ then $q \vdash P_\mu$.  

\begin{lem} \label{galoisorbits} Let $a_l \in K^{alg}$.  Then the $a_l$ approximate $P_\mu$ iff
\begin{enumerate}
\item They are generic, i.e. only finitely many $a_l$ lie on any given proper subvariety of $X$.
\item They are near-integral away from $S$.  In other words, 
 \[  \sum_{p \in \Pi \m S} \int_{v| p} \max_{i=1}^n v^-(x_i) \to 0\]
\item The Galois orbit measures approximate $\mu_p$ for each $p$.
\end{enumerate}
\end{lem}

 The following theorem applies in particular to $K=k(t)[1]$ and finite extensions thereof, by \thmref{ec}.

By the characteristic measure of a finite set $O$, we mean $\mu_O= \frac{1}{|O|} \sum_{a \in O} 1_a$,
where $1_a$ is the characteristic function of $a$.

\begin{thm}  \label{fs9} Assume the {\em symmetric extension} $K^{alg}$ is existentially closed as a GVF.    Then 
 \cha$(A)  \geq 1$ iff there exists a sequence of Galois orbits $O_n$ that are near-integral away from $S$,
 generic on $X$, and such that for each $p$, if we fix an extension $K_p^{alg}$ of the 
 valued field $(K,P)$ to $K^{alg}$, the 
   characteristic measures of the $O_n$ on  the Berkovich space $X(K_p^{alg})$ approach a measure on $A_p$.  
  \end{thm} 
\prf This follows from \lemref{fs6} and \lemref{galoisorbits}. \eprf

 \end{section}

\begin{section}{Divisors}   
\label{divisors}
  Let $X$ be an irreducible normal projective variety over a field $k$,
$K=k(X)$ the function field.

Let $Div(X)$ be the group of Cartier divisors on $X$.  If $j: X' \to X$ is a surjective morphism, pullback
induces an injective map $Div(X) \to Div(X')$.     By $Div(K)$ we denote the direct limit of the groups $Div(X')$
over all birational morphisms $X' \to X$.  

A {\em global section} of a Cartier divisor $D$ on $X$ is a rational function $f \in k(X)$
   such that for any open $U$, if $D$ is represented by $(d)$ on $U$, then $f = r d$ on $U$ for some regular
   function $r$ on $U$.   Equivalently, $(f)+D $ is an effective Cartier divisor.

\ssec{Stable joins and meets of Cartier divisors}  \label{stablemeet}

Given   effective  Cartier divisors $E,E_1,E_2$ on $X$, represented by hypersurfaces with locally principal ideal sheaves $I,I_1,I_2$,
we say that $E$ is the {\em stable meet} of $E_1,E_2$ if $I$ is the ideal sheaf $I_1+I_2$ generated by $I_1 \union I_2$.  In general, $E_1,E_2$ may not have
a stable meet on $X$, since $I_1+I_2$ may not be locally principal.  But if we pull back to  divisors $E_1',E_2'$ on the blowup $X'$
at $I_1+I_2$, then $E_1',E_2'$ have a stable meet $E'$, and this remains the case in any further  pullback.
We write in this situation:  $E = E_1 \wedge E_2$.

 A similar definition with  sheaves of submodules of $K$ extends to all Cartier divisors.  It is moreover easy to see
 that $E = E_1 \wedge E_2$ iff $E+D = (E_1 +D) \wedge (E_2 + D)$, so one can reduce to the effective case.
 
Define $E = E_1 \vee E_2$ iff $-E = (-E_1) \wedge (-E_2)$; we then say that the {\em stable join} of $E_1,E_2$ exists
and equals $E$.   

It is clear that if $E$ is the stable meet  of $E_1,E_2$, then it is characterized as the greatest Cartier divisor
below both $E_1$ and $E_2$; dually if $E$ is the stable join of $E_1,E_2$, it is the smallest Cartier divisor
above both $E_1,E_2$.  Here we refer to the ordering:  $X \geq Y$ iff there exists an effective Cartier divisor
$Z$ with $X=Y+Z$.

There need not exist any Cartier divisor $E$
on $X$ with $E=E_1 \wedge E_2$.    For example, two distinct lines on $\Pp^2$ never have a stable  meet,
since their intersection is a point and not a divisor.   In fact even a stable meet with $0$ may not exist.  

However,  if $E_1,E_2$ are effective Cartier divisors, with ideal sheaves $I_1,I_2$, and $f:X' \to X$ is the result of blowing
up the ideal $I=I_1+I_2$, then it is clear that $f^*E_1 \wedge f^*E_2 = E$ where $E$ is the exceptional divisor.
By twisting with an appropriate very ample divisor, we can reduce to the effective case; so any two Cartier divisors do have
a stable join, once pulled back to an appropriate blowup.   Moreover once we have
   $E=E_1 \wedge E_2$,   the same relation persists in pullbacks to variety $X'$ dominating $X$.
Thus going to the limit, we obtain a statement that is very essential to the theory of GVF's:  

\begin{prop}  \label{lattice} $Div(K)$ forms a lattice under stable meet and join.  \qed
\end{prop}
 
 
 \ssec{The Picard group}

   $Pic(X)  $ is the group of   Cartier divisors on $X$, 
  modulo the {\em principal divisors}  $(f)$.     An element of $Pic(X)$ determines a line bundle on $X$.
  This shows that 
     $Pic$ is covariantly functorial; given a morphism $f: Y \to X$ of varieties, there is a natural map
$f^*: Pic(X) \to Pic(Y)$.

\ssec{Intersections of subvarieties with divisors}
An important  special case   occurs when $X$ is a subvariety of $Y$ and $f$ is the inclusion.  In this case,
we obtain a map from divisors on $Y$ to divisors on $X$.  If we further represent the divisor on $X$ using subvarieties
of $X$ of codimension $1$ and then view them as subvarieties of $Y$ of codimension $1+codim(X)$, this
map becomes the operation of {\em intersecting a subvariety with a divisor}.

Specifically,  a {\em very ample} divisor on $X$ can be represented by a subvariety $D$ on $X$ that intersects 
$Y$ transversally; further when $\dim(Y)\geq 2$, the intersection is an irreducible variety itself.  So \[ [D] \cdot [Y] = [D \meet Y] \]
When $\dim(Y)=1$, $D \meet Y$ is a finite set of isolated points, and we are interested in their number:
\[  [D] \cdot [Y] = |D \meet Y|  \in \Nn\]

\ssec{Weil divisors} 

We will generally work   with Cartier divisors, but sometimes the Weil viewpoint is useful.

 An {\em  integral Weil divisor}
  is a $\Zz$-linear combination of irreducible hypersurfaces of $X$.    More generally, the group $Z_k(X)$ of  $k$-dimensional  Weil cycles is the free Abelian group generated by the   irreducible subvarieties of dimension $k$.  $W_k(X)$ is the quotient by
  the $k$-dimensional divisors that are rationally equivalent to $0$.   $W^k(X)=W_{n-k}(X)$.  
  
   Principal Cartier divisors correspond to principal Weil divisors $(f) =  \sum_D {v_D(f)} D$.  (Here $v_D$ is the valuation associated with $D$;
  in an open affine intersecting $D$, $D$ is the zero locus of a minimal prime ideal; localizing, we obtain a Dedekind domain.)
  
There exists a canonical biliear map $Pic(X) \times W^k(X) \to W^{k+1}(X)$, the intersection product.
   
   \ssec{Global sections}
   
 With each real Weil divisor $D=\sum \alpha_i D_i$ ($\alpha_i \in \Rr$)  we have an associated $k$-subspace of   $K$, namely
  $L(D)= L(X,D) := \{0\} \union \{f \in K \m (0): \bigwedge_i v_{D_i} (f) \geq -\alpha_i \}$.  If $D$ is a Cartier divisor,
 then $L(D) = \{0\} \union \{f \in K:  (f) + D \geq 0 \}$.  The elements of $L(D)$ can be viewed as rational sections of the associated line bundle, and called sections of $D$.   When $D$ is a Cartier divisor, with corresponding sheaf
 $\tilde{D}$, one can identify $L(D)$ with $H^0(X,\tilde{D})$.

\ssec{Intersections of curves with divisors}  The $0$-dimensional cycles, up to algebraic equivalence, can be identified
with the integers $\Zz$ via the degree map.  At any rate we have the map $\deg: W_0(X) \to \Zz$,
$\deg(\sum m_i c_i) = \sum m_i$ where $c_i$ are points of $C$.  Composing with $\deg$ we obtain a  numerical
intersection product $W^1(X) \times W_1(X) \to \Zz$.   

The fact that the intersection pairing is well-defined amounts here to the {\em product formula}:  when $D=(f)$ is a {\em principal} divisor, we have $D \cdot C = 0$.

\ssec{Basic positivity property}  Let $D$ be an effective Weil divisor and $C$ a curve, not contained in the support of $D$.  
Then $D \cdot C \geq 0$, and equality holds iff $C \meet D = \emptyset$.  (The intersection multiplicities are non-negative.)

\ssec{Numerical equivalence}.  
     We consider several subgroups on the  group of 
  Weil-divisors, described via their sets of generators.   For us, numerical equivalence is the essential equivalence
  relation; we mention rational equivalence mostly in order to have a well-defined intersection product.    
  
  \begin{itemize} \item Rational equivalence:  $g: X \to \Pp^1$ a rational map; $g \inv(0) - g \inv(\infty)$.  
 (We can view $g$ as a morphism on  $X' = X \m W$ for  some codimension $2$ subvariety $W$; hypersurfaces of $X$
 and of $X'$ can be identified.)
\item Algebraic equivalence:  $g: X \to C$ for some algebraic curve; or more generally, pushforwards of such under
generically finite morphisms.  
\item  Numerical equivalence;   $Num_0 = \{D: (\forall C ) D \cdot C = 0 \}$.
\end{itemize}
%
 \begin{example}  \rm Let $X$ be a smooth surface, $a \in X$, and $f: Y \to X$ the blowing up of $X$ at $a$.
 Let $E=f\inv(0)$ the exceptional  divisor.   Let $C$ be a  a curve on $X$ with $a$ a nonsingular point of $C$.
 Then \begin{enumerate}
 \item $f^*C = E+C'$, where $C'$  is the strict transform of $C$.  
 \item  $E^2 <0$, i.e. for the inclusion morphism $i: E \to Y$ we have $i^* E <0$.  
 \end{enumerate}    
 To see   (2), let $D$ be a divisor on $X$ rationally equivalent to $C$, with $0 \notin D$.  
 As $f \inv(D) \meet E = \emptyset$,  $[f^* D] = [f^* C]$, and by (1) we have
 $(E+C' ) \cdot E = f ^* (D) \cdot E = 0$,  so $E^2 = -C' \cdot E$.  But $C',E$ intersect properly with nonempty intersection,
 so  $C' \cdot E >0$ and $E^2 < 0$. 
\end{example}

\ssec{The N\'eron-Severi group} 

The N\'eron-Severi group  is $W^1(X)$ modulo the subgroup of divisors algebraically equivalent to zero.

\begin{thm}[Lang-N\'eron, 1951]  $NS(X)$ is finitely generated.  In particular, $N^1(X) := NS(X) \tensor \Rr$ is
a finite-dimensional real vector space. \end{thm}

Note that the   Mordell-Weil theorem on finite generation of the rational points of a simple non-isotrivial
Abelian variety  over a function field   is a special case.  For instance, if $X$ is a surface admitting a morphism to $\Pp^1$, such that the generic
fiber is an elliptic curve $E_K$, that does not descend to $k$,  then any $K$-rational point of $E_K$ has a Zariski closure
over $k$ which is a curve; this defines a homomorphism from $E_K(K) $ to $NS(X)$ which is injective.   Hence finite generation 
of $NS(X)$ implies finite generation of $E_K(K)$.

The Mordell-Weil theorem is usually proved using heights (in any book on diophantine geometry).  \\
Given this special case, the N\'eron-Severi theorem yields, with some care, to an  induction on dimension.

\begin{rem}  A fourth equivalence relation   is obtained by a homological representation.
In characteristic zero, the finite dimensionality of $N^1(X)$ can also be seen via the representation in the homology group,
along with a proof based on real geometry (e.g. o-minimality) of the finite topological type of $X(\Cc)$.  
\end{rem}

\begin{defn} \label{nef} The closure of the ample cone in $N^1(X) $ will be called the {\em nef} cone.  \end{defn}

\begin{rem}  It is clear from the above definition that if $D$ is nef, then $D \cdot Y \geq 0$ for any curve $Y$ on $X$,
and moreover $D^k \cdot Y \geq 0$ for any $k$-dimensional subvariety $Y$ of $X$.   In fact, the converse is also true:
if  $D \cdot Y \geq 0$ for any curve $Y$ on $X$, then $D$ is nef.  This follows from Kleiman's theorem, and builds
upon Nakai-Moishezon; see \cite{debarre} for a very nice exposition.   We will however simply use the characterization 
as the closure of the ample cone.  \end{rem}

\begin{lem}[\cite{bfj} Prop. 2.3] \label{monotone}
 If $d_1,\ldots,d_k$ and $c_1,\ldots,c_k$ are nef 
classes in $N^1(X)$, 
 $k=\dim(X)$, $c_i \leq d_i$, 
then $c_1 \cdots c_k  \leq d_1 \cdots d_k  $.
\end{lem}

\prf  By continuity and multilinearity we may assume the $d_i$ are ample, or even very ample. 
  We may replace
the $c_i$ by $d_i$ one by one;   in this case, $d_1 - c_1$ is pseudo-effective, so it  is the limit of
effective classes $e$; as $e$ is effective and $c_i$ is very ample for $i>1$, it is clear that 
$e \cdot c_2 \cdots c_k \geq 0$ and so    $(d_1-c_1)c_2 \cdots c_k \geq 0$.  \eprf

Let us also define other cones we will use later.

\begin{notation} \label{notation}  Let $W$ be a subset of $\Rr^n$.   \begin{itemize}
\item The interior of a set $W$ will be denoted $W^o$. 
\item $\Rr^{\geq 1} W := \{rw: r \in \Rr, r \geq 1, w \in W \}$.
\item  The smallest closed, convex set containing $W$ will be denoted $[W]$.  
\item 
 A function $\phi$ on a subset of a real vector space is {\em  positively 1-homogeneous} if for all real $\a>0$ and $x \in \dom(\phi)$ 
  we have  $\a x \in \dom(\phi)$ and $f(\a x)=\a f(x)$.
  \item  For a real vector space $V$, $V^*$ is the dual space, and the natural pairing $V \times V^* \to \Rr$
  is denoted by $x \cdot y$ or $(x,y)$.
\end{itemize}
 \end{notation}


$\psf(X)$ = closure of the effective cone. 

$\psf(X)^o$ = the big cone,  interior of $\psf(X)$ = classes with positive volume.

 $N_1^+(X)$  is defined to be the closed dual of $\psf(X)$.  
 
 $N_1^{eff}(X)$ is the cone in $N_1(X)$ generated by the classes of  irreducible curves.  
 It may be bigger than $N_1^+(X)$.

\ssec{Functoriality of the cycle groups}   \label{pullbacks} Let $\pi: X \to Y$ be a morphism of  varieties.
If $\pi$ is {\em dominant}, there is a natural pullback map of Cartier divisors.  If $\pi$ is not dominant,
the trouble is that a Cartier divisor represented by $(f)$ may not pull back to a Cartier divisor, since $f \circ \pi$
may be zero.   But there is still a natural pullback of line bundles, or Cartier divisors up to principal divisors.  
This induces a  pullback map $N^1(Y) \to N^1(X)$.

In the case of a generically finite morphism $\pi: X \to Y$, we can also define a pullback map
$\pi^*: N_1(X) \to N_1(Y)$, by dualizing the pushforward map $\pi_*: N^1(X) \to N^1(Y)$.   
We have $\pi_* \pi^* (c) = \deg(\pi) c$.  

When $X,Y$ are smooth we can more generally define a pullback map $\pi^*:  N^k(Y) \to N^{k}(X)$   by the formula:
\[ \pi^*(u) = \pi_X ((X \times u) \cdot \pi) \] 
where $X \times u$ is the obvious class in $X \times Y$,
and $\pi$ here denotes the graph of $\pi$ as a subvariety of $X \times Y$.
Note that if $w \in N^l(X)$ is represented by $W$, and  $u$ is represented by a cycle $\sum \alpha_i D_i$ with  
  $\pi \inv(D_i) \meet W=\emptyset$, then $\pi^*(w) = \emptyset$.  (The intersection 
  $(Y \times u) \times \pi$ is proper, and the intersection of that with $(X \times D_i)$ is 
  empty, hence proper.)  
 
 We will often use the {\em  projection formula}.
 
\begin{lem}[Projection formula,   \cite{debarre} 1.9,1.10]  \label{projection}
 Let $\pi: Y \to X$ be a surjective morphism between proper varieties. Let $D_1,\ldots,D_r $ be 
 Cartier divisors on $X$,   $r \geq \dim(Y)$.  Then 
 \[ \pi^*D_1 \cdot \ldots \cdot \pi^* D_r = \deg(\pi) (D_1 \cdot \ldots \cdot D_r) \]
More generally if $W$ is a subvariety of $Y$ and $r \geq \dim(W)$, then 
 \[ \pi^*D_1 \cdot \ldots \cdot \pi^* D_r \cdot W = \deg(\pi|W) (D_1 \cdot \ldots \cdot D_r \cdot \pi W) \]
 \end{lem}
When $r> \dim(X)$, we have $D_1 \cdot \ldots \cdot D_r=0$; so both sides are zero unless $r=\dim(Y)=\dim(X)$.
Note that the `moreover' follows from the first statement applied to $\pi|W: W \to \pi W$.  
When cycle classes on $X$ are considered, notably when $W$ is a curve and has a class
$[W] \in N_1(Y)$, one defines $\pi_*[W] = \deg(\pi|W) [ f_* W]$, and 
the right hand side is written $D_1 \cdot \ldots \cdot D_r \cdot \pi_*[W]$.  

When $f: X \to Y$ is flat, we   have a pullback map of cycles, see \cite{fulton}.  

\ssec{Weighted curves}  It will sometimes be convenient to speak of weighted curves.
These are just formal pairs $(\alpha,C)$ with $C$ a curve on $X$ and $\alpha >0$.
The corresponding class in the cycle class group is $\alpha [C]$.  
 
\end{section}

\begin{section}{Convexity}  \label{convexity-s}

Recall the notation of \ref{notation}.   The natural pairing 

\begin{lem}  \label{convexity1} Let $\phi$ be a convex differentiable  function on an open set $U \subset \Rr^n$.  Let 
 $F=d \phi: U \to (\Rr^n)^*$ the differential.  Then $F$ is increasing in the sense that
for all $x_1,x_2 \in U$ we have
 \[(F(x_1)-F(x_2), x_1-x_2) \geq 0 \]
 If $x_1 \neq x_2$ and 
 $\phi$ is strictly convex on the line through $x_1,x_2$, then strict inequality holds.  
\end{lem}

\prf    By restricting to the line through $x_1,x_2$ we reduce to the one-dimensional case, where the statement is
that a (strictly) convex differentiable function has a (strictly) increasing derivative. \eprf

\begin{thm} \label{convexity}  Let $U$ be a (nonempty) open convex cone  in $E=\Rr^n$.
Let $\phi$ be a  continuous real-valued function on $E$,  differentiable on  $U$.    Let 
\[ F=d \phi: U \to E^* \]
 be the differential.  Let
\[ C_1 =  \{u' \in E^*: (\forall u \in U  ) (u,u') \geq \phi(u) \}\]  
and let ${C_1^o}$ be the interior of ${C_1}$.  
 Assume:
  \begin{enumerate}
\item $\phi$ is concave on $U$. 
  \item   $\phi$ is positively 1-homogeneous
  \item  $\phi > 0$ on $U$, and $\phi=0$ on $E \m U$.
   \end{enumerate}
Then \[  C_1^o \subset \Rr^{\geq 1} F(U) \subset C_1 \]

 \end{thm}
 
\prf    We first note a few facts about $F$:
\begin{enumerate}
 \item  $\phi(x+y) \geq \phi(x)+\phi(y)$.
  \item    $\phi$ is strictly increasing on $U$ with respect to $<_U$.
\item $F(U) \subseteq C_1$.  
 \item Let $p \in cl(U)$, $p \neq 0$, and $l \in {C_1}^o$.   Then $l(p)>0$.  
 
\end{enumerate}
 
Proof of (1-4):  
  \begin{enumerate}
 \item This is immediate from concavity and 1-homogeneity.
\item  If $x,y \in U$ and $x <_{U} y$, then $y=x+z$ for some $z \in {U}$; so as $\phi(z)>0$, 
$\phi(y) \geq \phi(x)+\phi(z) > \phi(x)$. 
\item   For $u,u' \in U$ and $t>0$ we have by (2) $ t\inv( \phi(u+tu') - \phi(u) ) \geq \phi(u')$.  Thus 
\[ u' \cdot F(u) = \frac{d}{dt} \phi(u+tu') (u)= \lim_{t \to 0^+} t\inv( \phi(u+tu') - \phi(u) ) \geq \phi(u') \]

Here is a more visual version of the proof of (3):  if $c = F(u)$ with $u \in {U}$, then for appropriate $b$, the graph of $y=c \cdot x +b  $   is tangent to the graph of $\phi$ at $u$, and 
  so by concavity of $\phi$ lies above it, i.e. $cx +b \geq \phi(x)  $ for $x \in U$.  By homogeneity,
 for any $x \in U$ and $t>0$ we have $ctx+ b \geq \phi(t x) = t \phi(x)$, so $cx + b/t \leq \phi(x)$; letting $t \to \infty$
 we find that $cx \geq \phi(x)$ on $U$, so $c \in {C_1}$. 
%

\item   
Let $p \in cl(U)$, $p \neq 0$.   Then $l \cdot p \geq \phi \geq 0$ for $l \in {C_1}$.  
By \lemref{interiorofdual}, $l \cdot p > 0$ for $l \in {C_1}^o$.  
 \end{enumerate}
 
 \

Since $\Rr^{\geq 1} {C_1} = {C_1}$, it follows from (3) that  $\Rr^{\geq 1} F({U}) \subseteq {C_1}$.   It remains to show
that  $ {C_1}^o \subseteq \Rr^{\geq 1} F({U}) $.

Let $l \in {C_1}^o$.
By (4) above, $\phi/l$ is continuous on $cl(U) \m (0)$, and invariant for the multiplicative action of $\Rr^{\geq 1}$.
Let $B$ be a compact ball in $E$,  with boundary sphere $S$,  and $p$ a point of $S \meet cl(U)$ with $(\phi/l) (p)=\gamma$ maximal;
we have $\g \leq 1$ since $\phi \leq l$ on $U$ and hence on $cl(U)$; also, $\phi>0$ on $U$ implies $l>0$
on $U$ so $\phi/l >0$ on $U$, and thus $0 < \g$ since $S \meet U$ is nonempty.  We have $\phi \leq \gamma l$
on $S \meet U$ and hence by $\Rr^{>0}$-invariance, on $U$.  

Now  $\phi(p)=\gamma l(p)$; as  $l(p) \neq 0$ and $\gamma >0$ we have $\phi(p)>0$, so $p \in U$, and  
 $p$ is a local minimum of $\gamma l - \phi$.  Hence  $\phi,\gamma l$ are tangent at $p$, i.e.   $d (\gamma l - \phi(p))=0 $
so $\gamma l = d \phi(p) = F(p) \in F(U)$.    Thus $l = \gamma \inv (\gamma l) \in \Rr^{\geq 1} F({U})$.

 Thus $ C_1^o \subset \Rr^{\geq 1} F(U) \subset C_1$.
\eprf
\begin{cor}\label{convexity-c}    With the hypotheses of \thmref{convexity}, let ${C}$ be the closed dual cone:
\[ {C} = \{u' \in E^*: (\forall u \in U  ) (u,u') \geq 0 \}\] 
and let $C^o$ be the interior of $C$.  Then $C = \Rr^{>0} C_1$.  Thus 
\[ C^o \subset \Rr^{>0} F(U)  \subset  {C}  \]  
 \end{cor}
 \prf     
By  \lemref{interiorofdual}, 
If $l \in C^o$ then $l>0$ on $\bar{U}$; in this case by compactness of the sphere, $\kappa l > \phi $ on $\bar{U}$
for some $\kappa >1$; so $\kappa l \in C_1^o \subset \Rr^{\geq 1} F(U)$.  Thus $l \in \Rr^{>0} F(U)$.    The other inclusion
is immediate from \thmref{convexity}. 
\eprf  

\begin{rem} \begin{enumerate}  \item  We saw in (3) of the proof of \thmref{convexity} 
that $v \cdot F(u) \geq \phi(v)$.  
It  follows from the 1-homogeneity of $\phi$  that equality holds for $u=v$:  
\[ u \cdot F(u) = \frac{d}{dt} \phi(u+tu) = \phi(u) \]
\item Let $a \in U^o$.  Then 
  $F(a)$ is the unique element $w \in C_1$ satisfying $a \cdot w = \phi(a)$.  Indeed assume $w$ is such an element.
For $v \in E$ and for  sufficiently small real $t$ we have $a+tv \in U^o$, so 
\[ (a+tv) \cdot w \geq \phi(a+tv)  = \phi(a) + tv \cdot F(a)   + o(t) \] 
Since $a \cdot w = \phi(a)$ this gives for small $t>0$ the inequality
$t( v \cdot w ) \geq t v \cdot  F(a) $, and for small $t<0$ the opposite one.  We thus
 obtain $v \cdot w = v \cdot F(a)$, and thus $F(a) = w$.  
\ 
 If for $b \in E^*$ we define $\phi^*(b) = \inf \{ (u,b) / \phi(u) : u \in U\}$, so that $b \in C_1$ iff $\phi^*(b) \geq 1$,
then $a \cdot b \geq \phi(a) \phi^*(b)$, and for fixed $a \in U^o$, equality holds on a unique ray of $b$.


%
%
\item  We will apply the theorem with $\phi= \vol ^{1/n}$ for a certain function $\vol$, and
set $\psi = \frac{1}{n} d \vol$.   Then (1) translates to $u \cdot \psi(u) = \vol(u)$; in 
  the setting of \cite{bdpp} and \cite{bfj} this is the `orthogonality relation'.   

\item  Similarly, (2) yields that $\psi(a)$ is the unique element $b$ of $C_1$ such that
\[   a \cdot b = \vol(a) \]

   \end{enumerate}
 \end{rem}

\begin{example}\label{convexnotopen}  Let $C$ be a compact convex subset of $\Rr^n$.  Let $L: \Rr^n \to \Rr^m$
be a linear map.  Then $L|C: C \to L(C)$ need not be an open map.  \end{example}
. \\
\prf  $L$ will be the projection from $\Rr^3$ to $\Rr^2$, $L(x,y,z)=(x,y)$.  Let 

\[ C=\{(x,y,z) \in [0,1]^3:   x^2 \leq y,   x^4  \leq zy^2  \} \]
   Let $p=(0,0,0)$.  Then $p$ projects to 
$(0,0)$, but a neighborhood of $p$ will have no images strictly on the curve $y=x^2$ with $x>0$.   
Convexity can be checked by intersecting with lines $l$.   If $L(l)$  passes through $(0,0)$, 
i.e. it has the equation $y=mx$, $m \geq 1$, we have above $l$ the set $ x^2 \leq m^2 z, \, x,y, z \leq 1$ which is clearly convex.  Away from $(0,0)$, differentiating  $x^4y^{-2}$ twice we see that it has positive definite Hessian, so the set of points above it
is convex.


\eprf

\begin{lem}  \label{openconvex} Let $U$ be an open subset of $\Rr^n$, and let $C$ be a dense convex subset of $U$.
Then $C=U$.  In particular, an open convex set is equal to the interior of its closure.    \end{lem}

\prf  Say $0 \in U$; we have to show that $0 \in C$.  Let $(a_i)_{i=1}^{2^n}$ be points of $C$
in each of the $2^n$ quadrants; this is possible by density.  Then $0$ lies in the convex hull of
the points $a_i$.  Indeed otherwise there is a hyperplane $c \cdot x \sum_{k=1}^n c_k x_k = 0$ through $0$, such that
all the points $a_i$ lie to one side, i.e. satisfy say $c \cdot x_i <0$.  However there is some $i$ in the same
quadrant as $c$, i.e. $(a_i)_k>0$ iff  $c_k<0$.  Then $\sum_k c_k (a_{i})_k >0$, a contradiction. 

For the last statement, let $C$ be open and convex, and let $U$ be the interior of the closure of $C$; then $C$ is dense in $U$, so equals $U$.
 \eprf

  \begin{lem} \label{interiorofdual} Let $D$ be a closed cone in $E=\Rr^n$, $C$   the closed dual cone in the dual space $E'$, i.e. 
  \[ c \in C \iff (\forall d \in D) c \cdot d \geq 0 \]
   Let $C^o$ be the interior of $C$.  Then 
   \[ c \in C^o \iff  (\forall d \in D \m (0)) \, c \cdot d > 0 \]      \end{lem}
 
   \prf   Let $c \in C^o$ and $0 \neq d \in D$.  As $d \neq 0$, there exists $a \in E'$ with $a \cdot d > 0$.   
 As $c   \in C^o$, for small enough $\e>0$   we have $c - \e a \in C$;  so $(c-\e a) \cdot d \geq 0$,
 or $c \cdot d \geq \e a \cdot d >0$.
 
 Conversely, assume $a \in E'$ and $a \cdot b >0$ for all $b \in D \m (0)$.  Fix some Euclidean structure
 on $E$, let $S$ be the unit sphere and $P= S \meet U$.  Then $P$ is compact, so for some $\e >0$
 we have $a \cdot b > \e$ for all $b \in P$.  Similarly there exists a neighborhood $U$ of $0$ in $E$ such that $|c \cdot b| < \e$ for all $c \in U, b \in P$.  Then clearly $(a+c)\cdot b >0$ for all $c \in U$, so $a \in C^o$.
\eprf 

    \begin{lem} \label{coneproject} Let $L: \Rr^n \to \Rr^m$ be a surjective linear map; let $C$ be an open convex cone in $\Rr^n$.
  Then $L(C)=\Rr^m$ if $C \meet \ker L \neq \emptyset$.  Otherwise, $L(C)$ is an open convex cone in $\Rr^m$;
  and $L(cl(C))=cl(L(C))$.  Moreover, $L| cl(C): cl(C)  \to cl(L(C))$ is an open map at $0$.  \end{lem}
  \prf    If $b \in C \meet \ker L$ and $c \in C$, then since $C$ is open, there exists $a \in C$ on the line from 
  $c$ to $b$, but beyond $b$, so that $b$ lies between $a$ and $c$.  Then $0$ lies between $L(a)$ and $L(c)$,
  so that as $L(C)$ is a cone, the full line lies in $L(C)$, and in particular $-L(c) \in L(C)$.  This shows that $L(C)$
  is a subspace of $\Rr^m$.  
  
  Assume now that $C \meet \ker(L)= \emptyset$.  It is clear that $L(C)$ is an open cone in $\Rr^m$,
  and $L(cl(C)) \subseteq cl(L(C))$.  The other direction can be sen using   \lemref{interiorofdual}.

   Next let $D=cl(C)$, and $\lam=L|D$; we show that $\lam: D \to  L(D)$ is an open map at $0$.  Put some Banach space structure on
 $\Rr^n$, and the quotient structure on $\Rr^m$; let $B_\alpha$ be the ball of radius $\alpha$.
 Then for some $\alpha>0$, $L(B_\alpha \meet D) $ contains $B_1 \meet LD$.  The reason is that the unit sphere
 $S_1$ in $\Rr^m$ is compact; each point $p$ of $LD \meet S_1$ is the $L$-image of some point of $D$, and  hence
 some neighborhood of $p$ lies in $L( D\meet S_{\alpha(p)})$ for some $\alpha(p)$; by choosing a finite subcovering we can find $\alpha$ such that $LD \meet S_1 \subset L(D \meet S_\alpha)$.  
 But then $LD \meet S_{\beta/\alpha} \subset L(D \meet S_\beta)$ for all $\beta >0$.  This shows  that $\lam$ is open at $0$.

   \eprf
    
\begin{lem} \label{convex-measure} Let $C$ be a closed convex subset of $\Rr^n$.  Then the boundary of $C$ has Lebesgue measure zero.
If $C$ has nonempty interior and finite measure, it is compact.
\end{lem}

\prf For the first point, we may assume $C$ is bounded, since the boundary of $C$ is contained in the union of
the boundaries of $C \meet B_n$, where $B_n$ is the ball of radius $n$.  Also, we may assume $0 \in C$.
Then for any $0<\alpha <1$ we have $\alpha C \subset C^o$, i.e. $\alpha C $ is contained in the interior of $C$.
As $\vol(\alpha C) = \alpha^n \vol(C)$, we see that $\vol \bd C \leq (1-\alpha^n) \vol(C)$.  Letting $\alpha \to 1$
we see that $\vol \bd C = 0$.

Assume now $C$ has non-empty interior; so it contains $a+B_\e$ for some $a$ and some open ball $B_\e$ of $0$.
If $a_1,a_2,\ldots$ are points of $C$ then $C$ contains $(a_i+a)/2+B_{\e/2}$; if the volume of $C$ is finite then
two of these balls must intersect, so $d(a_i,a_j) < \e$ for some $i,j$.   It follows that $C$ is bounded. \eprf

\begin{thm} [Khovanskii] \label{khovanskii}   Let $F$ be a finite subset of a $\Qq$-vector space $B$.  Let 
$S,A,C$ be respectively  the subsemigroup, subgroup, and rational convex cone  generated by $F$:\[ S = \{ \sum m_i a_i:  m_i \in \Nn, a_i \in F \} \]
\[ A =  \{ \sum m_i a_i:  m_i \in \Zz, a_i \in F \} \]
\[ C = \{\sum m_i a_i:  m_i \in \Qq^{>0}, a_i \in F \} \] 
  Then for some $s \in S$,
\[  A \meet (s+C) =  s + A \meet C \subseteq S \subseteq A \meet  C \]  \end{thm}

\prf  We may assume $B$ is also generated by $F$.  
Let $E = A \tensor \Rr = B \tensor_\Qq \Rr$.  Then $A$ is a discrete subgroup of the  finite-dimensional $\Rr$-space $E$.   Let  $Y = \{ \sum_{i=1}^m r_i f_i:  -1 \leq r_i \leq 0\}$.  Then $Y$ is compact.  As 
 $A$ is discrete, $A \meet Y$ is a finite set; we can find $s \in S$ such that $(A \meet Y)+s \subset S$.   
Let $a \in A \meet C$; so
$a = \sum \a_i f_i$  for some $\a_i \geq 0 $.   
Let $[\a_i]$ be the least integer $\geq \a_i$, and  $r_i = \a_i - [\a_i]$. 
Then $\sum r_i f_i \in A \meet Y$; so $s+ \sum r_i f_i \in S$.  Certainly $\sum [\a_i] f_i \in S$.
So $s+ a \in S$.   
\eprf   
%

We remark that if we identify $A$ with $\Zz^n$, then
$S$ is definable in Presburger arithmetic;  
hence globally it is defined by a combination of inequalities and congruences.

\begin{prop}  \label{sgroups2}  Let $\Lam$ be a finitely generated, torsion free Abelian group.  For $k \in \Nn$ let $S_k$ be   subsets of $A$ with  $0 \in S_1$ and $S_k+S_l \subset S_{k+l}$.  Let $S=\union_{k=1}^\infty S_k $.  
  Let $C_1=[\union_k S_k/k]$ be the closed convex hull of 
$\{ s/k: s \in S_k \}$ in $E=\Lam \tensor \Rr$.  If $K$ is a compact subset of $C_1^o$, then for all large enough $m$,
\[ K \meet \frac{\Lam}{m} = K \meet \frac{S_m}{m} \]
\end{prop}

\prf 
Let $\bar{S}=\{(n,x): x \in S_n  \} $; we view $\bar{S}$ as a subset of $\Rr \times E$.  It is clear that $\bar{S}$ forms a semigroup  of $\Zz \times \Lam$.   The group generated by $\bar{S}$ includes $\Zz \times 0$ (since $0  \in S_1$), 
and maps onto $\Lam$, so it equals  $\Zz \times \Lam$.  Since  $\Zz \times \Lam$ is a Noetherian $\Zz$-module, 
it is generated by a finite set $F$; let $W$ be the subsemigroup of  $\bar{S}$ generated by $F$.  Let $C$ be the convex
 cone generated by $F$.    By enlarging $F$ we may arrange that $\{1\} \times K \subset C$.  (The interior $C_1^o$ of $C_1$ equals the union of the interiors of the convex hulls $[F']$ of finite subsets $F'$ of $\union_k S_k/k$, using 
 \lemref{openconvex}; as $K$ is compact, $\{1\} \times K \subset  [F']^o$ for some such $F'$.)
 By \thmref{khovanskii}, there exists $w \in W$ with 
 \[   (\Zz \times \Lam) \meet C+w  \subseteq W\]

Let $K'= \{1\} \times K$.  
  As $K'$ is compact and $K' \subset C$, there exists a symmetric neighborhood $U$ of $0$ in $\Rr \times  E$ such that 
$ K' + U \subset C$.  
Let $m$ be large enough so that $w/m \in   U$.    If $b \in K \meet \Lam/m$, then $(1,b) - w/m \in C$, so $(m,mb) \in w+C$
and hence $(m,mb) \in W \subset \bar{S}$.  So $mb \in S_m$ and thus $b \in S_m/m$.  \eprf


\begin{lem}  Let $C$ be a nonempty  open cone in $\Rr^n$, with closure $\bar{C}$.  Then $a \in \bar{C}$ iff $a+C \subseteq C$. 
 \end{lem}
 \prf  
  Pick $c \in C$.  If $a+C \subseteq C$, pick $c \in C$; then  $a+tc \in C$ for all $t \in (0,1)$, so $a \in \bar{C}$.   
 Conversely, if $a \in \bar{C}$, let $U$ be a neighborhood of $0$ with $c+\bar{U} \subset C$; then
 $a+c+\bar{U} \subset \bar{C}$; so $a+c$ lies in the interior $C$ of $\bar{C}$.  \eprf

\begin{lem}[Tropical  Stone-Weierstrass]  \label{stonew} Let $Z$ be a compact Hausdorff topological space.  Let $A$ be a 
subset of the space of continuous functions $Z \to \Rr$ containing the constant
function $1$, closed under addition,   scalar multiplication and $\min$, and separating points.   Then $A$ is dense in
$C(Z)$ in the uniform topology.  \end{lem}

\prf    Let $B$ be the uniform closure of $A$.   If $a \in A$ and $p: \Rr \to \Rr$ is a convex function of the form
$p(x) = \max_{k=1}^n (a_k x + b_k)$,  it is clear that $p \circ a = \max_{k=1}^n (a_k a + b_k) $ lies in $A$.  It follows that $B$ has the same closure property.  

 On any compact
subset $W$ of $\Rr$, let $s_W: W \to \Rr$ be the function $s_W(x)=x^2$; then $s_W$   is the uniform  limit of a sequence of piecewise linear maps $p^j_W$, 
  obtained as
the maxima of  of certain  affine functions $ax+b$.   If $g \in B$, letting $W=g(Z)$, we see that  $p^j_W \circ g \in B$
and hence $g^2= lim_j p^j_W \circ g$ lies in $B$.  Writing  $xy=((x+y)^2-x^2-y^2)/2$, we see that $B$ is closed under multiplication.
By the usual  Stone-Weierstrass, $B=C(Z)$.  \eprf

  \end{section}
\begin{section}{Higher-order hyperbolic spaces and Hodge index theorem}  
\label{hodge-s}

We will encounter a generalization of the basic linear algebra of Hilbert spaces to higher order maps; we provide
here an elementary setting for this.

\begin{defn} \label{hyperbolic-def} An $n$-th order hyperbolic space is an $\Rr$- vector space $V$, an open cone $A$ and
    a symmetric  map $\cdot: A^n \to \Rr$  denoted $(a_1,\ldots,a_n)$, such that: \begin{itemize}
    \item  $  (a_1,\cdots,a_n) >0$ 
  \item Positive homogeneity in each variable: for fixed $a_2,\ldots,a_n$, and $t>0$, we have $(ta_1,a_2,\ldots,a_n) = t (a_1,\ldots,a_n)$.
  \item For fixed $a_2,\ldots,a_n$, the map $x \mapsto (x,a_2,\ldots,a_n)$ is concave.
       \item  For any $a_1,\ldots,a_{n-2},b,c \in A$, we have:
\[ (a_1,\ldots,a_{n-2},b,c)^2 \geq (a_1,\ldots,a_{n-2},b,b)(a_1,\ldots,a_{n-2},c,c) \]\
\end{itemize}
  \end{defn}
 
Let us write $|a| := |(a,\cdots,a)|^{1/n} \geq 0$.
\begin{example}[The bilinear, order $2$ case]  \label{negdef}
  \begin{enumerate} \ 
\item  Consider  a real vector space $V$  with a symmetric bilinear form of signature $(1,-1,\cdots,-1)$, or more generally $(1,-1,\cdots,-1,0,\cdots,0)$ (with precisely one
$1$ and at least one $-1$).  
The {\em future cone} $A$ is one of the two components of  $ \{x: (x,x) >0 \}$; say the one with $x_1>0$.  

On $A$, Cauchy-Scwhartz takes a `reverse' form,  
$(x,y)^2 \geq (x,x)(y,y)$;  whereas at any fixed time, i.e. on the orthogonal complement to any $a\in A$,
we have the usual Cauchy-Scwhartz inequality.

 
 To see that $(x, y) ^2 \geq (x,x)(y,y)  $ whenever $(x,x) >0$:  we may assume $(y,y) >0$ as well; renormalizing we may take $(x,x)=1$, and $y=x+n$ with $(x,n)=0$; as the orthogonal space to $x$ is negative definite, we have $(n,n) \leq 0$.   Now $(x,y)^2=1=(x,x)$ while $(y,y)=1+(n,n) \leq 1$.
\item Conversely, let  
   $(V,A,\cdot)$ be hyperbolic of order $2$, with   $\cdot$   bilinear.
Pick $a \in A$ with $|a|=1$.  Consider \[H=\{b: (a,b)=0 \}\]  
Then $H$ is negative-semidefinite:  otherwise we can find   $b \in H$ with $(b,b)>0$ and
(replacing $b$ by a small positive scalar multiple, so that $b$ is close to $0$)  with 
 $a+b \in A$,  yet 
$(a,a+b)^2 = 1 < 1+ (b,b) = (a,a)(a+b,a+b) $.  

\item  Going back and forth, from (2) to (1) and back,
we see that in (2) we may take any $a \in H$ with $(a,a)>0$; the orthogonal complement
is still negative-semidefinite.
\item   It is sometimes convenient to formulate Cauchy-Scwharz for two classes $x,y$, assuming $(a,a) >0$
and $y \cdot a = 0$, but without any assumption on $x$.  Then by
considering $x'=x- \alpha a$, where $a \cdot x = \alpha a^2$, one obtains 
\[  (x,y)^2 \leq (x^2 - \alpha^2  a^2) y^2 \] 
\item   A nondegenerate symmetric bilinear map $(\, ,\, )$ is hyperbolic
iff the function $|x|=(x,x)^{1/2}$ is concave on $A=\{x: (x,x)>0\}$.  Indeed assuming concavity, we have 
$(a+b,a+b) ^{1/2} \geq (a,a)^{1/2} + (b,b)^{1/2}$, or squaring,
$(a,a) + 2(a,b) + (b,b) \geq (a,a) + 2 (a,a)^{1/2}(b,b)^{1/2} + (b,b)$ which yields $(a,b)^2 \geq (a,a)(b,b)$.
\end{enumerate}
\end{example}
%
%
%

\begin{lem} \label{hh2}  For $x_1,\ldots,x_n \in A$ we have
 \[  (x_1,\ldots,x_n) \geq |x_1| \cdots |x_n|   \]
More generally, for $0 < k < n$, 
\[ (x_1,\ldots,x_k,x_{k+1},\ldots,x_n)^{n-k} \geq (x_1,\ldots,x_k,x_{k+1},\ldots,x_{k+1}) \cdot \ldots \cdot (x_1,\ldots,x_k,x_{n},\ldots,x_{n}) \]  \end{lem}

\prf By induction on $n-k$, being given  the case $n-k=2$.   Let us consider for example the case $n=3$.  
First we have:
\[ (a,b,b)^2 \geq (b,b,b)(a,a,b) \]
so
\[ (a,b,b)^4 \geq (b,b,b)^2 (a,a,b)^2 \]
Plugging in the dual of the first equation, i.e. $(a,a,b)^2 \geq (a,a,a)(b,b,a)$, and dividing by $(a,b,b)$ we obtain:
\[ (a,b,b)^3 \geq (b,b,b)^2 (a,a,a) \]
Now this in turn can be plugged into the relation
\[ (a,b,c)^2 \geq (a,a,c)(b,b,c) \]
to give
$(a,b,c)^6 \geq (a,a,a)^2 (b,b,b)^2(c,c,c)^2 $
and taking square roots we obtain the required formula. 

(See   \cite{lazarsfeld}.)

\eprf   
 
\ssec{log-concavity}  Let $f$ be a function on a cone in a real vector space, into the positive reals, and assume
$f(tx)=t f(x)$ for $t >0$.  Then $f$ is concave iff $\log f$ is concave.  The left to right implication is clear since $\log$ is
concave.  For the other, assume $\log f$ is concave; then whenever $0 < \a,\b$ and $\a+\b=1$,
we have $\log f(\a x+\b y) \geq \a \log f(x) + \b \log f(y) $; in particular when $f(x)=f(y)$ this implies that
$f(\a x + \b y) \geq f(x)$.  By positive homogeneity of $f$, for any $\g >0$ we have 
$f(\g \a x + \g \b y ) \geq f(\g x) = \g f(x)$; in other words for any $0<\a,\b$, if $f(x)=f(y)$ we have
 $f(\a x + \b y) \geq (\a+\b) f(x)$.   Given arbitrary $x,y$, find $\g$ such that $\g f(x)=f(y)$; let $x' = \g x$;
 let $\a' = \a/\g$; then 
    $f(\a  x + \b y) = f(\a' x'+\b y) \geq (\a' + \b) f(y)  = \a f(x) + \b f(y)$.  
    
\paragraph{{\bf The volume function on $A$.} } We define $\vol: A \to \Rr$ by $\vol(a)=(a,\ldots,a)$.  
\begin{lem} \label{logconcave1}
  $\vol^{1/n}$ is concave on $A$.\end{lem}
  \prf
  
Using concavity in each variable, we have
\[   \vol(a+b) =  (a+b,\cdots,a+b) \geq  \sum_k {{n}\choose{k} } (a,\cdots,a,b,\cdots,b) \]
(with $k$ instances of $a$, and $n-k$ of $b$.)
By the main statement of \lemref{hh2}, applied term-by-term, we obtain:
\[ \vol(a+b) \geq \sum_k {{n}\choose{k}}  \vol^{k/n}(a) \vol^{(n-k)/n}(b)  = (\vol^{1/n}(a) + \vol^{1/n}(b))^n \] 
or 
\[ \vol^{1/n}(a+b) \geq \vol^{1/n}(a) + \vol^{1/n}(b) \]
Using the positively homogeneity of $\vol^{1/n}$, this gives concavity of $\vol^{1/n}$  \eprf

  Thus   $- \log \vol$ is   convex on $A$. 
  
  \begin{example}   \label{bm1}  Let  $\Rr^n$ have basis $e_1,\ldots,e_n$.
  Consider the symmetric $n$-multilinear form with $(e_1,\cdots,e_n) = 1$, and $e_i^2=0$
(meaning $(e_i,e_i,x_3,\ldots,x_n)=0$ for all $x_3,\ldots,x_n$).  Let $A=\{\sum \alpha_i e_i:  \alpha_1,\ldots,\alpha_n>0\}$.  This is a generalized hyperbolic space; if we fix any $a_1,\ldots,a_{n-2} \in A$,
then $(x,y) \mapsto (a_1,\ldots,a_{n-2},x,y)$ has a kernel of dimension $n-2$ and signature $0,\cdots,0,1,-1$.
Note that if $c=\sum \alpha_i e_i \in A$, then $\frac{c^n}{n!}$ is the Euclidean volume of the rectangle
with sides $(\alpha_1,\ldots,\alpha_n)$.  Thus in this case \lemref{logconcave1} gives the concavity
of $\vol^{1/n}$ for rectangles.  This is a special case of Brunn-Minkowski.  \end{example}
 
 .   

 \ssec{The Hodge index theorem. } 
  
 \begin{thm} \label{hodge}[\cite{lazarsfeld} 1.6.1]  Let $X$ be a   projective   variety of dimension $n$.   
 Then $N^1(X)$ with intersection product, and the ample cone is an   order-$n$ multilinear hyperblic space. 
\end{thm}

\prf  

Let  $a_1,\ldots,a_{n-2}, a_{n-1}=b,a_n=c$ be ample divisors.  We have to show that 
$\vol(a_1) = a_1 \cdot \ldots \cdot a_1 > 0$ (this is clear), concavity in each variable (we have linearity), and that 
\[ (a_1 \cdot \ldots \cdot a_{n-2} \cdot b \cdot c)^2 \geq (a_1,\ldots,a_{n-2},b,b)(a_1,\ldots,a_{n-2},c,c) \]\
Now by Bertini, $a_i$ can be represented by an irreducible hypersurface $D_i$, such that
$D_1 \meet \cdots \meet D_k$ is a transversal intersection, resulting in an  irreducible subvariety of $X$
(for $k \leq n-2$).  Let $S= D_1 \meet \cdots \meet D_{n-2}$.  Then $S$ is a   surface.  Let $b' = D_{n-1} |S, c' = D_n|S$.
Then $(a_1,\ldots,a_{n-2},b,c) = b' \cdot_S c'$ and similarly for $(b,b)$ and $(c,c)$.  Thus the inequality reduces
to the case of a {\em surface}. 
 In this case, as noted above, it suffices to show that the signature is $(1,-1,\cdots,-1)$.
 
 If $\tS \to S$ is a resolution, then $N^1(X)$ embeds into $N^1(\tS)$ in a way that respects the intersection form,
 and maps an ample class to a  nef    class, so that it respects the closure of the ample cone, which suffices for the weak inequality.  Thus it suffices to show that $\tS$ has hyperbolic signature.
This is the classical Hodge index theorem for surfaces,  cf. \cite{hartshorne}.     If one uses the Okounkov approach, \corref{logconcavevolume}, 
this 
  requires
only the 2-dimensional Brunn-Minkowski.
  \eprf

 \begin{rem} \label{hardlef}
 When $X$ is a smooth projective variety of dimension $n$, $V=N^1(X)$, $A=$ the ample cone, the $n$-fold intersection pairing on divisors is   {\em nondegenerate}, i.e. whenever $v \in V$ and $(a_1,\ldots,a_{n-1},v)=0$ for all $a_1,\ldots,a_{n-1} \in A$,
then $v=0$.

 Note that the pairing  $N^1(X) \times N_1(X) \to \Rr$ is nondegenerate by definition; thus
 $\dim(N^1(X))=\dim N_1(X)$, and the content of the above statement is that every class in $N_1(X)$
 is in the image of $\cdot^{n-1}$.

   This can be seen as   a consequence of the Hard Lefschetz theorem
 for algebraic classes, asserting that even for a fixed ample $a$, the map $(u,v) \mapsto (a,\cdots,a,u,v)$
 is nondegenerate.   This can be interpreted as the Hessian of volume.
   We will see a  cone-adapted version of this   further down.   \end{rem}

  \begin{example} \label{castelnuovo} Let $S$ be a smooth projective surface admitting 
  morphisms $\pi_i: S \to C_i$ into curves, mapping onto $C_1 \times C_2$.  
     Let $P_i$ be the pullback of a point under $\pi_i$; 
  then $A=P_1+P_2$ is ample on $S$.  
The Hodge index theorem includes, and is easily seen to be equivalent to, the negative definiteness
of the intersection form on $V$,  the orthogonal space to $P_1$ and $P_2$.  
 For   any divisor $D$,  we have $D^* := D- (D \cdot P_1) P_2 - (D \cdot P_2) P_1 $ 
 in $V$.  
   Thus the Hodge index theorem asserts that 
$ (D^*,D^*) \leq 0$.    Explicitly, this reads:
\[ (D,D) \leq 2 (D \cdot P_1) (D \cdot P_2)  \]
This is the  Castelnuovo  inequality \cite{gcn}.  Weil developed the foundations of algebraic geometry between 1940-1948 precisely in order to provide a statement and proof of this inequality valid in 
  positive characteristic.     \end{example}

  \begin{example} Here is Weil's 1941 deduction of the Riemann hypothesis for curves from Castelnuovo,
  or as we will present it, directly from the Hodge
  index theorem for a product of curves.  Let $C$ be a curve over $\Ff_q$, $\Delta \leq C \times C$   the diagonal, 
  $\phi: C \to C$ the $q$-Frobenius morphism, 
   $\Phi=\Phi_q$      the graph of the  $\phi$.  
  Let   $^*$ be the operator of \exref{castelnuovo}.
    By \exref{negdef} we have 
  \[  (\Phi^* \cdot \Delta^* )^2 \leq  (\Phi^*,\Phi^*) (\Delta^*,\Delta^*) \]
  Here $(\Delta^*,\Delta^*)$ does not depend on $q$; 
  $(\Phi,\Phi) = q (\Delta,\Delta) $ using $\Phi = (\phi,Id)^* \Delta$; 
  $\Phi \cdot \Delta$
  is the number of points of $C(\Ff_q)$ that we seek; all other   terms are easily evaluated,    yielding Weil's theorem,
$|C(\Ff_q)|   = q + O(q^{1/2})$.     
 \end{example}

\begin{lem} \label{hodge3}  Let $S$ be a irreducible projective surface.  Let $\tS \to S$ be a surjective,
generically finite morphism,
with $\tS$ smooth.  
Let $e$ be an irreducible curve on $S$, 
 and let 
  $A,D \in N^1(S)$.  Assume $D \cdot A =0$ while $(A,A)=1$.
Then
\[ (e \cdot D)^2 \leq ( \te^2 -   (A \cdot e)^2 )  D^2 \] 
where $\te$ is any weighted curve on $\tS$ lifting $e$.  In particular if $D^2=0$ then $e \cdot D = 0$.

 \end{lem}
\prf   
  If we replace $D,A$ by their pullbacks, and $e$ by $\te$,
the same conditions hold and all the quantities in the displayed inequality are the same, by the projection formula (\cite{debarre} 1.9,1.10).    Thus we may assume $S$ is smooth 
and connected,
and $e=\te$.  By \remref{negdef} (4) we have
$ (e \cdot D)^2 \leq (e^2 - (A \cdot e)^2 ) D^2 $.
  \eprf

 The material from here to the end of the section is adapted from \cite{yuan-zhang} (part of Theorem 1.1; they   acknowledge   an idea of Blocki's).
  We begin with   a {\em relative} Hodge index theorem.    (In case $\pi^*(c)$ is represented
  by a subvariety of $X$, it is easy to deduce it from the Hodge index theorem on a
  plane section of that subvariety.    For a surface  mapping to a curve, 
let $l$ be an ample divisor with $(f,l)=1$, and note that $v_i \mapsto v_i - (v_i,l)f$
 is an isometry into the space orthogonal to $l$, with kernel $\Rr f$.  But we give a direct proof.)

\begin{prop} \label{rel-hodge} Let  $X,U$ be   smooth varieties with a morphism $\pi: X \to U$, of relative dimension $n$.  Let 
 $D$ be an effective divisor  of $U$;   and let $c \in N_1(U)$ be a movable class, represented
 by an effective curve $C$ disjoint from any specified subvariety of $U$ of codimension $>1$.  
Let $V$ be the subspace of $N^1(X)$ generated by the classes of the   irreducible components of $\pi \inv(D)$.
  Then  for nef $a_1,\ldots,a_{n-1}$, the pairing $V^2 \to \Rr$,
 \[ (v,v') = a_1 \cdot \ldots \cdot a_{n-1} \cdot v \cdot v'  \cdot \pi^*(c)   \]
 is negative semidefinite. 
 \end{prop}
  
 \prf 
   By approximation,  we may take the $a_i$ to be ample, with rational coefficients;
  multiplying by an integer, we may take them integral, and very ample.
  Let $v_1,\ldots,v_k$ be the distinct  irreducible components of $\pi \inv(D)$. 
  Let $D'=\sum \beta_j D_j$ be a divisor on $U$, numerically equivalent to $D$,
  such that $D_j \meet D$ has codimension $\geq 2$ in $U$ for any $j$. (We use smoothness of $U$ to move $D$, to a not necessarily effective divisor.)  Choose a representative
  $C$ of $c$ which is disjoint from $D \meet D_j$, for any $j$.  
  Let $\alpha_i$ be the multiplicity of $v_i$ in $\pi \inv(D)$; then $f=\sum \alpha_i v_i$ is numerically equivalent 
 to $\pi \inv(D')$.  Now $v_i \meet \pi \inv(D_j) \meet \pi \inv(C) = \emptyset $ for each $j$.
 Thus (see \secref{pullbacks}) $ v_i \cdot \pi^*(D') \cdot \pi^*(c) = 0 $; so $v_i \cdot f \cdot \pi^*(c)=0$,
hence  $(v_i,f)=0$.    For $i \neq j$,  $v_i$ meets $v_j$
 properly in an effective cycle; thus $a_1 \cdot \ldots \cdot a_{n-1} \cdot v_i \cdot v_j$
 is effective, hence so is $\pi_*(a_1 \cdot \ldots \cdot a_{n-1} \cdot v_i \cdot v_j)$; by the moveability property of $c$, it follows that $c \cdot \pi_*(a_1 \cdot \ldots \cdot a_{n-1} \cdot v_i \cdot v_j) \geq 0$; so $(v_i,v_j) = a_1 \cdot \ldots \cdot a_{n-1} \cdot v_i \cdot v_j \cdot \pi^*(c) \geq 0$.   Thus \lemref{pdc} applies.  
 \eprf

 \begin{lem} \label{pdc} 
  Let $V$ be a vector space over $\Rr$ with a 
symmetric bilinear form.   Assume $V$ is  generated by vectors $v_1,\ldots,v_n$, 
\[f=\sum_{i=1}^n \alpha_i v_i\] for some $\alpha_1,\ldots,\alpha_n>0$, and:
 \[v_j \cdot f =0, v_i  \cdot v_j \geq 0 \]
 for   $i \neq j \leq n$.  Then $V$ is negative semidefinite.  The kernel 
 $\{x: (\forall y) x \cdot y =0\}$ is generated by the vectors $\sum_{i \in c} \alpha_i v_i$,
 where $c$ is a connected component of the graph whose  edges are the pairs $(i,j)$ with    $v_i \cdot v_j \neq 0$.  \end{lem}

 \prf  Replacing $v_i$ by $\alpha_i v_i$, we may assume $f=\sum_i v_i$.     
Let $b_i \in \Rr, i=1,\ldots,n$, $g=\sum b_i v_i$; we have to show that $g^2 \leq 0$.   
As $v_j \cdot f = 0$ we have also $b_j ^2 v_j (\sum_i v_i) = 0$; 
summing over $j$  we find that 
\[\sum_{i,j} b_j^2 v_iv_j = 0 \] 
As $v_iv_j \geq 0$ for $i \neq j$, we have:
\[  g^2 = -\frac{1}{2}  \sum_{i,j=1}^l (b_i-b_j)^2 v_iv_j  = -  \sum_{i \neq j} (b_i-b_j)^2 v_iv_j \leq 0 \]
If equality holds, then $b_i=b_j$ whenever $v_i v_j \neq 0$, i.e. the coefficients are constant on each 
connected component. (And the converse is clear.)
  \eprf
  
   \begin{rem} \label{hodge2c} This also shows (assuming  $\pi \inv(D)$ is connected) that there are no unexpected linear  relations among the classes of the $v_i$ 
   in $N^1(X) / \pi^*( N^1(U))$ .   In fact, assuming for simplicity that  $\union_i v_i$ is connected, there are  no unexpected linear inequalities among the classes $[v_i]$.
   
 In the disconnected case, say when $D= D' \union D''$ with $D',D''$ disjoint,  connected, and numericaly equivalent in $N_1(U)$, the pullback to $\oplus_i \Rr v_i$ of the effective cone   is generated by the $v_i$ themselves, and in addition by expressions  $  \sum m_i' v_i'  - \sum m_i'' v_i$, the difference of the irreducible component decompositions of $\pi \inv(D'), \pi \inv(D'')$.  
 
      Suppose $[\sum \alpha_i v_i] =[e]$   represented by an  effective divisor $e$ on $X$; we have to show all $\alpha_i \geq 0$.  We may assume $e$ does not include any of the $v_i$ as components, or else we may subtract.  Let $J^+ = \{j: \alpha_j \geq 0$, $J^-  = \{j: \alpha_j < 0\}$, 
   $u=\sum \max(0,\alpha_i) [v_i] = \sum_{j \in J^+} \alpha_j [v_j]$.  
 Then $(u,u) \leq 0$  by negative semidefiniteness.  However $(u,e) \geq 0$, and  $(u,v_j)\geq 0$ 
 for $j$ with $\alpha_j<0$,   since $e,v_j$ are effective and these intersections are proper.  This implies
that $(u,\alpha_j v_j) \leq 0$ for  $j \in J^=$, and hence 
 $(u,u) = (u,e) - \sum  \geq 0$, so $(u,u)=0$;  in fact
   by connectedness of the graph and the description of the kernel of $(\cdot,\cdot)$, for $j \in J^+$,
$\alpha_j$  must be proportional to the multiplicity of $u_j$ in $\pi \inv(D)$; so either all $\alpha_j \geq 0$,
or all $\alpha_j \leq 0$.  In the latter case, pushing forward to $U$ we find that a negative multiple of $D$ is 
effective, a contradiction.  

   \end{rem}
 
 \ssec{Calabi uniqueness}
 
 \begin{lem}  \label{calabi} Let $W$ be a real vector space endowed with a symmetric,
 nondegenerate $n+1$-multilinear map $\cdot$.
  Let $V$ be a subspace of $W$,  $A$ a   subset of a coset of $V$.  
For $a=(a_1,\ldots,a_{n-1}) \in A$, let $b_a(v,v') = a_1 \cdot \ldots \cdot a_{n-1} \cdot v   \cdot v' $.  
Assume each $b_a$ is  negative semidefinite on $V$ ($a \in A$.)   
 
 Let $a_1,a_2 \in A$, 
 and assume  $v \cdot a_1^n = v \cdot a_2^n$ for all $v \in V$ (or even just for $v=a_1-a_2$).
  Then $a_1-a_2$ lies in the kernel of $b_a$ for any $a \in A^{n-1}$. 
 \end{lem} 
 
 \prf    Let $c=a_1-a_2$.    Consider first the case $n=2$.  
 
 We assume $c a_1^{2} = c a_2^{2} $.  So
 \[ 0 = c(a_1^2-a_2^2) = c^2 (a_1+a_2) \]
 But $c^2 a_i \leq 0$,  so $c^2a_1=c^2a_2=0$.
 
 By Cauchy-Schwarz (for $b_{a_i}$),  we have $vca_i  =0$ for any $v \in V$.
 
  Hence $c^2 v=c(a_1-a_2)v = 0$ for any $v \in V$.   As $c^2 a_1=0$, we have $ c^2 x= 0$ for $x \in a_1+V$
  and in particular $ c^2 a=0$ for all $a \in A$.  Again by negative semidefiniteness, $c$ lies in the kernel of $b_a$.
  
Now assume $n>2$.    
 We have by assumption:  
 \[  c a_1^{n} = c a_2^{n} \] So
 \[ 0 = c(a_1^n-a_2^n) = c(a_1-a_2)  (a_1^{n-1} + \cdots a_1^ia_2^{n-1-i} + \cdots a_2^{n-1}) \]
 But $c^2 a_1^ia_2^{n-1-i} \leq 0$;  so each term in the sum must vanish, in particular $c^2a_i^{n-1}=0$.   
 
 By Cauchy-Schwarz  we have $vca_i^{n-1}  =0$ for any $v \in V$.  As $a_i c a_i^{n-1}  =0$, we hve
 $e c a_i^{n-1}  =0$ for any $e \in A$.  Thus fixing any $e \in A$, the  $n$-multilinear map $x_1 \cdot \ldots \cdot x_{n-1} \cdot e$ satisfies the same conditions.   By induction, $c$ lies in the kernel of $b_{e_1,\ldots,e_{n-2},e}$ for 
 any $e_1,\ldots,e_{n-2} \in A$, which finishes the proof.

 \eprf

 \end{section}
 
\begin{section}{ Volume via section growth}  \label{volume-s}

So far, we defined volume {\em for ample divisors} via the intersection product.  Our official  definition 
of volume will instead use the asymptotic dimension of the space of sections of multiples of the divisor. 
This definition is also more general; an equally general definition in terms of positive  intersections
will be considered later (\secref{naive}), see  \cite{lazarsfeld}, 11.4.10, 11.4.11.  
 We begin by
relating the section growth with the intersection theory definitions  for ample divisors.

If $X$ is $0$-dimensional scheme,  let $|X|$ be the sum over all points of $X$ of the dimension
of the local ring.  

If $X$ is closed subvariety of $\Pp^N$, of dimension $d$, we define 
\[ \deg(X) = |X \meet L_1 \meet \cdots \meet L_{N-d}| \]
 where $L_i$ are hyperplanes chosen generically enough, so that
$\dim(X \meet L_1 \meet \cdots \meet L_{N-d})=0$.  (It will follow from \lemref{hilbertp} that the choice of the $L_i$
is irrelevant.)

On the other hand, let $A= \oplus H_n$ be the graded ring $k[X_0,\ldots,X_N]$; $H_n$ the homogeneous polynomials of degree $n$.    For a homogeneous ideal $I$, say 
 $ vol(I) = v$ if $\dim (H_n / H_n \meet I) = v m^d / d! + O(m^{d-1})$.   
 
 \begin{lem}  \label{hilbertp} Let $X=Proj(A/I)$.  Then $\deg(X), \vol(I)$ are both defined and are equal.  In fact, for large enough $n$, 
$\dim  (H_n / H_n \meet I)=p(n)$ for some    polynomial $p(n)$ of degree $d$, with leading term equal to $\deg(X) n^d / d!$.  \end{lem}

\prf  The case $\dim(X)=0$  follows easily from the definitions.  Assume the result
known below $\dim(X)$, and choose any $x \in H_1$ such that $x$ is not a 0-divisor in $A/I$.    
Let $I_n=I \meet H_n$, $I'=I+xA$, $I'_n=I' \meet H_n$, $X'=Proj(I')$.  
Note   that $(I+xA) \meet H_n = I_n + xH_{n-1}$ since the ideal $I$ is  homogeneous.  Thus
 \[ \dim (I'_n / I_n)  = \dim (I_n+xH_{n-1} / I_n ) =  \dim(xH_{n-1} / I_n \meet xH_{n-1}) 
= \dim (H_{n-1}/ I_{n-1})\]
using in the last step the assumption that $x$ is not a 0-divisor in $A/I$.  
So 
\[\dim H_n / I'_n  = \dim H_n / I_n - \dim H_{n-1} / I_{n-1} \]
Let $f(n)= \dim H_n / I_n$.   Then $\dim H_n / I'_n=f(n)-f(n-1)$.  By induction,
as $\dim(X')<\dim(X)$,  for large $n$ we have $\dim H_n / I'_n=q(n)  $ with $q$  a polynomial in $n$ of degree $d-1$, and leading
 term $\deg(X')=\deg(X)$.    We can find a polynomial  $p(n)$ of degree $d$, with the same  leading term $\deg(X)$, and with
 $p(n)-p(n-1) = q(n)$.  
 Then $p-f$ is eventually constant; adjusting by this constant we may assume $p-f$ is eventually $0$, i.e. $p(n)=f(n)$
 for large $n$.   
\eprf

Let $X$ be a variety of dimension $d$.  
We define the volume of an integral divisor $D$ to be 
\[vol(D) = d!  \limsup_m \dim L(X,D^m) / m^d \]

In particular $\vol(D)=v$  if $\dim H^0(X,D^m) = v {{m}\choose{d}} + O(m^{d-1})$.  
 
$D$ is {\em big} if it has nonzero volume.

\begin{cor} \label{vol-ample} If $D$ is very ample, then $\vol(D) = D^d$.   Hence for any $D'$, $L(X,mD') = O(m^{\dim(X)})$.   \end{cor}
\prf  Embed $X$ into projective space $\Pp^N$ via $D$.  Then $D^d$
equals the intersection number of $X$ with $d$ generically chosen hyperplanes on $\Pp^N$, and the 
statement follows from \lemref{hilbertp}.    The second statement follows by choosing an ample $D \geq D'$.  
\eprf

\begin{prop}[\cite{lazarsfeld} 2.2.15] \label{2.2.15}  Let $A,B$ be ample Cartier divisors on $X$.  Then 
\[\vol(A-B) \geq \vol(A) - n (A^{n-1}B) \]  \end{prop}
 
\prf  We may assume $A,B$ are very ample. Fix a countable base field of definition $k_0$ for the data.   
Using Bertini, let $E=E_1,E_2,\ldots$ be irreducible  hypersurfaces of $X$,  $E_i$ defined over some extension $k_i$ of $k_0$,
such that $(k_i,X,E_i)$ are all isomorphic over $k_0$, and $E_i$ represents $B$ as a Weil divisor.     Let
$A_i = A|E_i$; then the $n-1$-fold intersection product of $A_i$ on $E_i$ equals $ A^{n-1}  \cdot B$; by  \lemref{hilbertp} applied to $A|E$  on $E$, 
\[ h^0( E_i,mA_i) = h^0(E,m A|E) = A^{n-1}  \cdot B m^{n-1} /(n-1)! + O(m^{n-2}   ) \]
Again by 
  \lemref{hilbertp} for $A$ on $X$, 
\[ h^0(X,m(A-B)) = A^n m^n/n! +O(m^{n-1}) \]
We have an exact sequence
\[ 0 \to H^0(X, mA-mB) \to H^0(X,mA) \to \oplus_{i \leq m} H^0(E_i, mA_i) \]
Thus
\[ h^0(X,m(A-B)) \geq h^0(X,mA) - \sum_{i \leq m } h^0( E_i,mD|E_i) =(\vol(A) - n (A^{n-1}B)) {\frac{m^n}{n!}} +O(m^{n-1})  \]
 \eprf
 
 \begin{rem}[Continuity of volume] \label{vol-conts}  If we assume only that $B$ is   ample, with $B$ represented by $E$ as above,    the proof of \propref{2.2.15} still yields the inequality
 in the form: 
$\vol(A-B) \geq \vol(A) - n \vol_E (A|E) $.    Thus if $\beta=1/b$, we have
 $ | \vol(A-\beta B) - \vol(A) |  =  b^{-n} | \vol(b A-  B) - \vol(b A) | \leq b^{-n}( n b^{n-1} \vol_E (A|E)) = nb \inv \vol_E (A|E)$, 
 so  $ | \vol(A-\beta B) - \vol(A) | \leq n \beta \vol_E(A|E)$.  By choosing $n$ linearly independent very ample divisors,
 we see in this way that the volume function is continuous on $N^1(X)$; see \cite{lazarsfeld} 2.2.44.  
  \end{rem}
 
 \begin{rem} \label{15nef} From the continuity of volume \remref{vol-conts}, it is clear e.g. that \propref{2.2.15} is valid for {\em nef}
 and not only ample divisors. \end{rem}

\begin{defn}\label{big}   The open cone of elements of $N^1(X)$ of positive   volume is called the {\em big cone}.  \end{defn}
 A big integral divisor clearly has an effective multiple, so the big cone is contained in the interior of the pseudo-effective cone.
 Conversely if $D$ lies in the interior of the  effective cone, equivalently by \lemref{openconvex} of the pseudo-effective cone,  it can be written as ample + effective.  (Indeed for any ample $Q$, for some $\e>0$,
  $D - \e Q$ is effective.)  So some $mD$ lies above a very ample divisor, and hence has at least as many sections;
  thus $mD$ and hence $D$ are big.      We denote the big cone, interior of the pesudo-effective cone, by $\psf(X)^o$.  In particular ample classes lie in
 $\psf(X)^o$.
 
 The following corollary of \propref{2.2.15}, a lower bound on volume near a nef divisor, is  used in \cite{bfj} to prove differentiability of volume on the big cone.
 The point is that the constants in front of the quadratic error term $t^2$ do not change upon blowing up $X$ or diminishing $\beta$.
 
 \begin{cor}  \label{15cor}  \cite{bfj} Let $X$ be a   projective $n$-dimensional variety.   Let $\beta \in N^1(X)$ be nef, and $\gamma \in N^1(X)$.  Let $\omega$ be nef and big,
 and assume $\beta \leq \omega$ and $\omega \pm \gamma $ is nef.  Then for any $t \in [-1,1]$ we have
 \[ \vol(\beta+t \gamma) \geq \beta^n + n t \beta^{n-1} \cdot \gamma - 8^n \omega^n t^2 \]
  \end{cor}
 
\prf     As the inequality for $t$ and $\gamma$ is the same as for $-t$ and $-\gamma$,
we may assume $t \in [0,1]$.  
Note first that with $\gamma'=\omega-\gamma$, we have $\gamma' \leq 2 \omega$, so 
\[ \beta^i \gamma^{n-i} = \beta^i (\omega -  \gamma')^{n-i} \geq    \beta^i (\omega-2 \omega)^{n-i} \geq - \beta^i \omega^{n-i} 
\geq - \omega^n \]
   Thus (using $|t|\leq 1$)
\[ (\beta+t \gamma)^n - \beta^n -   n t \beta^{n-1} \cdot \gamma =  \sum_{i \geq 2} {{n}\choose{i}} t^i \beta^i \gamma^{n-i}
\geq  - \omega^n t^2 \sum_{i \geq 2}  {{n}\choose{i}}  \geq - \omega^n 2^n t^2 \]
Now let $A=\beta+t(\gamma + \omega)$ and $B=t \omega$.  Then $A \leq 3 \omega$ so 
\[ (A-B)^n = A^n - n A^{n-1}B + \sum_{i\geq 2} \pm {{n}\choose{i}}  t^i A^{n-i}  \omega^i \leq A^n - n A^{n-1}B  +t^2(3 \omega)^n   2^n \]
But $(A-B) = \beta + t \gamma$, so combining the last two displayed equations,
\[ A^n - n A^{n-1}B \geq \beta^n + n t \beta^{n-1}\gamma - 8^n \omega^n t^2 \] 
By \propref{2.2.15} and \remref{15nef}, we have 
\[ \vol(\beta+t\gamma) = \vol(A-B) \geq \beta^n + n t \beta^{n-1} \gamma - 8^n \omega^n t^2 \] 
\eprf

We include here a basic lemma on the effect of morphisms on the pseudo-effective cone of divisors, and its dual.
%

 \begin{lem}\label{functorialcone}  Let $\pi: X \to U$ be a surjective morphism of projective varieties, $U$ normal.  Then 
    \begin{enumerate}
    \item 
    $(\pi^*) \inv \psf{X} = \psf{U}$. 
    \item    If $\dim(X) > \dim(U)$, then $\pi^*(\psf(U))$ lies on the boundary of $\psf{X}$; and only then.
    \item 
    $ \pi_* N_1^+(X) = N_1^+(U)$,
    and $ \pi_* N_1^+(X)^o= N_1^+(U)^o$.  
    \item     $ \pi_* N_1^{eff}(X) = N_1^{eff}(U)$,
    and $ \pi_* N_1^{eff}(X)^o= N_1^{eff} (U)^o$.    
    \item  If $U$ is normal and $D$ is a Cartier divisor on $U$, then $\Oo(X,\pi^* D)  \cong \Oo(U,D)$.
     \end{enumerate}
    \end{lem}
\prf  
We first prove (4).  If $C$ is an irreducible curve on $U$, then it is the image under $\pi$ of an irreducible curve $C'$ on $X $,
 and we have $[C]= d \pi_* [C']$ for an appropriate $d >0$.   It follows that
 $\pi_* N_1^{eff}(X)=N_1^{eff}(U)$.  In particular
  $\pi_* : N_1(X) \to N_1(U)$ is surjective, and hence (being linear) open.  In particular  $\pi_* N_1^{eff}(X) ^o \subset N_1^{eff}(U)^o$.  By continuity $(\pi_*) \inv N_1^{eff}(U)^o \subset \pi_* N_1^{eff}(X) ^o$ and so
  $\pi_* N_1^{eff}(X) ^o = N_1^{eff}(U)^o$. 

As for (5), consider a rational function $h$ on $X$ which is regular on $X \m \pi^*D$.  Let $h$ be a section of $\pi^*D$.
Then $h$ has no zeroes or poles on a fiber of $X \to U$ above $U \m D$,
so it must be constant on the fibers of $\pi$ above $U \m supp(D)$.   So $h$ equals a rational function on $U$. 
Moreover $\pi^*(D + (h) ) \geq 0 $ so $D+(h) \geq 0$ - any poles of $D+(h)$ as a Weil divisor would show up
in the pullback.

Towards (1), consider  Cartier divisors $D$ on $U$ (with integer coefficients). If $D$ is effective, it is clear that $\pi^*D$ is effective.  Going to rational coefficients and by continuity of $\pi$, it follows that   $(\pi^*) \inv \psf{X} \subseteq \psf{U}$.  
 
If $D \in N^1(U) \m \psf{U}$, then $D$ is represented say as $\sum \alpha_i D_i$ with   $\alpha_i \in \Rr$ 
and $D_i$ a Cartier divisor, and (by \corref{bdpp}) there exists a curve $C$ on $U$ with proper intersection with each $D_i$,
and with $\sum \alpha_i D_i \cdot C < 0$.  By the proof of (4) we may take $C=\pi_* C'$ for some curve
$C'$ on $X$.  It follows that $C'$ is not contained in any $D'_i= \pi^* D_i$, and that 
$\sum \alpha_i D'_i \cdot C'  = \sum \alpha_i D_i \cdot C_i < 0$.  Thus $\pi^* D \notin \psf{X}$.  This finishes
the proof of (1).  

A more elementary proof of (1) can be obtained as follows.   We may find a common generically finite cover
of $X$ and of $U \times \Pp^m$, for some $m$.  As the pullback to $U \times \Pp^m$ preserves non-effectiveness,
we may assume we have a generically finite cover.  

Let $D \in N^1(U) \m \psf(U)$.  Then $D$ is a positive
real multiple of some Cartier divisors, and we may take $D$ to be a Cartier divisor, no integer power of which is effective.   
  Now if a Cartier divisor $D$ on a normal variety $U$ is {\em not} effective, the pullback is not effective either.
Indeed on some affine neighborhood $Spec(A)$, $D$ is represented by $(f)$ with $f \notin A$. 
By normality, $f$ is not integral over $A$ so there exists a valuation $v$ of the function field $K$, non-negative on
$A$ but with $v(f)<0$.  The valuation $v$ extends to the function field of $X$, and has nonempty
center on $X$ since  $X$ is projective.  Around any point of this center, $f$ cannot be
a regular function; so the pullback of $D$ is not effective.    This persists with $\Qq$-coefficients.  
 So the pullback of the effective cone of divisors on  $X$ is the effective cone of divisors on $U$.

For (2), suppose $\pi^*(\psf(U))$ meets the interior of $\psf{X}$.  Then by density
there exists a $\Qq$-Cartier divisor $D$ with $\pi^*(D) $ big, and multiplying by an integer we may
assume $D$ is a Cartier divisor.  But by (5),  
$\Oo(X,\pi^* D)  \cong \Oo(U,D)$.  Thus $\pi^*D^n$ has a dimension $O(n^{\dim(U)})$,
and cannnot be big if $\dim(U)<\dim(X)$.  Conversely, if $\dim(X)=\dim(U)$ then the pullback of an ample divisor
of $U$ is big on $X$, and thus lies in the interior of $\psf{X}$ and not on the boundary.  
This proves (2).  
 
For (3),  Let $c \in N_1^+(X)$.  Then for any $d \in \psf(U)$, we have $\pi^*(d) \in \psf(X)$,
so $\pi_*c \cdot d = \pi_*( c \cdot \pi^*(d)) \geq 0$.  Thus $ \pi_* N_1^+(X) \subseteq N_1^+(U)$. 
By the M. Riesz theorem, since $(\pi^*) \inv \psf{X} = \psf{U}$, any positive linear map on $\psf{U}$
extends to one on $\psf{X}$; in other words, $\pi_* N_1^+(X)  =N_1^+(U)$.  The rest of (3)  follows from
 \lemref{coneproject} and \lemref{openconvex}.

\eprf

 \end{section}

\begin{section}{Okounkov bodies}
 
\label{okounkov-s}
 
\ssec{Valuations of maximal rank}  Let $K/k$ be a field extension with $tr.deg._k K = d$.   A 
 {\em maximal rank} valuation $v$ on $K/k$ is a valuation on $K/k$ with   value group $\Lam=v(K^*)$ of rank $d$;
 i.e. $E_v:= \Rr \tensor \Lam$ is a $d$-dimensional real vector space.   Later, we will use the Lebesgue (or Haar)
 measure on $E_v$, normalized so that $\Lam$ has covolume $1$.
 
When $K/k$ is finitely generated, a maximal rank valuation necessarily has finitely generated value group, isomorphic to $\Zz^d$;
though we will not need this fact here.  

We note  that any field extension $K/k$ of finite transcendence degree $d$ admits a maximal rank valuation:
we can embed $K$ into the algebraic closure
of $k[X_1,\ldots,X_d]$, and extend the lexicographic valuation $v(\sum \a_\nu X^\nu) = \min \{\nu: \a_\nu \neq 0 \}$.  
This gives a valuation on $K$ with value group a finite extension of $\Zz^d$, so it must also be isomorphic to $\Zz^d$.

\begin{lem} \label{dim-c} Let $v$ be a  valuation on $K/k$.  Let $H$ be a finite-dimensional $k$-subspace of $K$, of dimension $e$.  
Then $v(H \m (0)) $ is a set of cardinality at most $e$.    If $v$ has residue field $k$, in particular if $k=k^{alg}$ and $v$ has maximal rank, then $|v(H \m (0))| = e$.
\end{lem}

\prf  For $\g \in v(H \m (0))$, let $H_\g = \{h \in H:  v(h) \geq \g \}$.  Then $H_\g$ is a nonzero subspace of $H$; if
$\g' < \g \in v(H \m (0))$ then $H_{\g'}$ contains $H_\g$ properly.  Thus $\g \mapsto \dim H_\g$ is an injective map
from $v(H \m (0))$ to $\{1,\ldots,e\}$.

Now assume $v$ has residue field $k$.  Pick $\g \in \G$, and let $\g'$ be the least element of $v(H)$ greater than $\g$.
Pick $h_\g \in H$ with $v(h_\g)=\g$.  Define
$f: H_\g \to k$ by $f(x) = \res(x/h_\g)$.  Then $f$ is $k$-linear, with kernel $H_{\g'}$.  This shows that $H_\g/H_{\g'}$ is 
$1$-dimensional.  So $e=\dim(H)$. 
\eprf

\ssec{Okounkov bodies ( \cite{okounkov}  )}  

To define Okounkov bodies, we will need the following data:  a function field $K/k$, i.e. a finitely generated extension field of $k$, of
transcendence degree $d$; a   maximal rank  valuation $v$ on $K/k$, and a   $k$-subalgebra $S$ of $K$ generating $K$ as a field, along with a filtration $S=\union_n S_n$, with $1 \in S_1$ and $\dim(S_n) =O( n^d)$ (i.e. for some $C \in \Rr$, for $n \geq 1$ we have $\dim(S_n) \leq C n^d$.)

We will view $v$ as  auxiliary and omit it from the notation; we thus write $Ok(S)$.  

Typically, $S_n$ will be $L(X,nD) \cong H^0(X,nD)$, where $X$ is a  $k$- variety, $D$ an effective divisor on $X$
(so $1 \in S_1$).  In this case we write $Ok(X,D)$ for the Okounkov body.

 \begin{defn}  Let $V_n=\{ v(x):  x \in S_n \m (0) \}$,  $V=\union_n V_n$.  
 $Ok(S)$ is the closed convex hull in $E_v$ of $\union_n \frac{V_n}{n}$.     We will also write $Ok(V)$ for the same set.  \end{defn}

 We could also define $Ok(S)$ to be the closure of $\union_n \frac{V_n}{n}$, since this set is closed under the operations
 $(u,v) \mapsto \alpha u + (1-\alpha ) v$ for rational $\alpha \in (0,1)$.  
 
 
Since $S$ generates $K$, there exist $x_1,\ldots,x_n$ with $v(x_i)$
linearly independent;  the convex hull $[0,v(x_1),\cdots,v(x_n)]$ of $0=v(1),v(x_1),\ldots,v(x_n)$ is an $n$-simplex, so $Ok(S)$  has nonempty interior $S^o$.  

We let $\mu_n$ be the counting measure on $\frac{V_n}{n}$,
normalized so that each point has measure $n^{-d}$.  This can be compared to $\mu'_n$, the (infinite) counting measure
on all of $\Lambda/n$, with the same normalization.

If $S'$ is a $k$-subalgebra of $S$, $S'_n:= S_n \meet S'$, then $V'_n \subset V_n$ and $Ok(S') \subset Ok(S)$.

 \claim{}  \label{ok-main} Let $C \subset Ok(S)^o$ be compact.  Then
 there exists a finitely generated subalgebra $S'$ of $S$,
with $C \subset Ok(S')^o$; in fact there exists a finitely generated subsemigroup $V'$ of $V=v(S \m (0))$
with $C \subset Ok(V')^0$.  
  Moreover,  for large enough $n$, 
\[ C \meet \frac{\Lam}{n} \subset V'_n \]

\prf  Immediate from \propref{sgroups2}.  \eprf

\begin{thm} \label{ok} The measures  $\mu_n$ converge weakly to Lebesgue measure $\mu$ on $Ok(S)$.  We have
$\mu(Ok(S)) = \vol(D) / \dim(X)! $.  
\end{thm}

\prf
 Let $\mu'_n $ be the counting measure on $\Lambda / n$, normalized so that each point has measure $n^{-d}$.  
Then $\mu_n \leq \mu'_n$, and $\mu'_n$ tends weakly to $\mu$ on any compact set.  This shows that any limit
of the $\mu_n$ is bounded below $\mu$.

Note that $\mu(\bd Ok(S))=0$ (\lemref{convex-measure}).  Thus it suffices to show that $\mu_n \to \mu$ on any compact
subset $C$ of the interior $Ok(S)^o$.  

But on $C$, for large $n$, by \lemref{ok-main}, $\mu'_n = \mu_n$.   As $\mu'_n$ tends weakly to $\mu$, so  does $\mu_n$.  

Let $v_n=\mu_n(Ok(S)) $; so $v_n \to \mu(Ok(S))$; and we have by definition $\vol(D) = (\dim(X)!) \limsup v_n  $.

\eprf

\begin{cor}  $Ok(S)$ is compact.  \end{cor} 
\prf By \thmref{ok}, and since $\mu_n$ has bounded total mass, $Ok(S)$ has finite Lebesgue measure.  Compactness follows from \lemref{convex-measure}.\eprf

\begin{cor}[bigness, volume properties]  Let $D$ be an integral divisor on $X$ such that $tr. deg. (\oplus _n H^0(X,nD)) = \dim(X)$.
Then $\vol(D)>0$, and this  volume is attained as a limit, not just a lim sup.  \end{cor} 

\prf As the $\mu_n$ tend weakly to a limit, namely to Lebesgue measure on $Ok(S)$, the total volume 
$v_n=\mu_n(Ok(S)) $ tends to a limit $\alpha>0$; we have by definition $\vol(D) = (\dim(X)!) \limsup v_n = (\dim(X)!) \lim v_n$.
\eprf

Here is Okounkov's proof of log-concavity of volume, in the setting of big divisors.   
 (Another proof that will be accessible a little later  reduces to 
 the ample case \lemref{logconcave1} by Fujita approximation; but our proof of Fujita approximation will
 also use Okounkov bodies.)   As shown in \cite{bfj}, this latter proof lends itself to
 an analysis of the case of equality (sharp log-concavity) for nef big divisors.  It would be interesting 
 to determine when two divisors have the same Okounkov bodies, or translates thereof, and thus 
 deduce the sharp log-concavity  using Okounkov bodies. 

\begin{cor}[log concavity of volume] \label{logconcavevolume} $\vol^{1/n}$ is concave:  for  $D_1,D_2$ with $\vol^{1/n}(D_i) >0$ we have $\vol^{1/n}(D_1+D_2) \geq \vol^{1/n}(D_1) + \vol^{1/n}(D_2) $ 
\end{cor}

\prf  $Ok(X,D) $ contains the Minkowski sum $Ok(X,D_1)+Ok(X,D_2)$.   Denoting these bodies by $O,O_1,O_2$ we have by Brunn-Minkowski (\cite{gromov} 3.1 ):
\[ \mu(O)^{1/n} \geq \mu(O_1)^{1/n}  + \mu(O_2)^{1/n}  \]
where $\mu$ is  Lebesgue measure, normalized so that the   group $\G$ generate by $v(\oplus_m L(X,mD))$ has covolume $1$.  In this case the corresponding groups $\G_i$ of $D_i$  have covolume $\geq 1$, so if
we let $\mu_i$ be the normalization giving $O_i$ covolume $1$ we have $\mu_i = c_i \mu$ with $c_i<1$.  Thus
\[ \mu(O)^{1/n} \geq \mu_1(O_1)^{1/n}  + \mu_2(O_2)^{1/n}  \]
and we are done by \thmref{ok}.   \eprf

\end{section}

\begin{section}{Fujita approximation}

\label{fujita-s}

 \ssec{ Global generation}

   \begin{defn}   
   
   A Cartier divisor $D$ is {\em globally generated} or {\em base-point free} or simply {\em free}
    if for any point $p$ there exists a global section that does not vanish at $p$.  
 
  A $\Qq$-divisor $D$ is said to be 
semi-free if $mD$ is globally generated, for some $m \in \Nn$.  \end{defn}

 Since a one-dimensional vector space is generated by any subset containing
 a nonzero element, a Cartier divisor $D$ is free iff    associated invertible sheaf is
 generated by its global sections.  

\begin{rem}\label{freenef}
If $D$ is semi-free, then $D$ is nef, i.e.  $D^k \cdot Y \geq 0$
for any $k$-dimensional subvariety $Y$.   In particular, $D \cdot Y \geq 0 $ for every curve $Y$. 
   \rm To see this we may assume $D$ is  free.
Given a variety $Y$ on $X$, we can find a global section $s$ of $D$ that 
does not vanish on $Y$; then $D \cdot Y$ is represented by the zero scheme of $s$ on $Y$. so it is a positive sum
of effective $[Y']$ of dimension $k-1$.  %
\end{rem}

 \begin{rems} \label{freerems}\begin{enumerate}  Let $D$ be a free Cartier divisor on $X$.  
 \item  Any pullback $g^*D$ via $g: X' \to X$ is free on $X'$.  
 \item  If $D,D'$ represent the same line bundle, then $D'$ is free.   (Global generation is a property of the line bundle.)
 \item   
 We have a natural morphism $f:X \to PH^0(X,D)$ (Grothendieck convention).  We have $D = f^* \Oo[1]$
 (the global sections agree by definition, and both line bundles are generated by them).
\item  $D$ is ample iff $D \cdot C >0$ for any curve $C$ on $X$.   Indeed,  $D \cdot C = f^* \Oo[1] \cdot C = \Oo[1] \cdot f_*C$ so 
if $D \cdot C  > 0$ then 
$f_*C$ cannot be a point.  Thus $f$ does not collapse any curve to a point, so it is finite-to-one.  By a theorem
of Serre, $f^* \Oo[1]$ is ample. 
\item  $D$ has at least one section, so it is effective.
\item  If $D'$ is an ample Cartier divisor, then $D+D'$ is ample.  (This follows immediately from (4), aplied to $D+D'$.)
\end{enumerate}
   \end{rems}
   
   \begin{lem}[ Corollary 1.32 of \cite{debarre}]  \label{freeample} Assume $D$ is a nef Cartier divisor
   with $\vol(D)>0$.   Then there exist ample Cartier divisors
   $D_t$ (say for $t$ a small positive rational) such that $D_t \leq D$ and $D_t \to D$ numerically.  In particular, 
for any $\e>0$   there exists an ample $A \leq D$ with $\vol(A) \geq (1-\e) \vol(D)$.  \end{lem}

 \prf  Let $H$ be any ample divisor, represented by a hypersurface $E$.   As $\dim(E) =n-1$, by \corref{vol-ample} we have
  $h^0(m D|E) = O(m^{n-1})$, while by definition $h^0(mD) $ grows like $m^n$.  The exact sequence 
  \[ 0 \to L(X, mD-H) \to L(X,mD) \to   L(E, mD | E) \]
shows that  for large enough $m$,  $L(X,mD-H) \neq 0$. 
     Thus $mD = H + E'$ with $E' \geq 0$.
For any $t \in (0,1)$, let $D_t = \frac{t}{m} H + (1-t) D$.   
But $D_t$ is ample (the boundary of the ample cone stabilizers the ample cone).   Now $D= D_t +    \frac{t}{m} E'$
so $D_t \leq D$, and $D_t \to D$ as $t \to 0$;  by continuity of volume at $D$, 
 $\vol(D_t) \to \vol(D)$.  
  \eprf

  \begin{rem}  This gives a quick way of seeing using \defref{nef} that the pullback of the nef cone is contained in the nef cone:  the pullback of an ample is free,
 of positive volume, and so can be approximated by ample divisors.  
 \end{rem}

 \begin{lem} \label{amalg-gg}    (cf. \cite{bfj} Lemma 2.6). Let $C,D$ be Cartier divisors on $X$, and assume 
 their stable join $C \vee D$ exists.  If $C,D$   are free, then so is  $C \vee D$.
  
 \end{lem}
 
 \prf  By assumption, $C,D$ are generated by global sections; we have to show the same for $C \vee D$.
 In other words given a point $p \in X$, we have to find a global section of $C \vee D$ that does not
 vanish at $p$.  Let $f$ be a global section of $C$ not vanishing at $p$, and $g$ a  global section of $D$ not vanishing at $p$.  On some neighborhood $U$ of $p$, $C$ is represented by $(f \inv)$ and $D$
 by $(g \inv)$, while $C \vee D$ is represented by $h \inv$ say.   Then $-C$ is represented by $(f)$,
 $-D$ by $(g)$, and $-(C \vee D) = (-C) \wedge (-D)$ by $(h)$; so we have $(h) = (f,g)$.  Thus
 $f=rh$, $g=sh$ with $r,s$ regular at $p$. Not both $r,s$ can be $0$, since otherwise $(h) = \Mm_p (h)$
 where $\Mm_p$ is the maximal ideal of the local ring $\Oo_p$, contradicting Nakayama.   Now 
 both $f,g$ are global sections of $C \vee D$, and at least one of them fails to vanish at $p$.    
 \eprf

  \begin{lem}[\cite{bfj} Lemma 2.6]   \label{amalg-ggnef}  $D_1,D_2$ be $\Qq$-Cartier divisors, and assume  $D=D_1 \vee D_2$ 
  is the stable join of $D_1,D_2$.  If $D_1,D_2$ are semi-free, then so is $D$.  If $D_1,D_2$ are nef, then so is $D$.  
 \end{lem}
 
 \prf  For semi-free this is immediate from \lemref{amalg-gg}.  Assume $D_1,D_2$ are nef.  Let $A$ be a very ample
 Cartier divisor.  
   Then for any rational $t>0$, $D_i+tA$ is ample, hence semi-free.   Moreover $D+tA$ is the 
   free join of $D_1+tA$ and $D_2+tA$.   Thus $D+tA$ is semi-free, hence nef.  Letting $t \to 0$ we see
   that $D$ is nef.   
 \eprf
  
  \begin{lem} \label{fujita-lem} Let $D$ be an effective Cartier divisor on $X$, and let $f_1,\ldots,f_n \in L(D)$ be sections of $D$.
  Then there exists a blow-up $\pi: X' \to X$ and a free divisor $D'$ on $X'$ with $D' \leq \pi ^*D$,
    such that each $f_i \circ \pi $ is a section of $D'$.  \end{lem}

\prf  First we may find a  birational morphism  $\pi^*: X^* \to X$, such that each rational function $f_i$ extends to a morphism $f_i^*: X^* \to \Pp^1$.      Let $N_i =   (f_i^* )\inv(\infty) $.  Clearly $f_i \circ \pi^*$ is a section of $N_i$.
The Cartier divisor $N_i$ is represented by $f_i \inv$
away from $f_i \inv(0)$, and by $1$ away from $f_i \inv(\infty)$; as $D \geq 0$ and $(f_i)+D \geq 0$, 
we have $N_i \leq \pi^* D$.    Being the pullback of the free divisor $\infty$
of $\Pp^1$,   $N_i$ is free.

Passing to a further blowup $X'$, we may assume the $N_i$ have a stable join $N$ (\secref{stablemeet}).  Each $N_i \leq   \pi^* D$, so $N \leq  \pi^* D$.   Being the stable join 
of free  divisors, $N$ is also free, by \lemref{amalg-gg}.  
\eprf

\thmref{fujita1}  was   proved   in \cite{fujita} in characteristic zero, developing, according to \cite{lazarsfeld},
ideas of Demailly on positive intersections.    In positive characteristic it is due to \cite{cutkosky}.  
   
\begin{thm}[Fujita approximation]  \label{fujita1}   
  Let $D=\sum \alpha_i D_i$   
  with $D_i$ an effective Cartier divisor, and $\alpha_i \in \Rr, \alpha_i >0$.  Assume $\vol(D)>0$, and let 
 and let $\e>0$.
Then there exists a birational morphism $\phi: X'  \to X$, $D' = \phi^* D$, with $D' \geq N$ for some
ample  $\Qq$-divisor $N$,
satisfying $\vol(N) \geq (1-\e) \vol(D)$.  
   \end{thm}
  
\prf  Note that there exists a very ample $A$ with $D \leq A$; it follows that $\dim L(X,kD) \leq \dim L(X,kA) \leq O(k^d)$, $d =\dim(X)$.

 Let $K=k(X)$, and $S_k = L(X,kD)$.  Let 
$C$ be a compact subset of $Ok(X,D)^o=Ok(S)^o$ with Lebesgue measure $\mu(C) > (1-\e/2) \mu(Ok(X,D))$.
Let $f_1 \in L(X,k_1D),\ldots,f_r \in L(X,k_rD)$ be elements, such that the multiplicative semigroup $S'$
they generate in $\oplus_k L(X,kD)$ has Okounkov body containing $C$ (see Claim before \thmref{ok}). 

  If we replace
each $f_i$ by $f_i^m \in L(X,mk_iD)$, the Okounkov body does not change; hence it does not change if we just change
$f_1$ by $f_1^{m_i}$.  In this way, taking  $m_i = \Pi_{j \neq i} k_j$, we may take all $k_i$ to 
be equal to the same integer $k$. 

 By \lemref{fujita-lem},   there exists a blow-up $\pi: X' \to X$ and a free divisor $D'$ on $X'$ with $D' \leq \pi ^*(kD)$,
    such that each $f_i \circ \pi $ is a section of $D'$.  Let $N=D'/k$; then $N$ is semi-free, and  $N \leq \pi^*(D)$.

Each $f_i \in L(X,kN)$.  So $S' \subset \oplus L(X,N)$, hence $C \subset Ok(X,N)$,
and thus $\vol(N) \geq (1-\e/2) \vol(D)$.    

As $\vol(D)>0$,  `semi-free' can be replaced by `ample' using \lemref{freeample}.   
\eprf

This shows in particular that volume is a numerical invariant.

\section{Positive intersection product and differentiability of volume} 

\label{bfj-s}

 The Zariski-Riemann approach of \cite{bfj} fits very well indeed 
with the model theory of GVF's (indeed the description of quantifier-free types led us to the same formalism).
Nevertheless we wish  to remain for a while longer in the finite-dimensional setting, and will describe
the main results of \cite{bfj} without this formalism.  The positive intersection product here is the restriction
of the one in \cite{bfj} to  domain $N^1(X)$ and range $N_1(X)$; moreover we consider only the $n-1$'st power, where $n=\dim(X)$.

\ssec{Definition of $\psi$, the $n-1$'st positive intersection power}

For $x \in \psf(X)^o$, define:

\[ A^{nef}(x) = \{(f,X',y):  f: X' \to X \hbox { birational}, y \in N^1(X') \hbox{ a nef $\Qq$-Cartier divisor },    y \leq f^*(x) \}\]

\[ W_x =  \{ f_* y^{n-1}:   (f,X',y) \in A^{nef}(x) \} \]

 $W_x$ is precompact.  To see this pick an ample $a \geq x$.   
Also pick   ample divisors $b_j$ such that $x \mapsto (x \cdot b_1,\ldots,x \cdot b_m)$ is an isomorphism  $N_1(X) \to \Rr^m$.  Any nef $y \leq f^*(x)$ satisfies $y \leq f^*(a)$, $y^{n-1} \leq  f^*(a)^{n-1}$, and $f_* y^{n-1} \cdot b_i \leq a^{n-1} \cdot b_i$.   (We used \lemref{monotone}.)   Thus the image of $W_x$ is bounded.

Secondly,     $W_x$ is {\em directed} 
 with respect to the partial ordering of $N_1(X)$ given by the closed cone generated
by classes of curves on $X$.     In fact $A^{nef}(x)$ is directed, in the following sense:   if $y_1,y_2$ are nef $\Qq$-Cartier divisors below $f^*(x)$, 
let $D'$ represent $f^*(x)$; and represent $y_1,y_2$ by $\Qq$-Cartier divisors below $D'$.  
By pulling back
to a further blowup we may assume their stable join $y$ exists.  Then $y \leq D'$, and
by \lemref{amalg-ggnef}, $y$ is nef.

 Thus we can define $\psi: \psf(X)^o \to N_1(X)$ by:   
   \[ \psi(x) = \sup W_x \]

 \begin{rems} \begin{enumerate}
 \item    The definition of $\psi$ would not change if in the definition of    $W_x$ we replaced {\em nef}
 by {\em semi-free} or by {\em ample}.  Indeed, since $x$ has positive volume, by \thmref{fujita1} some $y_0 \in A^{nef}(x)$
 has positive volume.   Whenever $(f,X',y) \in A^{nef}(x)$ lies above $y_0$, it has positive volume, and so can be
 approximated by an ample divisor, as any big, nef divisor is approximated by smaller amples
 by  \lemref{freeample} .
\end{enumerate}\end{rems}

   \ 
 \ssec{A second definition} \label{naive}
 
 We also give a definition that does not require blowing up, and instead uses a naive intersection theory on $X$; compare \cite{lazarsfeld}, 11.4.10, 11.4.1,
 where the $n$-fold intersection is treated. 
 
 We have not yet defined $N^k(X)$ for $k$ other than $1,n-1$.  For our purposes we could easily do 
 so by numerical equivalence {\em against $k$ divisors} (a definition that differs from the classical one when $X$ is smooth, but does not really lose information when in the construction below).  At any rate, we will use this material    only
 for $k=\dim(X)-1$.

   Let $D$ be an effective Cartier divisor, with at least one nonzero section.   
 Let $k$ be an integer.  Let $X'$ be a Zariski open subset of $X$. 
 Let $f_1,\ldots,f_{k}$ be $k$ mutually generic sections of $L(X,D)$,   let $Z$ be their 
 common zero set {\em on $X'$}, viewed as a subscheme of $X'$; let $[Z]$ be the sum of components
 of $Z$ of codimension $k$, counted with multiplicity.  Each irreducible subvariety of $X'$ appearing in $[Z]$
 can be identified with its Zariski closure in $X$, so we can also view $[Z]$ as a cycle in $N^k(X)$;
 denote $\nint{D;X'}^k = [Z]$.  
 
 \claim{}   For some $z \in N^k(X)$, 
 for all Zariski open $X'  $ disjoint from the support of $D$, we have $\nint{D;X'}^k = z$.
 
 \prf
Let $B_D$  be the base locus of $D$:  the set of points at which every section of 
 $D$ vanishes.   Let $W$ be a subvariety of $X$.  If $W$ is not contained in $B_{D}$, then   a generic section $f$ of $L(X,D)$
 does not vanish on $W$, so cutting $W$ with the zero set of $f$ reduces dimension.   It follows
 that as soon as $X'  \meet B_D = \emptyset$, all components of $Z$ have codimension at least $k$.
 
 Similarly if we replace $X'$ by $X''=X' \m W$, where $W$ is some hypersurface of $X$, then for generically chosen
 $f_1,\ldots,f_k$,  the intersection of $W$ with the zero sets of the $f_i$ has codimension at least $k+1$; no $k$-dimensional
 component of $Z$ will be contained in $W$; thus $\nint{D;X'} ^k = \nint{D;X''} ^k$.   \eprf
 
 We can now define $\nint{D}^k_{naive} $ to be equal to $\nint{D;X'} ^k$, for any sufficiently small open $X' \subset X$.  We will assume in particular (the base field being perfect) that $X'$ is smooth, i.e. disjoint from the singular locus
 of $X$.  
 
 Note also that if $g_1,\ldots,g_k \in L(X,D)$ have common zero set $Z'$ on $X'$ of codimension $k$,   
 then $[Z'] \leq [Z] = \nint{D;X'}^k$ (the inequality holds with respect to the effective cone in $N^k(X)$).  In other words the generic choice of sections maximizes the  
 intersection on $X'$.  The reason is that $(f_1,\ldots,f_k)$ specializes to $(g_1,\ldots,g_k)$;
 by smoothness of $X'$, every component of $Z'$ is the specialization of some component of $Z$,
 and is covered with appropriate multiplicity. 
  
  Some components of $Z$ may specialize into $X \m X'$
 and not be counted in $Z'$; but at any rate we have the inequality.  We may summarize
 this symbolically in the formula:
 \[ \nint{D}^k_{naive} = \inf_{X'} \sup_{f_1,\ldots,f_k \in L(D)} [Z_{X'}(f_1) \meet \cdots \meet Z_{X'}(f_k)]_{n-k} \]
 

 If $D \leq D'$, let us compare 
   $\nint{D}_{naive}^k $ with $ \nint{D'}_{naive}^k$.  We have $L(X,D) \leq L(X,D')$  Choose a sufficiently small open $X'$,
 and  generic $f_1,\ldots,f_k \in L(X,D)$ and $g_1,\ldots,g_k \in L(X,D')$.  Then
 by the above specialization argument, the generic choice in $D'$ gives a bigger or equal result.

 Thus $\nint{D;X'}_{naive}^k \leq \nint{D';X'}_{naive}^k $, and 
 \[ \nint{D}_{naive}^k \leq \nint{D'}_{naive}^k \]

 Next we compare $\nint{D}_{naive}^k$ to $\nint{l D}_{naive}^k$.   
 Again choose a sufficiently small open $X'$,
 and  generic $f_1,\ldots,f_k \in L(X,D)$ and $g_1,\ldots,g_k \in L(X,lD)$.   Then $(g_1,\ldots,g_{k})$ specializes to $(f_1^l,\ldots,f_{k}^l)$.  By the linearity  of intersection with Cartier divisors, $\meet Z_{X'}(f_i^l) = l^k \meet Z_{X'} (f_i)$.  Hence
 it follows as above that $l^n \nint{D;X'}_{naive}^k \leq \nint{l D; X'}_{naive}^k$ and  so 
 \[ l^k \nint{D}_{naive}^k \leq \nint{l D}_{naive}^k \] 

Finally we define
\[ \nint{D}^k = \limsup_m \frac{1}{m^k}   \nint{m D} _{naive}^k \]

(in fact it follows easily from the above that this is also the $\liminf$ and hence the limit.)

Clearly $\nint{l D}^k = l^k \nint{D}$.  Observe, directly from this definition:  
\begin{prop} \label{sectionsonly}   $\nint{ D}^k$ depends on $D$ only via the section spaces $(L(mD):m \in \Nn)$.  
\end{prop}
  
 In the following lemma, we show equality of the $1$-cycles $\nint{D}^{n-1} $ and $ \psi(D)$ on $X$.
 For our purposes, this can be understood numerically, i.e. $  \psi(D) \cdot A  =  \nint{D}^{n-1} \cdot A$
 for all divisors $A$, equivalently for all very ample divisors $A$.  The proof will use inequalities 
 with respect to the pseudo-effective cone of curves; equivalently, $u \leq v$ iff $u \cdot A \leq v \cdot A$
 for all very ample $A$, since  the dual of the ample cone of divisors is the pseudo-effective cone of curves.
 
 \begin{lem}   Let $D$ be a big Cartier divisor.  Then  $\nint{D}^{n-1} = \psi(D)$.  \end{lem}
 
 \prf   We first show that $\nint{D}^{n-1} \leq  \psi(D)$.  Since $\psi$ is homogeneous, it 
 suffices to show that  $\nint{D}^{(n-1)}_{naive} \leq \psi(D)$.  (This for all $lD$ implies the claim, by taking limits.)
   Let $X'$ be small enough  to compute $\nint{D}^{n-1}_{naive}$, as in the Claim.
  Let $f_1,\ldots,f_{{(n-1)}}$ be ${(n-1)}$ mutually generic sections of $L(X,D)$.
 By   \lemref{fujita-lem} there exists a blow-up $\pi: X' \to X$ and a free divisor $D'$ on $X'$ with $D' \leq \pi ^*D$,
    such that each $f_i'= f_i \circ \pi $ is a section of $D'$; moreover by the proof, such that the zero locus $(f_i)^+$ of $f_i'$ is nef.  
    So for any very ample $A$
    \[  \nint{D}^{(n-1)}_{naive} \cdot A =  \nint{D;X'}^{(n-1)}  \cdot A \leq (f_1)^+  \cdot \ldots \cdot (f_{(n-1)})^+ \cdot \pi^* A \leq (D')^{n-1}\cdot \pi^*A \leq \psi(D) \cdot A  \]
   (As $(f_i)^+$ and $D'$ are nef, and $(f_i)^+ \leq D'$.)

Conversely  for any blow-up $\pi: {\tX} \to X$ and ample $\Qq$-Cartier divisor $B \leq \pi^*D$, we have to show that
\[ \pi_*(B^{n-1}) \leq \nint{D}^{n-1} =  \limsup_m \frac{1}{m^n}   \nint{m D} _{naive}^{n-1} \]
It suffices to show for some $m$ that  \[\pi_*(B^{n-1})  \leq \frac{1}{m^n}   \nint{m D} _{naive}^{n-1} \] 
Multiplying   $B$ and $D$ by a large integer,   we may assume $B$ is a very ample Cartier divisor (and then prove the last inequality with $m=1$.)

Note that (as $X$ is normal)  $L(X,D) = L(\tX,\pi^*(D) ) \supseteq L(X,B)$.  Let $f_1,\ldots,f_{n-1}$ be  sufficiently generic global sections   of $B$.
Let $Z_i$ be the zeroes of $f_i$, viewed as a divisor on   $\tX$.   Let $e$ be the blowing-up locuas of $\pi$, and $E = \pi^{-1}(e)$;
so $\pi$ restricts to an isomorphism $\tX \m E \to X \m e$.   Then $E$ has
codimension at least $1$ on $\tX$; as $B$ is very ample, $E \meet \meet_{i=1}^{n-1}  Z_i $ is finite.  Hence, $\pi_*(  \meet_{i=1}^{n-1}  Z_i)$ 
as a $1$-cycle has no components contained in $e$.  Viewing the $f_i$ now as elements of $L(X,D)$, denoting their zeroes by $Z_i^+$, we see that
no component of their intersection will be contained in $e$ (as all components are $1$-dimensional), and so $\pi_*(  \meet_{i=1}^{n-1}  Z_i) = \meet_i Z_i^+$.
 As we saw above, the intersection of the zero sets
on $X$ of generic $g_1,\ldots,g_{n-1} \in L(D)$ can onlty be bigger.  This proves the required inequality.  
 \eprf

 \begin{lem} \label{psicont}  $\psi$     
  is homogeneous of degree $n-1$, monotone and continuous on $\psf(X)^o$.
 \end{lem}
 
 \prf      Monotonicity  is immediate from the definition of $\psi$, and homogeneity is also easy.  
 
 Now any   function $v$, monotone with respect to the ordering of a cone $E$, and homogeneous of degree $d \geq 1$, must be continuous
 on the interior of $E$.  For let $e \in E^o$ and $\e>0$.  Then $\e a \in E^0$, so there exists a symmetric neighborhood $N=N_\e$ of $0$
 with $\e a + N \subset E^0$, hence $\e a \geq N$.  Thus $(1+\e) a \geq a+N \geq (1-\e)a$ so by homogeneity $(1+\e)^d f(a) \geq f(a+N) \geq (1-\e)^d f(a)$; hence for any linear functional $L$ with $L(E) \geq 0$ we have
 $(1+\e)^d Lf(a) \geq Lf(a+N) \geq (1-\e)^d Lf(a)$; 
   letting $\e \to 0$ we see that $Lf$ is continuous at $a$.    As such $L$'s span the dual space, $f$ is continuous.    \eprf

 Here is the main theorem of \cite{bfj}.  
  
\begin{thm} \label{bfj1}   The volume function $\vol$  is differentiable on the big cone of $N^1(X)$; and $d \vol = n \psi$.
\end{thm}
Here we use  the identification $N_1 = N^1(X)^*$ given by the intersection pairing.  Explicitly, if 
  $\alpha \in N^1(X)$  is big and $\gamma  \in N^1(X)$ is arbitrary, the statement is that 
\[ \frac{d}{dt}  |_{t=0} \vol(\alpha + t \gamma) = n \psi(\a) \cdot \gamma \]

 \prf  Fix an ample  $\omega$ such that $\alpha \leq \omega$ and $\omega \pm \gamma $ is ample.   Consider
birational morphisms $f: X' \to X$; let $\omega' = f^* \omega$, $\alpha' = f^* \alpha$, $\gamma' = f^* \gamma$.
Note that $\omega'$ and  $\omega' \pm \gamma'$  remain nef.  
Recall \corref{15cor}:  if   $\beta \in N^1(X')$ is  nef,   for any $t \in [-1,1]$ we have
$ \vol(\beta+t \gamma') \geq \beta^n + n t \beta^{n-1} \cdot \gamma' - 8^n \omega^n t^2$.   When $\beta \leq \alpha'$
it follows that $\vol(\alpha'+t \gamma') \geq \beta^n + n t \beta^{n-1} \cdot \gamma' - 8^n \omega^n t^2 $.
Taking the supremum over all $X'$ and all nef $\beta \leq \alpha'$, we find, using \thmref{fujita1} and the definition of $\psi$:
 \[ \vol(\alpha+t\gamma) \geq \vol(\alpha) + n t \psi(\alpha) \cdot \gamma - 8^n \omega^n t^2 \]
\, 
 Applying this to $\alpha^* = \alpha+ t \gamma$ ($t$ small enough so that this is still big), $\gamma^* = -\gamma$,
 so $\alpha= \alpha^* + t \gamma^*$, 
 \[ \vol (\alpha) \geq \vol(\alpha + t \gamma) - n t \psi (\alpha^*) \cdot \gamma - 8^n \omega^n t^2 \]
 As $\psi$ is continuous at $\alpha$, this implies that  $\frac{\vol(\alpha+t\gamma)-\vol(\alpha)}{t}$ is bounded on both sides by functions that approach $n \psi(\alpha)$
 as $t \to 0$.   
\eprf 

%
%
%
 
 Let $N_1^+(X)$ be the closed dual cone to $\psf(X)$;  thus $ c \in N_1^+(X)$ iff for every effective $e \in N^1(X)$
 we have $c \cdot e \geq 0$.

 \begin{defn}\label{movable-def} 
A class $c \in N_1(X)$ is called {\em movable} if for any hypersurface $H \subset X$, there exists
an irreducible curve $C$ not lying on $H$, and some $\alpha >0$ with $c = \alpha [C]$.  
  Call a class in $N_1(X)$  `strongly movable' if
it is a pushforward,   of   the intersection product of $n-1$ ample rational
 classes in some blowup of $X$; if these are the same ample class, call it `very strongly movable'.  
 In either case, note that   such a class is   a positive scalar multiple of the class of a curve, a complete intersection of very ample divisors.  
 \end{defn}
 
\begin{thm} \label{surj} Let $\dim(X)=n >1$.  Then 
\[ N_1^+(X)^o \subseteq \psi(\psf(X)^o) \subseteq N_1^+(X) \] \end{thm}

\prf  By definition of $\psi$ as a positive intersection product,  any element of $\psi(\psf(X)^o) $ is the supremum of a family of very strongly movable classes,
and in particular the limit of a sequence of such classes.  Now  if $a$ is ample in $X'$ and $f: X' \to X$ is a birational morphism,
then for any effective $e$ on $X$ we have $f_*(a^{n-1}) \cdot e = a^{n-1} \cdot f^*(e) \geq 0$.  Thus any very strongly movable class lies
in $N_1^+(X)$ and so $\psi(\psf(X)^o) \subseteq N_1^+(X)$.

To show the inclusion $N_1^+(X)^o \subseteq \psi(\psf(X)^o)$ we verify the conditions of \thmref{convexity}, with 
$U=\psf(X)^o$ the big cone, $\phi(x)=\vol(x)^{1/n}$, $F=\vol^{-\frac{n-1}{n}} \psi$, $C=N_1^+(X)$.  
The volume function is continuous and non-negative and hence so is  $\phi$.
  The concavity (1) follows from \corref{logconcavevolume} (or from the nef case along with \thmref{fujita1},
by continuity).

The positive homogeneity (2) is clear, and (3) comes from the definition  of the big cone.

Note that $\Rr^{\geq 1} F(U) \subseteq \Rr^{>0} \psi(U) = \psi(U)$ (this last equality, by the $n-1$-positive homogeneity of $\psi$).
By \corref{convexity-c}, $N_1^+(X)^o \subseteq  \psi(U)$.  

\eprf

\begin{cor} \cite{bdpp} \label{bdpp} The dual of the pseudo-effective cone, $N_1^+(X)$, is the closed convex cone generated by the (very) strongly movable classes, or by all movable classes.  \end{cor}

For  our purposes it is important to  know that the image of $\psi$ is dense, without taking convex closure.    

 \begin{cor} \label{eclem} Let $c \in N_1^+(X)$.
 Then $c$ can be approximated by (very strongly) movable classes.
 
 \end{cor}
 
 Indeed $[C]$ can be taken to have the form $f_*(D_1 \meet \cdots \meet D_{n-1})$
 where $D_i$ are Bertini / Weil representatives of an ample class on $X'$, $f: X' \to X $ birational.

 \begin{question}  Is any $\Qq$-Cartier divisor in $N_1^+(X)^o$ movable? \end{question}
 
\begin{rems}  \label{1.13}  \begin{enumerate}
\item   $  {\psi}$ is injective on $Nef^1(X) \meet \psf(X)^o$.  
\item If $x,y \in Nef^1(X) \meet \psf(X)^o$ and ${\psi}(x),{\psi}(y)$ are proportional then so are $x,y$.  
 \item  ${\psi}$ is a homeomorphism   on the ample cone, and continuous  on the nef cone.

\end{enumerate}
 \end{rems}
\prf  (1,2)  $  \vol^{1/n} $ is strictly concave on $Nef^1(X) \meet \psf(X)^o$ by \cite{bfj} Corollary E.  
Hence $ \log \vol $ is also strictly concave there.  By    \lemref{convexity1},  the differential $-n \vol \inv {\psi}$ is therefore strictly increasing,
and hence injective. 

If ${\psi}(x),{\psi}(y)$ are proportional, say ${\psi}(x)=c {\psi}(y)$.  Note $c>0$ since $N_1^+(X)^o$ is a strict cone.  Let $x' = \lambda x$ with  $\lambda = c \vol(y) / \vol(x)$.  
Then $\vol(x') \inv {\psi}(x') = \lambda^{-1} \vol(x) \inv {\psi}(x) = \vol(y) \inv {\psi}(y)$.  By injectivity
of $\vol^{-1} {\psi}$ we have $y = \lambda x$.  In case $c=1$,    by homogeneity of ${\psi}$ we have $\lambda^{n-1}=1$
so $\lambda =1 $, and $x=y$.     

(3) In fact $\psi$ is polynomial on the nef cone, and injective on the ample cone (by (1)). 

\eprf
\end{section}

\begin{section}{Quantifier-free types over a one-dimensional base}

\label{qft-s}

We will work over a constant field $k=k^{alg}$, and describe the quantifier-free GVF types over $k$ geometrically.  We will see later how to deduce a description of types over the non-constant GVF
$k(t)^{alg}$.

When $X=Spec A$ is affine, the Berkovich space of $X$ over $k$   can be defined as the set of all $\Rr$-valued
semi-valuations of $A$; these are just the $\Rr$-valued valuations of the fraction field of $A/P$, for any prime ideal $P$ of $A$.
The  topology on $\hX$ is induced from the Tychonoff topology on the functions from $K  $to $\Rr \union \{\infty\}$.   For general $X$,
the Berkovich space $\hX$ can be defined by glueing over an affine covering.

Let $X$ be an irreducible, normal projective variety over $k$, $K=k(X)$.  Then $\hX$ is compact.

Let $\hK$ be  the space of $\Rr$-valued valuations of $K/k$.   So $\hK= \hX \m \union_{Y} \widehat{Y}$, the union ranging over proper subvarieties of $X$.
 
We define a pairing 
$\beta:  \hX \times Pic(X) \to \Rr_\infty $
as follows:   given $v \in hX$ and a Cartier divisor $D$ pick a Zariski open set $U$ of $X$ such that $\Oo_X(U) \subset \Oo_v$
and $D$ is represented on $U$ by $f \in K$; and let $\beta(v,f) = v(f)$.  It is easy to check that this is well-defined and
linear in the $Pic(X)$- coordinate.  We also have $\beta(v,(f)) = v(f)$ where $(f)$ is the principal Cartier divisor represented
everywhere by $f$.    

This restricts to 

$\beta:  \hK \times Pic(X) \to \Rr $

 Extend it to an $\Rr$-linear map
\[ \beta:  \hK \times Pic_\Rr(X) \to \Rr \]

Let $\Ii_X$ be the set of birational morphisms $\pi: Y \to X$ (over $k$).  Given $\pi,\pi'$ there exists at most one $j_{Y',Y}: Y' \to Y$
with $\pi' = \pi \circ f$; in this case write $\pi \leq \pi'$; we obtain a directed partially ordered set.  Let $\pic (K)$
be the direct limit of $Pic_\Rr(Y)$ over $\Ii_X$.    By going to the limit we obtain:
\[ \beta:  \hK \times \pic_\Rr (K) \to \Rr \]

Finally, let $C(\hK)$ be the space of continuous functions $\hK \to \Rr$,
and   define $\beta_*:   Pic_\Rr(K)  \to C(\hK)$ by ${\beta_*}(D)(v)=\beta(v,D)$.

\begin{lem}\label{1100}  The image   $\beta_*(Pic_\Rr(K))$ separates
points on $\hK$, and is closed under addition, scalar multiplication and $\min$.  It separates points on $\hK$.
\end{lem}

\prf 
Closure under addition and scalar multiplication follows from linearity of $\beta$ in the second variable.  

Let $v \neq v' \in \hX$.   Then $v(f) \neq v'(f)$ for some $f \in K$.   
So  $\beta(v,(f)) \neq \beta(v',(f))$ for the associated principal divisor $(f)$. Thus $\beta_*$ separates points.

Closure under $\min$ follows from the existence of stable meets, in an appropriate
blowup.  By adding a large ample divisor, it suffices to prove this for two effective divisors; we
can assume they come from divisors $W_1,W_2$ on a blowup $Y$ of $X$, with ideal sheaves $I_1,I_2$.
Further by blowing up $I_1+I_2$ we may assume it is principal, and hence determines a Cartier divisor
$W$.  In this case, it is easy to see that $\beta_*(W) = \min \beta_*(W_1),\beta_*(W_2)$, as required.  
  
\eprf

If $Z \subset \hX$ is compact, let $C(Z)$ be the space of continuous functions $Z \to \Rr$,
and   define $\beta_Z $ to be the composition of $\beta_*$ with the restriction function $C(\hK) \to C(Z)$.
We endow $C(Z)$ with the  topology of uniform convergence. 

\begin{lem}  \label{110} Let $Z \subset \hK$ be compact, and assume $v(f_0)=1$ for some $f_0$ and for all $v \in Z$.  
  Then  $\beta_{Z}(Pic_\Rr(K))$ is   dense in $C(Z)$.  
  \end{lem}

\prf  By \lemref{stonew}, it suffices to show that the image of $\beta_Z$ separates
points on $Z$, and is closed under addition, scalar multiplication and $\min$, as well as containing the constant
function $1$.    The constant function $1$ is provided by $v(f_0)$.  The rest  was proved in \lemref{1100}.  
\eprf

\begin{rem}   \rm  We can similarly define $\beta: \hX \times Pic_\Rr(K) \to \Rr_\infty$, and
hence $\beta_* :Pic_\Rr(K) \to C(\hX, \Rr_\infty)$.  Let $Z \subset \hX$ be compact,
and $\beta_Z(y) =\beta_*(y)|Z$.    Let $B=\beta_* (Pic_\Rr(K)) \meet C(Z)$,
be the set of functions of the form $\beta_*(y) |Z$ that take only finite real values, and $\bar{B}$
the  uniform closure of $B$.  Then $\bar{B}=C(Z)$, provided that $1 \in \bar{B}$.  

The proof is the same as of \ref{110}.  For separation of points, if $z \neq z' \in Z$, we find $y \in Pic_\Rr(K)$
as above with $\beta(z,y) \neq \beta(z',y)$, then replace $y$ by $\max(-r,\min(\beta(y,r)))$ where $-r \leq \beta(z,y),\beta(z',y) \leq r$.
\end{rem}

Let $M_K$ be the set of regular Borel measures $\mu$ on $\hK$, such that for any $f \in K$,
$v \mapsto v(f)$ is integrable.  It follows that for any 
$D \in \pic (K)$,
$v \mapsto \beta(v,D)$ is $\mu$ integrable.    
Thus we can naturally extend $\beta$ to a map
\[ \beta:  M_K \times \pic (K) \to \Rr \]

\begin{cor}  \label{11} Let $\mu, \mu' \in M_K$ and suppose they define different GVF structures on $K$, i.e.
$\mu,\mu'$ are distinct even up to renormalization. 
   Then for some $D \in \pic (K)$, $\beta(\mu,D) \neq \beta(\mu',D)$.  \end{cor}
\prf    We may assume $k$ is countable \footnote{For this we use the equivalence of GVF structures and globalizing measures up to renormalization}.   We have  $\hK  \m \{v_{triv}\} = \union_{f \in K \m (0)} W_f$
where $W_f=\{v: v(f)>0\}$.    
Thus   for some $f$,   the measures $\mu,\mu'$ are distinct (even up to normalization) when restricted to $W=\{v: v(f)>0\}$.   By renormalizing we may assume $\mu,\mu'$
concentrate on $ \{v: v(f)=1\}$.   Moreover these are regular measures, approximated by their values on compact sets,
hence $\mu,\mu'$ have distinct restrictions to some compact $Z \subset \{v: v(f)=1\}$.
Now the statement follows from \lemref{110} (with $f_0=f$).
\eprf
 
\begin{prop}  Let $\mu \in M_K$ satisfy the product formula:
\[ \hbox{for all } f \in K \m (0), \ \ \int v(f) d\mu(v) = 0 \]
Then 
 $\beta(\mu, \ ): Pic_\Rr(Y) \to \Rr$ factors through the quotient map $Pic_\Rr(Y) \to NS(Y)$.  
\end{prop}
\prf  We have a 
normal, irreducible variety $Y$, and a homomorphism $\xi=\beta(\mu,\  ): Pic(Y) \to \Rr$. 
 Let $H$ be an ample divisor on $Y$.  Let $Pic^0(Y)$ be the group of divisors algebraically equivalent to $0$.  
 It is known (\cite{kleiman}, Proposition 5.3 or Theorem 5.4, and Proposition 5.10) that $Pic^0(Y)$ is an algebraic variety, in particular constructible.  For any $a \in Pic^0(Y)$, (as for any element of $Pic(Y)$),
 for large enough $m$, the line bundle corresponding to $mH - a$ has a nonzero section.  So $-a+mH \geq 0$.
 Since $Pic^0(Y)$ is constructible, by compactness, there exists $m$ such that for all $a \in Pic^0(Y)$,
 $-a+mH \geq 0$.  So $\xi(a) \leq m \xi(H)$.  Thus $\xi | Pic^0(Y)$ is a homomorphism into $\Rr$ with bounded image;
 so it must be zero.  Since $\Rr$ is torsion-free, if $mb \in Pic^0(Y)$ then $\xi(b)=0$.
 Hence $\xi$ vanishes on the divisors that have integral multiples algebraically equivalent to zero. 
   But by 
  \cite{kleiman}, Theorem 6.3, 
    these are the same as the divisors numerically equivalent to $0$.   
 So $\xi$ induces a homomorphism $NS(Y)_\Rr \to \Rr$.   and factors through that homomorphism.\eprf

Note that the homomorphism $\beta_\mu: Pic_\Rr(Y) \to \Rr$ obtained above is order-preserving with respect
to the effective cone;  i.e. an effective divisor has non-negative image.  This is because by definition,
$\beta(v,D)=v(g) \geq 0$ where $D$ is represented by $(g)$ on some affine open set $Y'$ with $\Oo_{Y'} \subset \Oo_v$,
and as $D$ is effective, $g \in \Oo_{Y'}$ so $v(g)\geq0$.  Taking into account that $\mu$ is a non-negative measure,
  the induced homomorphism $NS(Y)_\Rr \to \Rr$ is order-preserving with respect to the 
to the pseudo-effective cone.  

Let $S_K$ be the space of GVF
structures on $K$ that are trivial on $k$; it can be viewed as the space 
of quantifier-free GVF types whose restriction to ACF is the generic type of $X$.  If we pick any transcendence
basis $f_1,\ldots,f_n$ for $K$ over $k$, and $r \in \Rr$, then the subspace of types with $ht(f_i) \leq r$ is compact.
Each element $p$ of $S_K$ is induced by a measure $\mu \in M_K$ satisfying the product formula:
this can be proved using a homogeneous variant of the Riesz representation theorem, but will also follow
from the proof of \thmref{qft} below.  

It is clear that $\beta(\mu,D)$ depends only on $p$ and not on the choice of $\mu$; for instance
this may be checked for very ample $D$, with a section $f$, so $\beta(\mu,D)=\int \min(v(f),0) d\mu(v)$.
We thus write $\beta(p,D)$ for $\beta(\mu,D)$.

Let $N^1(K)$ be the direct limit of the $\Rr$-spaces $N^1(Y)$ along $\Ii_X$; let $N_1(K)$ be the dual space, the inverse limit
of the spaces $N_1(Y)$.  (Cf. \cite{bfj}.) Note that when $j: Y' \to Y$ is birational, $j^* N_1^+(Y) \subset N_1^+(Y')$.  
Let $N_1^+(K)$ be the inverse limit of the cones $N_1(Y)^+$, dual to the effective cones
of $N^1(Y)$.

We define a map $\alpha: S_K \to N_1^+(K)$ via the pairing $\beta: S_K \times N^1(K)\to \Rr  $.
So $\alpha = \lim \alpha_Y$ where $Y$ is similarly defined via $\beta:  S_Y \times N^1(Y) \to \Rr$,
with $Y $ a birational extension of $X$.

\begin{thm}\label{qft}  $\alpha: S_K \to N_1^+(K)$  is a homeomorphism.  \end{thm}

\prf Injectivity:  Let $p,p'$ be distinct elements of $S_K$.  They are represented by measures $\mu,\mu'$
that are distinct even up to renormalization.   By   \lemref{11}, for some $D$ we have $\beta(\mu,D) \neq \beta(\mu',D)$.
Thus $\alpha(p) \neq \alpha(p')$.

Continuity:   It suffices to prove that $\alpha: S_K \to N_1(Y)$ is continuous.   Now the topology on $N_1(Y)$
is generated by $O_{\beta,\gamma,D} = \{a: \beta < |a \cdot D | < \gamma\}$, where $\beta,\gamma$ are reals and $D$ is a very ample
divisor on $Y$.  Let $s_1,\ldots,s_n$ be global sections of $D$.  Then 
$\alpha \inv (O_{\beta,\gamma,D} ) $ is cut out by $\beta < ht(s_1,\ldots,s_n) < \gamma$.

A closed map:  

Let $f_1,\ldots,f_n$ be a transcendence basis for $K/k$; let $A$ be an ample divisor, such that each $f_i$ 
has poles at most on $A$.   Let $D_i$ be the divisor of poles of $f_i$. 
Let $C$ be a closed subset of $S_X$.  For $t >0$ let 
 $R_t =\{b \in N_1^+(K): (b,D_i) \leq t\}$.  Then $R_t$ is contained in the interior of $R_{t'}$ for $t<t'$, and
 $\union_t R_t = N_1^+(K)$.  Thus to show that a set is closed, it suffices 
to show that the intersection with each $R_t$ is closed.  As $\alpha \inv(R_t)$ is closed, to show that $\alpha$
is a closed map, it suffices to show that 
$\alpha(C)$
is closed when $C$ is a closed subset of  $\alpha \inv R_t  $ for some $t$.  But such a set is compact as $ht(f_i)$
is bounded on it.

In particular, the image of $\a$ is closed.   Hence to show that $\a$ is surjective, it 
suffices to prove that $\a$ is dense.  For this, it suffices to show that each $\a_Y$ is surjective.  
Using de Jong, find a smooth projective variety $\tY$ and a generically finite morphism $\phi: \tY \to Y$.
By \lemref{functorialcone}, $N_1^+(\tY) \to N_1^+(Y)$ is surjective.   This reduces the statement to the smooth projective case; here it  follows from \lemref{surj3}.   
\eprf

\begin{remark} \label{heightsuffices} The fact that continuity was proved using the open sets \[\beta < ht(f_1,\ldots,f_n) < \gamma\]
alone shows that the projective height functions, along with the field operations, generates the entire language of GVFs.  \end{remark}

\ssec{Discrete globalizing measures}
\label{dgm}
Let $X$ be a smooth projective variety.   
Let $\wW$
be the set of irreducible hypersurfaces of $X$.

Let $a$ be a class in $N_1(X)^+$.

 We  let  $\mu_a$ be the measure, supported on divisorial valuations $v_D$ of $K$
for each irreducible hypersurface $D$ of $X$, and with $\mu_a(v_D) = a \cdot [D]$.

The same definition extends additively to Cartier divisors, in such a way that
principle divisors map to $0$.  This amounts to the fact that a rational function on a curve
has as many zeroes as poles.

From this  it is clear that for $f \in K$, letting $(f)=\sum m_i D_i$ be the divisor of zeroes and poles of $f$,
we have $\sum m_i \mu(D_i) = 0$.   Indeed let $p \in C_t \meet D_i$; if $D_i$ is cut out by $(g)$ locally at $p$ 
then $(g^{m_i})=(f)$ on this local ring, so $m_i i(C_t,D_i) $ is the order of zero (or pole) of $f$ at $p$ on $C_t$,
and the sum of these is zero.   Thus $\mu_a$ is a GVF measure.  In particular, it defines an element of $S_X$,
depending only on the numerical equivalence class of $C$, that we denote $\delta_X([C])$.  This gives a map

\[ \delta_X: N_1^+(X) \to S_K \]

\begin{lem} \label{surj3} $\alpha_X \circ \delta_X = Id_{N_1^+(X)}$.  \end{lem}
\prf This amounts to the relations:  $v_D(D)=1$, $v_D(D')=0$ when $D,D' \in \wW$ are distinct, 
and $v_D$ is the valuation corresponding to $D$, evaluated on $D,D'$ as Cartier divisors.  Both are straightforward.
\eprf

 \begin{rem}
It is  not true that $ \delta_X \circ \alpha_X $ is the identity on $S_K$; but   $ \delta_X \circ \alpha_X (p) $
and $p$ do give the same   value to a formula $\int t(vf_1,\ldots,vf_m) dv$ if $f_1,\ldots,f_m$ extend to
morphisms $X \to \Pp^1$,   $t$ is a term formed out of $+,-, \min$, and for any  minimum occuring in 
  $t(vf_1,\ldots,vf_m)$, the corresponding minimum of Cartier divisors $(f_1),\ldots,(f_m)$ is 
is a stable meet on $X$. 

 In particular, the same height is computed for rational functions $f$ extending
to a morphism $X \to \Pp^1$.

This makes it possible to define $\alpha_K$ as the limit of $\alpha_Y$ over all blowups $Y$, and prove injectivity
in a  more explicit way.
\end{rem}

\begin{rem}  \label{boundary} Let $a \in N_1^+(X)$;  assume $a \cdot b =0$ for some  effective $b$
with nonzero section growth.  
   Then there exists a nonconstant rational function $g$ on $X$
such that   $F(g)$ is a constant subfield of the GVF $F(X)_a$.    The converse is also true.  \end{rem}

\prf  By assumption there exists a nonconstant rational function $g$
with poles only on  (a representative of) $b$.  So the height of $g$ in $F(X)_a$ is zero.  \eprf

\ssec{Quantifier-free types over the function field of a curve}

Here the picture is very similar.   Let $k$ be a constant field, and let $k_1=k(U)$ for a curve $U$.  Let $K=k_1(X)$ be a finitely generated field extension of $k_1$.  Let $\mu_1$ be a nontrivial
GVF structure on $k_1/k$.  
Then  $N_1^+(k_1) = \Rr^{ \geq 0}$, with $\mu_1$ corresponding to $1$;  and we have a natural map $\pi_*: N_1^+(K) \to N_1^+(k_1)$; it is obtained as the limit of the projections 
$\pi_*:  N_1(X) \to N_1(U)$ at the level of varieties.   Let $N_1^+(K/\mu_1) $ be the pullback of  $1$.  We   have an 
$\Rr^{>0}$-action on $N_1(X) $ and on $N_1(U)$ respecting $\pi_*| N_1(X)$; so 
$N_1^+(K/\mu_1)$ can also be identified with the set of elements of $N_1^+(K)$ with nontrivial images in $N_1^+(U)$, 
 modulo the $\Rr^{>0}$-action.  

Letting $S_X$ denote the space of quantifier-free types over $(k_1,\mu_1)$, we immediately obtain from \thmref{qft},
by restriction:

\begin{thm}\label{qft1}  $\alpha: S_X \to N_1^+(K/\mu_1)$  is a homeomorphism.  \end{thm}

\end{section}
\begin{section}{Existential closedness of   $k(x)^{alg}[r] $}

\label{xc-s}


Let $k=k^{alg}$.  

 Let $\mu$ be a globalizing measure for $k(V)$ over $k$.  Let $D$ be an ample divisor on $V$.

\begin{lem}[Artin-Whaples] \label{whaples}  Let $C$ be a curve of genus $g$ over $k$.    Let $x \in k(C) \m k$, and let $r>0$.    Then $k(C)$ admits a unique GVF structure over $k$ with $ht(x)=r$.  
\end{lem}

\prf  The nontrivial valuations of $k(C)$ over $k$ can be identified with the points of $C(k)$.   They form a discrete set; the compactification includes also the trivial valuation,
which plays no role here.     We have to show that if $\mu$ is a measure on $C(k)$ satisfying the product formula,
then $\mu$ gives equal weight to each element of $k(C)$.  
Let $a,b$ be two   points.   Then $(n+g)a - nb $ is effective, so $(n+g) \mu(a) \geq n \mu(b)$.  Thus $(1+g/n) \mu(a) \geq \mu(b)$.  Letting $n \to \infty$, $\mu(a)\geq \mu(b)$.
So $\mu(a)=\mu(b)$ for all $a,b \in k(C)$.  \eprf

A $c$-twisting automorphism $\si$ of $K$ is an isomorphism from $K$ to $K'$, where $K'$ is the renormalization of $K$ moving from $\mu$ to $c \mu$.  

 \begin{cor}  \label{homogeneity1}  There is a unique $GVF$ structure $k(x)^{alg}[r]$ on $k(x)^{alg}$ over $k$, with $ht(x)=r$.  
   Moreover if 
 $F$ is a finitely generated extension field of $k$ of transcendence degree $1$, and $g,g': F \to k(x)^{alg}$ are two field embeddings, 
 and $a \in F \m k$, there exists an $ht(f'(a))/ ht(f(a))$-twisting automorphism $\si$ of $k(x)^{alg}$ with $g'=\si \circ g$. 
 \end{cor}
 
 \begin{rem} \label{rationalrepa}   If $r/r' \in \Qq$ then $k(x)^{alg}[r] \cong k(x)^{alg}[r']$.   \rm
 Let $m \in \Nn$ and let $\si$ be an automorphism of the field $k(x)^{alg}$ with $\si(x)=x^m$.
 Then $\si$ is $m$-twisting so it gives an isomorphism $k(x)^{alg}[r] \to k(x)^{alg}[mr]$.  Similarly find
 an $m'$-twisting $\si'$, with $m/m'=r/r'$.  Then $(\si') \inv (\si)$ is an isomorphism $k(x)^{alg}[r] \to k(x)^{alg}[r']$. \end{rem}
 
 \begin{cor} \label{homogeneity} Let $K^*$ be an ultrapower of $K=k(x)^{alg}$.  Let $F$ be finitely generated over $k$,  of transcendence degree $1$,  and let $f,f': F \to K^*$ be two embeddings, with $ht(f(a))=ht(f'(a))>0$ for some $a \in F $.
 Then   there exists an automorphism $\si^*$ of $K^*$ with $\si \circ f  = f'$.   
 
  This 
 also holds if $K^*$ is an ultraproduct of GVF fields $k_i(x)^{alg}$ of increasing positive characteristics, and $F$
 is finitely generated over $\Qq$, of transcendence degree $1$.    \end{cor}
 \prf  
 $F$ is the field of fractions of a  finitely generated $k$-algebra $D$, with $a \in D$.  
 $f|D,f'|D$ can be represented as ultrapowers of homomorphisms $f_i,f_i': F' \to K$.  We have 
 $ht(f_i(a)) \to ht(f(a))$ and  $ht(f_i'(a)) \to ht(f'(a))$ along the ultrafilter $u$.   So $r_i=ht(f_i(a))/ht(f_i'(a))$ approaches $1$ along $u$.  Note that $f_i,f_i'$ extend to field embeddings (first on $k(a)$ since $f_i(a), f_i'(a) \notin k = k^{alg}$;
 then on $F$ since it is  a finite extension of $k(a)$).
 Using   \corref{homogeneity1}  let $\si_i$ be an $r_i$-twisting automorphism of $K$ with $f_i'=\si \circ f_i $.  
 Let $\si$ be the ultraproduct of the $\si_i$.  Then $f'=\si \circ f$ and $\si$ is a $1$-twisting automorphism, i.e.
 simply an automorphism.    The  statement with ultraproducts is proved similarly.
 \eprf 
 
 \begin{proof}[Proof of \thmref{ec}]
 
 We first show that any quantifier-free type $q$ over $k$ can be approximately realized in $K =k(t)^{alg}[1]$.  
 Consider first a smooth projective variety $X$ over $k$, and a family of curves $C_t$ on $X$ as in \secref{dgm},
 with class $[C] \in N_1^+(X)$.  
 We have the corresponding quantifier-free type $\delta_X(C)$.  Now the inclusion of $C_t$ in $X$ corresponds to 
 an element $a_t \in X(k(C_t)) $.   It is easy to see that any approximation to $\delta_X(C)$ involving finitely many
 divisors is realized by $a_t$, for almost all $t$ (namely as soon as $C_t$ is not contained in  the support of these finitely many  divisors).    Thus $\delta_X([C])$ is  approximately realized in $K$.   Since $K \cong K[m]$ for integral $m$
 by \remref{rationalrepa}, we aso have $\delta_X(\frac{1}{m}[C])$  approximately realized in $K$.  
 Now by    \corref{eclem},  \thmref{qft} and \lemref{surj}, the types $\delta_Y(\frac{1}{m}[C])$ for $Y$ a blow-up of $X$, and 
 $C$ as above on $Y$, are dense on $S_X$.   Thus all  types in $S_X$ are approximately realized in $K$.

 Now we must also consider types over $K$ and not just over $k$.  
   Let $K \leq L$, $L$ finitely generated over $K$, with a GVF structure.   We have to find a $K$- embedding of GVF's  $L \to K^*$, where $K^*$ is an ultrapower   of $K$.    By the above we do have a $k$-embedding   of GVFs $h: L \to K^*$.   
 Let $g$ be the inclusion of $K$ in $K^*$,  and $g' = h|K$.  Using \lemref{homogeneity}, find  $\si \in Aut(K^*)$ with $\si \circ g' = g$.   Then $\si \circ h$ is an embedding of $L$ to $K^*$ over $K$, as required.  

 \end{proof}

\end{section}

\begin{section}{A non-archimedean Yau theorem}
 
  \thmref{surj} can be viewed as an absolute version of  the non-archimedean Yau  theorem of \cite{bfj2} and   
  \cite{kt}.      To bring this out, we give a   proof along the same lines of a   version of the relative theorem. 
  (Our proof is valid in any characteristic, and addresses
  to some extent the `regularity' issue raised in \cite{bfj2} in the model case, in that we show that the metric
  is determined by a big divisor on the same model.)  
  
  We are given not only a smooth projective $X$, but also a morphism to $U$.   In the relative (Yau) version, the data is a measure
    on the vertical divisors,  as well as an ample line bundle $L$ on the generic fiber of $X$; the theorem yields 
 a big line bundle on $X$ compatible   with the measure (vertically)    and with $L$ (generically.)  By comparison     
     in \thmref{surj},  the data is a measure  on {\em all} divisors of $X$, including the horizontal ones and the outcome is again big divisor compatible with both vertical and horizontal measures.

\ssec{The relative setting}  Assume now that $\dim(X)= n+1$, and we are given not only $X$ but
also a morphism $\pi: X \to U$ and a point
$0 \in U$, where $U$ is a curve.  Let $K=k(U)$, and let $X_K$ be the generic fiber, and $X_0=f \inv(0)$ the special fiber.
Let $I$ be the set of irreducible components of $X_0$.  
 Within $N^1(X)$, we have the subspace  $\vdiv$ of {\em vertical divisors} generated by $I$.  
   Note that $[X_0]^2=0$ and moreover $[X_0] \cdot [E]=0$ for any  $E \in I$, since $X_0$ can be moved to any other fiber.
   By 
  \lemref{pdc} (applied to $\vdiv$ and the vector $f=[X_0]$, and taking into account the Grothendieck connectedness theorem), $\vdiv / (f)$ is negative definite and $f=0$
  is the unique relation on the generators $I$;     as $f \neq 0$, in fact $I$ is a basis for $\vdiv$.  
  
  Any ample $D \in N^1(X)$ defines, on the one hand, an ample line bundle $L$ on $X_K$, obtained by restriction of $D$ to 
  the generic fiber.  On the other hand, the point $0 \in U$ determines a valuation  $v_0$  on $K$, and by extending
  rational sections of $L$ to sections of $D$ we obtain a metrization of $L$ over the valued field $(K,v_0)$. 
  (Actually for any big $D \in N^1(X)$ with fixed ample generic part $L$, we can define a metrization of $L$
  over $(K,v_0)$, the limit of the metrics associated to ample divisors $D-E$ with $E$ effective and vertical
  and $D-E$ ample, over all base change morphisms $X' \to X$ that are isomorphisms away from $X_0$.)
  
 In Arakelov theory, this is viewed as analogous to a smooth  metrization of $L$, assuming 
 given an embedding of $K$ into $\Cc$ in place of the valuation $v_0$.  
 
 In the K\"ahler case, a metrized line bundle gives rise to a measure on $X$,  obtained by taking the $n$'th power
 of the curvature $2$-form.   A non-archimedean analogue of this construction was defined in \cite{acl}, following work of Gillet-Soul\'e, Gubler, Zhang  and others.
   In the basic geometric case we are considering, the measure is very simple to define.  It concentrates on $I$;
  any  $c \in I$ is given weight $m_c c \cdot D^n$, where $m_c$ is the multiplicity with which $c$ occurs in $X_0$.
  
  More generally, if $D$ is a divisor on a blowup $b: X' \to X$ we can still define the measure
  by $\mu_D(c) = m_c c \cdot b_*(D^n) = m_c b^*(c) D^n$.  
  
       For us it is natural to extend the definition to any big Cartier $\Rr$-  divisor $D$
  (with generic part $L$):  
  \[ \hat{\mu}_D(c) = m_c c \cdot \psi(D) \]

\begin{lem}  \label{bigtobig}
Any big line bundle $L$ of $X_K$ extends to a big divisor $\bL$ of $X$.\end{lem}
\prf   We can  represent $L$ by an effective divisor $D_K$ on $X_K$.  
There exist $g_1,\ldots,g_n \in k(X)=K(X_K)$
that are algebraically independent over $K$, and whose poles are at most $D_K$.  We can view the $g_i$
as elements of $k(X)$, and by adding  some vertical divisors,  can find a divisor $D$ on $X$ with generic part $D_K$,
such that the poles of the $g_i$
are bounded by $D$.  On the other hand we can find $f_1,\ldots,f_r \in K = k(U)$, $r=tr.deg._k K = \dim(U)$,
algebraically independent over $k$;   
  by adding to $D$ the pullback
of the polar divisors of $f_1,\ldots,f_r$ to $D$ we may assume each $f_i   \in L(D)$
as well.    Finally since $D$ is effective,  we have $1 \in L(D)$.  It follows that $D$ is big.  \eprf

\begin{thm}\label{nacy}  Let $L$ be a big line bundle on $X$, represented by an irreducible Cartier divisor.
  Let $I_t$ be the set of irreducible components
of $X_t$, $I = \union_t I_t$.  
 Let $\mu$ be a 
 nowhere vanishing positive measure on the vertical divisors, such that the total mass of each fiber $X_t$ is
  $\mu(I_t)= D^n = \deg(D)$.    Then there exists a big $\Rr$-Cartier divisor $B$ on $X$ 
  with generic part $L$,
  $\hat{\mu}_B = \mu$. 
 
\end{thm}

Note that for   all but finitely many $t \in U$, $I_t$ has a single point, so $\mu$ is fully constrained on $I_t$. 

\begin{rem} 
 If only given a measure on $I_0$ of total mass $L^n$, it can be extended to a nowhere vanishing positive measure on $I$, as the proof below shows; thus  \thmref{nacy} includes the local version usually considered. \end{rem}

\prf   Let $D_K$ be an irreducible Cartier divisor on $X_K$, representing $L$.

Let $\gamma: Pic_\Rr(X) \to Pic_\Rr(X_K)$ be `restriction to the generic fiber'.   
We will also write $\gamma$ for restriction of $\Rr$-Cartier divisors.  Fix an $\Rr$-Cartier divisor $D$
with generic part $D_K$.
 
Let $E \leq Pic_\Rr(X)$ the the    vector  space  generated by  the class of $D$ and by the 
vertical divisors.   
Let $E^+ $ be the cone in $E$ generated by effective divisors.  
   
 Consider  the codimension one subspace $E_0:=\ker(\gamma)$ of $E$. 
 Let $E_0^+$ be the   cone of effective Cartier divisors in $E_0$.  $E$ is generated by the finite set $I$, 
 using \remref{hodge2c};   so it is finite dimensional, and closed; the natural homomorphism
 $Pic_\Rr(X) \to N^1(X)$ is, again by  \remref{hodge2c},  injective on $E_0$, and $E_0^+$ is the pulback of
 the effective cone in $N^1(X)$.  
 
  Let $U$ be the intersection of $E$ with the pullback to $Pic(X)$ of the cone of big divisors.  By \lemref{bigtobig}, this 
 intersection is nonempty, so it is an open cone on $E$.  The hypotheses (1)-(3) of \thmref{convexity}
 are satisfied for the function $\vol^{\frac{1}{n+1}} |E$, since this is shown for $N^1(X)$ in \thmref{surj}, and they remain true upon restriction to the image of $E$.   For $B \in U$, write $\psi(B)$ for $\psi$ of the image of $B$ in $N^1(X)$.
   
 Now $\mu$ defines a linear map on $E_0$, positive on $E_0^+$.
   By the strict positivity assumption,   $\mu$  lies in the interior of the dual of $E_0^+$
 by \lemref{interiorofdual}, by  \lemref{coneproject} it lifts to an element $\mu'$ in the interior of the dual
 of $E^+$.   
 
By \corref{convexity-c},  there exists $B'  \in E \meet U$ such that
$\psi(B') \cdot c = \mu'(c)$ for any $c \in I$.  We may represent $B'$ by  an $\Rr$-Cartier divisor
    $B'=\alpha' D + F'$.
    
It remains to show that  $\alpha'$ may be chosen to be $1$.   We do not claim that it necessarily equals
$1$, but rather that there exists an $\Rr$-Cartier $B \leq B'$ with the same sections, with $\psi(B)=\psi(B')$,
so that $B = D+F$.   Indeed take $B$, as an $\Rr$-Weil divisor, be $\sup \{ (f)/n:  (f) \leq  nB, n =1 ,2,\cdots\}$.   Then
by definition $B,B'$ have the same sections, hence same Okounkov bodies,   the same volume, and
by \propref{sectionsonly}, $\psi(B)=\psi(B')$.  
  
Restricting
to the generic fiber, and viewing $B,B'$ as $\Rr$-Weil divisors,  we have
$\gamma(B) \leq \gamma(B')= \alpha' D_K$ so using the irreducibility of $D_K$,
 $\gamma(B')= \alpha D_K$ for some $\alpha \leq \alpha'$.    

By definition, $\psi(B)$ is the limit of $H^n$ over nef divisors $H \leq b^*(B)$, where
$b:X' \to X$ is a blowup of $X$.  Working in $X'$, let $D' = b  ^* D, F'=m^*F$, $\gamma':N_1(X') \to N_1(X'_K)$
the restriction.   Let $f$ denote the class of a generic fiber of   $\pi$.  Note $f \cdot F = 0$ for any vertical
divisor $F$.  

$ H^n \cdot b^*(f) \leq (\alpha D+ F)^n \cdot f = \alpha^n L^n = \alpha^n \mu(f)$.  Thus taking the limit,
$\mu(f) \leq \alpha^n \mu(f)$ so $\alpha \geq 1$.

On the other hand  if $\alpha >1$ then some section witnesses this;   on some blowup $b: X' \to X$ there exists a semi-free (hence nef)
$\Qq$- divisor $H \leq \alpha D + F$ with generic part $H_K \geq \beta D$, for some $\beta >1$.   But
  in this case for some vertical $F'$ we have $H \geq  \beta D + F'$, so 
  $\mu(f)=\psi(B) \cdot f \geq H^n  \cdot b^*(f) \geq (\beta D +F')^n \cdot f = \beta^n D^n \cdot f =
  \beta^n 
  \mu(f)$, a contradiction.  Thus $\alpha=1$.   So $\hat{\mu}_B = \mu$.
  \eprf  

\begin{question} (and remarks). Clarify Calabi uniqueness in this setting. 
For semi-positive metrics on an ample   line bundle, supremum of the nefs on blowups dominated by $B$, is proved in \cite{yuan-zhang}.    In the global case, uniquenesss follows from
a the strict log-concavity of volume in non-radial directions, \cite{bfj}, \cite{cutkosky}.     

As far as I could see nothing is known outside the nef cone, and it is possible that the
collapse from $B$ to $B'$ above is responsible.  

Consider $B$ with generic part $D$,   $\hat{\mu}_B=\mu$.

\begin{itemize}
\item 
 The freedom in the choice of $\mu'$ corresponds to additive translation of $B'$ by a scalar multiple of $f$.
 At least in the multilinear domain, $B'+ \alpha f$ represents the same measure on $I$.  
  Thus we can only expect to determine the image of $B$ 
in $N^1(X)/N^1(U)$.  This is analogous to the fact that in the usual Calabi theorem, uniqueness is up to a scalar multiple.  
\item In the non-multi-linear region of $\psi$, is addition of $f$ still the same as addition of a constant
to the Weil function, after collapse?  (Seems so.)
\item  An interesting way to prove   uniqueness of the asymptotic section spaces:  use the strict convexity in Brunn-Minkowski.  What can be concluded about $B,B'$ if the Okounkov body is the same for {\em every} choice of $v$?
\item Is this image determined by the associated metrization ($v \mapsto \lim_n 1/n \min v(sections)$).
\item  Is it the case that if  
 $\hat{\mu}_B= \hat{\mu}_{B'}$ and $B|D  = B'|D$ is ample then the nef divisors on blowups below $B$ are cofinal
 with those below $B'$ ?   (issue of f).
 \item   Assuming this cofinality, we have $B=B'$.  For every section of $B$ is a section of some free $N$
 on a blowup, $N \leq \pi^*B$; hence of $\pi^*B'$ hence of  $B'$.  (issue of f.)
 \item  Assuming this and that the limit is uniform, so that Yuan-Zhang Calabi applies, we do have $B= B'$.(issue of f.)
 \item  Does a big divisor defines the same metrization as a limit over blow-up of nef divisors approaching it?
 This apparently does not imply that it is nef itself.  Fujita says yes in terms of volume.  Pointwise, for 
 maximal rank $v$, the answer is also yes: Okounkov bodies must approach if volume does.  
\end{itemize}\end{question}
 
  \begin{question} What about higher dimensional $U$?  Analogue of splitting divisors is clear. 
 \end{question}

 \end{section}

  \end{document}